\documentclass{article}
\usepackage{graphicx}
\usepackage{amsmath}
\usepackage{amssymb}
\usepackage{amsfonts}
\usepackage{amsthm}
\usepackage{comment}
\usepackage[backend=bibtex,maxbibnames=99,style=alphabetic]{biblatex}
\usepackage{xcolor}
\usepackage{tikz-cd}
\usepackage[margin=1in]{geometry}
\usepackage{mathrsfs}
\usepackage[title]{appendix}
\usepackage{mathtools}
\mathtoolsset{showonlyrefs}
\usepackage{quiver}
\usepackage{todonotes}
\usepackage[colorlinks=true,linkcolor=blue,citecolor=blue,urlcolor=blue]{hyperref}

\usepackage[normalem]{ulem}

\usepackage{enumitem}

\newtheorem{theorem}{Theorem}
\newtheorem{proposition}{Proposition}[section]
\newtheorem{lemma}[proposition]{Lemma}
\newtheorem{corollary}[proposition]{Corollary}
\theoremstyle{definition}
\newtheorem{definition}[proposition]{Definition}
\newtheorem{example}[proposition]{Example}
\newtheorem{remark}[proposition]{Remark}
\newtheorem{convention}[proposition]{Convention}

\newcommand{\R}{\mathbb{R}}

\newcommand{\M}{\mathcal{M}}
\newcommand{\MZpNoSim}{\widehat{\M}_p(Z)}
\newcommand{\MZp}{\M_p(Z)}
\newcommand{\MZ}[1]{\M_{#1}(Z)}
\newcommand{\dgwz}{\mathrm{GW}_p^Z}
\newcommand{\dgw}[1]{\mathrm{GW}_{#1}^Z}

\newcommand{\WZ}{\mathcal{P}_p(Z)}
\newcommand{\WZp}[1]{\mathcal{P}_{#1}(Z)}
\newcommand{\wz}{\mathrm{W}_p^Z}
\newcommand{\wzp}[1]{\mathrm{W}_{#1}^Z}

\newcommand{\Lp}{\mathcal{L}_p}
\newcommand{\Lpp}[1]{\mathcal{L}_{#1}}
\newcommand{\dpz}{\mathrm{D}^Z_p}
\newcommand{\dz}[1]{\mathrm{D}^Z_{#1}}

\newcommand{\curv}{\underline{\operatorname{curv}}}
\newcommand{\curvupper}{\overline{\operatorname{curv}}}

\newcommand{\angk}[2]{{\sphericalangle_K}\bigl(#1;#2\bigr)}
\newcommand{\disp}{\mathrm{dis}_p}

\title{Metric Geometry of Lebesgue, Wasserstein, and Gromov-Wasserstein Spaces: Submetries, Curvature, and Geodesics}
\author{Martin Bauer, Facundo M\'{e}moli, Tom Needham, Mao Nishino}

\hypersetup{
  pdftitle={Metric Geometry of Lebesgue, Wasserstein, and Gromov-Wasserstein Spaces:
            Submetries, Curvature, and Geodesics},
  pdfauthor={Martin Bauer, Facundo Mémoli, Tom Needham, and Mao Nishino}
}

\begin{document}

\maketitle

\begin{abstract}
A metric space $Z$ gives rise to three natural classes of infinite-dimensional metric spaces associated to $Z$: $p$-Wasserstein spaces of probability measures on $Z$, nonlinear Lebesgue $L^p$-spaces of $Z$-valued maps, and $p$-Gromov-Wasserstein spaces of $Z$-valued kernels. The latter class, referred to as $Z$-Gromov-Wasserstein ($Z$-GW) spaces, extends the classical Gromov-Wasserstein framework from metric measure spaces to more general, possibly attributed, network-like structures, and unifies many GW-type distances that nowadays play a significant role in metric geometry, data science and machine learning. In this article we develop a unified metric-geometric theory of these three classes of spaces, with a particular focus on the $Z$-GW spaces. Our first main result identifies a fundamental submetry structure linking them: the nonlinear Lebesgue space maps via a submetry onto the $Z$-GW space, which in turn maps via a submetry onto the Wasserstein space. This structure provides a mechanism for transferring geometric information among the three spaces. We apply this framework to geodesics and Alexandrov curvature. For $1<p<\infty$, we prove that geodesicity of $Z$ is equivalent to geodesicity of each of the three associated spaces; in the endpoint case $p=1$, all three associated spaces are geodesic, even when $Z$ is not. We also characterize geodesics in the $Z$-GW space as generalized interpolations, extending a known characterization in the classical setting due to Sturm. Finally, we give a complete classification of Alexandrov curvature bounds for these spaces in terms of the curvature of $Z$. Thus, while the main focus of the paper is a new metric-geometric theory of $Z$-GW spaces, the submetry framework also extends classical theorems for Wasserstein and Gromov-Wasserstein spaces and yields new geometric consequences for nonlinear Lebesgue spaces.
\end{abstract}

\tableofcontents

\section{Introduction}
A metric space $Z$ gives rise to three natural classes of infinite-dimensional metric spaces associated to $Z$: Wasserstein spaces of probability measures on $Z$, nonlinear Lebesgue spaces of $Z$-valued maps, and Gromov-Wasserstein-type spaces of $Z$-valued kernels. These spaces appear in optimal transport, metric geometry, and geometric data analysis, and their geometry is expected to reflect, in different ways, the geometry of the underlying space $Z$. A central theme of this paper is that these three classes should not be studied in isolation: they are linked by a natural submetry structure, and this structure allows geometric properties to be transferred among them.

Our main focus is the $Z$-Gromov-Wasserstein ($Z$-GW) space, the least understood member of this triad. The classical Gromov-Wasserstein distances form a family of metrics on the space of isomorphism classes of metric measure spaces. They were introduced and studied in the foundational papers~\cite{Memoli_2007,Memoli2011GromovWasserstein} and have found numerous applications in recent years to, for example, geometry processing~\cite{peyre2016gromov,solomon2016entropic,chowdhury2021quantized}, network analysis~\cite{chowdhury2019gromovwassersteindistancenetworksstable,xu2019scalable,Vayer2020FGW,chowdhury2021generalized}, computational biology~\cite{demetci2022scot,demetci2022scotv2} and data visualization~\cite{vandistributional, clark2025generalized}. Several of these applications call for generalized notions of GW distances which are designed to compare more complicated data objects, such as networks endowed with node or edge attributes~\cite{Vayer2020FGW,kawano2024multi}. With a view toward establishing a general theory of GW-type distances, our recent article~\cite{bauer2025zgromovwassersteindistance} introduced a family of GW distances which encompasses many of the variants in the literature. This family is parameterized by a choice of separable metric space $Z$ and number $p \in [1,\infty]$;  these distances are therefore referred to as \emph{$Z$-Gromov-Wasserstein $p$-distances} (or just \emph{$Z$-GW distances}). 

In the present paper, we study the metric geometry of the resulting $Z$-GW spaces and their relation to the corresponding Wasserstein and nonlinear Lebesgue spaces. We address three fundamental questions:
\begin{enumerate}
\item Under which conditions on $Z$ and $p$ is the $Z$-GW $p$-distance a geodesic metric?
\item How do Alexandrov curvature bounds on $Z$ translate to curvature bounds for the $Z$-GW $p$-distance?
\item How do the above properties relate to those of the other spaces associated to $Z$, the Wasserstein and nonlinear Lebesgue spaces?
\end{enumerate}
Our main contributions are full resolutions of these questions for arbitrary $p < \infty$ and Polish, geodesic $Z$. Theorems \ref{thm:geodesics_main} and \ref{thm:curvature} resolve the geodesic and curvature questions, respectively, while Theorem \ref{thm:submetries_for_GW} provides the submetry framework relating the three classes of $Z$-spaces. These results generalize seminal work of Sturm in the setting of classical GW distances~\cite{sturm2020}, which had already been adapted to several specific variants of GW distance~\cite{chowdhury2020gromov,Chowdhury2023HypergraphCOT,zhang2024geometry,zhang2025topological}  via extensions of Sturm's methods; in contrast, the level of abstraction considered in this paper requires several genuinely new techniques.  They also illuminate connections to many known results for the Wasserstein and nonlinear Lebesgue spaces (appearing in, e.g.,~\cite{korevaar1993sobolev,jost1994equilibrium,sturm2001nonlinear,serieys2025nonlinear,serieys2026nonlinear}), and yield new results on the Alexandrov geometry of these spaces.

The proofs in this paper are largely based on Theorem \ref{thm:submetries_for_GW}, which establishes a novel connection between the $Z$-GW space, the Wasserstein spaces over $Z$ and the nonlinear Lebesgue spaces of maps valued in $Z$, expressed in terms of submetries between the spaces. In order to state this connection precisely, and to describe our main results in more detail, we recall the definitions of these spaces below.  

\subsection{Metric Spaces Associated to \texorpdfstring{$Z$}{Z}}

\begin{figure}
    \centering
    \includegraphics[width=0.95\linewidth]{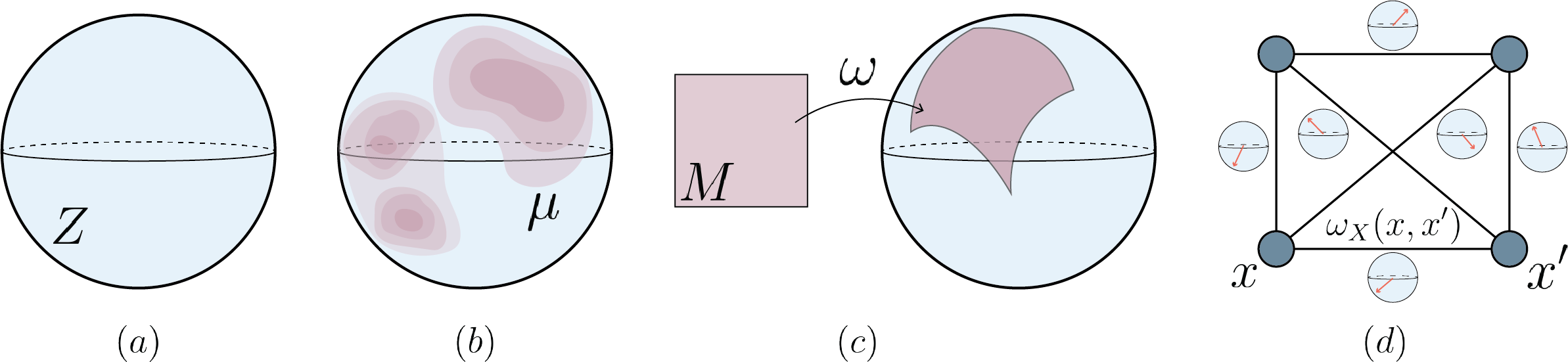}
    \caption{Examples of $Z$-Spaces. \textbf{(a)} A metric space $Z$, which is typically assumed to be geodesic and Polish (e.g., a 2-sphere with geodesic distance). \textbf{(b)} A point in the $Z$-Wasserstein space $\WZ$ is a distribution $\mu$ over $Z$, illustrated here by a heat map. \textbf{(c)} A point in the nonlinear Lebesgue space $\Lp(M,Z)$ is a ($p$-integrable) map $\omega$ from the probability space $M$ into $Z$. \textbf{(d)} A point in the $Z$-Gromov-Wasserstein space $\MZp$ is a triple $(X,\omega_X,\mu_X)$, where $(X,\mu_X)$ is a probability space and $\omega_X \in \Lp(X^2,Z)$. Here, the space $(X,\mu_X)$ is variable, so that there is no given correspondence between the underlying measure spaces of two points in $\MZp$. We visualize a point in $\MZp$ as a dense network, where $X$ is a node set and $\omega_X(x,x')$ assigns a point in $Z$ to the edge joining $x$ to $x'$; this is intended to evoke the intuition that comparing two points in $\MZp$ can be viewed as a generalized graph alignment problem.}
    \label{fig:ZSpaces}
\end{figure}

Let $(Z,d_Z)$ be a separable metric space. 
This paper is concerned with geometric properties of, and relations between, several natural  metric spaces associated to $Z$  which arise in optimal transport theory and various areas of applied mathematics. We generally refer to these as \emph{$Z$-spaces}; the $Z$-spaces under consideration fall into three classes, which we now briefly describe. Examples of these spaces are illustrated in Figure \ref{fig:ZSpaces}.

\smallskip

\noindent {\bf Wasserstein Spaces.} Likely the most familiar class of $Z$-spaces is the class of \emph{$Z$-Wasserstein Spaces} $(\WZ, \wz)$. For $p \in [1,\infty)$, let $\WZ$ denote the space of probability measures with finite $p$th moments. This is  endowed with the \emph{$p$}-Wasserstein distance, defined for $\mu,\mu' \in \WZ$ as 
\begin{equation}\label{eqn:wasserstein_distance}
\wz(\mu,\mu') \coloneqq \inf_\pi \left( \int_{Z \times Z} d_Z(z,z')^p d\pi(z,z')\right)^{1/p},
\end{equation}
where the infimum is over the space of \emph{couplings} of $\mu$ and $\mu'$, or measures $\pi$ on $Z \times Z$ whose marginals are $\mu$ and $\mu'$, respectively. The $Z$-Wasserstein spaces are central objects of study in optimal transport, whose geometry is quite well understood; see~\cite{villani2008optimal,villani2021topics} as canonical references.

\smallskip 

\noindent {\bf Lebesgue Spaces.} A relatively simpler, but comparatively less studied, class of $Z$-spaces is the class of \emph{$Z$-valued (nonlinear) Lebesgue spaces} $(\Lp(M,Z),\dpz)$. Let $(M,\mu)$ be a fixed Polish (separable and completely metrizable) probability space and let $p \in [1,\infty)$. Then $\Lp(M,Z)$ is the space of measurable maps $\omega:M \to Z$ such that $\int_M d_Z(\omega(m),z_0)^p d\mu(m) < \infty$ for some fixed $z_0 \in Z$. The distance between two such maps $\omega,\omega' \in \Lp(M,Z)$ is given by 
\begin{equation}\label{eqn:Lp_distance}
\dpz(\omega,\omega') \coloneqq \left(\int_M d_Z(\omega(x),\omega'(x))^p d\mu(x) \right)^{1/p},
\end{equation}
and maps are only considered up to the $\dpz=0$ equivalence relation. The study of nonlinear Lebesgue spaces goes back to at least~\cite{korevaar1993sobolev}, and these spaces are frequently considered as useful examples in Alexandrov geometry, particularly when $Z$ is assumed to be a Hadamard space~\cite{jost1994equilibrium,sturm2003probability,bacak2014convex}. The recent article~\cite{serieys2025nonlinear} includes a thorough survey of uses of nonlinear Lebesgue spaces in applied mathematics and medical imaging.

\smallskip

\noindent {\bf Gromov-Wasserstein Spaces.} Finally, we introduce the $Z$-Gromov-Wasserstein spaces, the main objects of study of this paper, by beginning with motivation from the more classical setting. Let $(X,d_X,\mu_X)$ and $(Y,d_Y,\mu_Y)$ be metric measure spaces (i.e., $(X,d_X)$ is a compact, Polish metric space and $\mu_X$ is a fully-supported Borel probability measure on $X$) and let $p \in [1,\infty)$. The \emph{Gromov-Wasserstein (GW) $p$-distance} between $X$ and $Y$ is given by 
\begin{equation}\label{eqn:GW_distance_intro}
\mathrm{GW}_p(X,Y) \coloneqq \frac{1}{2} \inf_\pi \left(\iint_{(X \times Y)^2} |d_X(x,x') - d_Y(y,y')|^p d\pi(x,y)d\pi(x',y')\right)^{1/p},
\end{equation}
where, as in the setting of the classical Wasserstein distance, the infimum is over couplings of $\mu_X$ and $\mu_Y$~\cite{Memoli_2007,Memoli2011GromovWasserstein}. The GW distances define a family of metrics on the space of metric measure spaces, which can be viewed as relaxations of the well-known Gromov-Hausdorff distance. The expression \eqref{eqn:GW_distance_intro} makes sense when the metrics $d_X$ and $d_Y$ are replaced with arbitrary $p$-integrable kernels on the respective spaces; the observation that this provides a meaningful comparison between arbitrary kernels goes back at least to \cite{peyre2016gromov}, and metric properties of this generalized GW distance were formalized in~\cite{chowdhury2019gromovwassersteindistancenetworksstable}---therein, triples of the form $(X,\omega_X,\mu_X)$, where $\omega_X:X \times X \to \R$ is an arbitrary $p$-integrable kernel, are referred to as \emph{measure networks}, and we use similar terminology for the general structures we consider below. The GW distance and its measure network relaxation have found numerous applications to problems in data science which require comparisons between unregistered objects, e.g.,~\cite{peyre2016gromov,xu2019scalable,chowdhury2020gromov,chowdhury2021generalized,demetci2022scot}.

In order to handle data with specialized structure, such as attributed graphs, the GW framework has been further extended in several ways over the years; see~\cite{Vayer2020FGW,Chowdhury2023HypergraphCOT,memoli2023ultrametric,kawano2024multi,yang2024exploiting}, and Section \ref{sec:variants_of_GW}, which describes some of these variants in detail. With the goal of developing a general framework for capturing these variants simultaneously, we introduced the $Z$-Gromov-Wasserstein distance in our previous paper~\cite{bauer2025zgromovwassersteindistance}. Therein, we define a \emph{$Z$-network} to be a structure of the form $(X,\omega_X,\mu_X)$, where $(X,\mu_X)$ is a Polish probability space and $\omega_X:X \times X \to Z$ is an element of $\Lp(X^2,Z)$. Intuitively, one can think of $X$ as a set of nodes in a dense network, with $\omega_X(x,x')$ assigning a $Z$-valued feature to the edge joining $x$ to $x'$; when comparing two structures of this form, there may be no known correspondence between node sets, which, from a metric perspective, distinguishes this situation from that of nonlinear Lebesgue spaces. A mild adaptation of the classical Gromov-Wasserstein distance \eqref{eqn:GW_distance_intro} defines a notion of distance between $Z$-networks $X$ and $Y$, which we refer to as the \emph{$Z$-Gromov-Wasserstein distance}:
\begin{equation}\label{eqn:ZGW_intro}
\dgwz(X,Y) \coloneqq \frac{1}{2} \inf_\pi \left(\iint_{(X \times Y)^2} d_Z(\omega_X(x,x'),\omega_Y(y,y'))^p d\pi(x,y)d\pi(x',y')\right)^{1/p}.
\end{equation}
This defines a pseudometric on the space of $Z$-networks, and we let $(\MZp,\dgwz)$ denote the induced metric space, which we refer to as a \emph{$Z$-Gromov-Wasserstein space}; that is, $\MZp$ is the space of $Z$-networks, considered up to the equivalence relation $X \sim Y \Leftrightarrow \dgwz(X,Y) = 0$, and we abuse notation and denote the induced metric by $\dgwz$. Several basic properties of this metric space were established in~\cite{bauer2025zgromovwassersteindistance}, some of which are recalled in Section \ref{sec:GWZ_distance}.

\subsection{Main Results}

The main results of this paper relate the properties of the metric space $(Z,d_Z)$ to those of the $Z$-spaces $(\WZ,\wz)$, $(\Lp(M,Z),\dpz)$ and $(\MZp,\dgwz)$. The focus is on deriving properties of the $Z$-Gromov-Wasserstein spaces, which are the least studied among the three classes, but the connections between the various spaces also yield new results for $Z$-Wasserstein spaces and $Z$-valued Lebesgue spaces. We now informally describe the three main results of the paper.

\smallskip
\noindent {\bf Submetries Between $Z$-spaces (Theorem \ref{thm:submetries_for_GW}).} A submetry is a map between metric spaces which preserves metric balls. The existence of such a map implies close connections between the geometric structures of its domain and codomain; see Section \ref{sec:submetry_main_result} for a precise definition and further discussion. Our first main result, Theorem \ref{thm:submetries_for_GW}, establishes the existence of a chain of submetries of the form
\[\begin{tikzcd}[ampersand replacement=\&,cramped]
	{\left(\Lp(\mathbf{I}^2,Z),\frac{1}{2} \cdot \dpz\right)} \&\& {\left(\MZp,\dgwz\right)} \&\& {\left(\WZ,\frac{1}{2} \cdot \wz\right),}
	\arrow["Q", from=1-1, to=1-3]
	\arrow["E", from=1-3, to=1-5]
\end{tikzcd}\]
where $\mathbf{I}^2$ is the unit square, endowed with the Lebesgue measure. This result will be used repeatedly throughout the paper to relate the properties of $(Z,d_Z)$ to the $Z$-spaces. For example, we immediately deduce that $Z$ is complete if and only if each of the $Z$-spaces $\Lp(\mathbf{I}^2,Z)$, $\MZp$ and $\WZ$ is complete (Corollary \ref{cor:completeness_shared}), and more involved geometric conclusions are considered in the subsequent sections. The submetry
$Q:\left(\Lp(\mathbf{I}^2,Z),\frac{1}{2} \cdot \dpz\right)\to \left(\MZp,\dgwz\right) $
can be viewed as a Gromov-Wasserstein analogue of Otto’s formal Riemannian submersion picture in classical optimal transport~\cite{otto2001geometry}. The analogy is especially suggestive when $p=2$: just as Wasserstein geometry can be formally obtained as a quotient of an $L^2$-type geometry, the total space in the quotient/submetry picture is equipped with the standard $L^2$-metric.

\smallskip
\noindent {\bf Geodesic properties (Theorem \ref{thm:geodesics_main}).} We prove in Theorem \ref{thm:geodesics_main} several results on geodesics in $Z$-spaces. In the following statements, $(M,\mu)$ is an arbitrary nonatomic Polish probability space and $(Z,d_Z)$ is a Polish metric space:
\begin{itemize}
    \item For $p \in (1,\infty)$, the space $Z$ is geodesic if and only if each of the $Z$-spaces $\Lp(M,Z)$, $\MZp$, and $\WZ$ is geodesic. 
    \item For $p \in (1,\infty)$, every geodesic in $\MZp$ can be expressed in a particular, somewhat explicit form.
    \item For $p =1$, the $Z$-spaces $\Lpp{1}(M,Z)$, $\MZ{1}$ and $\WZp{1}$ are always geodesic (even if $Z$ is not).
\end{itemize}

The first part of the theorem has been established previously in the literature for Wasserstein~\cite{villani2008optimal} and nonlinear Lebesgue spaces ~\cite{jost2012nonpositive,sturm2003probability,serieys2026nonlinear}. Our proof is different, as it relies on the submetry result. One direction of this result was proved in the $Z$-GW setting (that $Z$ being geodesic implies $\MZp$ is geodesic) in our previous paper~\cite{bauer2025zgromovwassersteindistance}, and we conjectured there that the converse should hold. The second point on the structure of geodesics in the $Z$-GW space generalizes seminal work of Sturm in the classical GW setting~\cite{sturm2020}, and subsumes some previous generalizations~\cite{chowdhury2020gromov,zhang2024geometry}. The final point on geodesicity in the $p=1$ setting was previously known for the Wasserstein distance~\cite{villani2008optimal}, but is new for nonlinear Lebesgue spaces, to our knowledge, and is certainly new in the $Z$-GW setting. See Remark \ref{remark:related_work_geodesics} for more detailed explanations of connections to the existing literature.

\smallskip
\noindent {\bf Alexandrov curvature bounds (Theorem \ref{thm:curvature}).} Alexandrov curvature bounds for $Z$-spaces (for arbitrary $Z$ and $p \in [1,\infty)$) are characterized completely in terms of curvature bounds for $Z$ in Theorem \ref{thm:curvature}. The full statement of the theorem is rather lengthy, but it includes the following, where $(M,\mu)$ is a nonatomic Polish probability space and $(Z,d_Z)$ is a geodesic Polish space containing at least two points:
\begin{itemize}
    \item The $Z$-spaces $\Lpp{2}(M,Z)$, $\MZ{2}$ and $\WZp{2}$ admit lower curvature bounds if and only if $Z$ has curvature bounded below by zero, in which case the optimal lower curvature bound for each $Z$-space is zero. 
    \item For $p \neq 2$, the $Z$-spaces $\Lp(M,Z)$, $\MZp$ and $\WZ$ do not admit lower curvature bounds.
\end{itemize}
The lower curvature bounds characterizations for $\WZp{2}$ and $\Lpp{2}(M,Z)$ have previously appeared in the literature~\cite{Sturm2006geometryI,serieys2026nonlinear}, but our proof strategy is distinct, in that it derives the bounds from the submetry theorem. The result for $Z$-GW spaces is new, and provides a vast generalization of a main result of~\cite{sturm2020}, which covers the classical case of $Z=\R$. The lack of curvature bounds in the $p\neq 2$ setting seems like a basic result, at least in the Wasserstein setting, but we were not able to find a reference in the literature, and the proof uses a nontrivial construction involving the concept of ultratangent spaces. Although we do not describe them here in the introduction, for the sake of brevity, Theorem \ref{thm:curvature} also characterizes upper Alexandrov curvature bounds for the $Z$-spaces, connecting to results in~\cite{jost1994equilibrium,jost2012nonpositive,sturm2001nonlinear,bertrand2012geometric,serieys2026nonlinear,calissano2024populations}. More detailed comments on related work are collected in Remark \ref{remark:related_work_curvature}.

\subsection{Outline of the Paper}

Preliminary concepts, including formal definitions of the $Z$-spaces described above, as well as some basic results regarding them, are provided in Section \ref{sec:preliminaries}. The three main results summarized above are presented in Section \ref{sec:submetry} (on submetries), Section~\ref{sec:geodesic_properties} (on geodesics), and Section~\ref{sec:curvature} (on curvature). Each section is arranged with a similar structure: the section begins with a statement of the main result, which is then proved in pieces throughout the ensuing subsections. 

\subsection{Acknowledgements}

During the preparation of this work, the authors used ChatGPT 5 and its versions (from 5.1 through 5.6) in order to explore possible proof strategies and to identify potential gaps in preliminary arguments. The authors  take full responsibility for the full content of the  article and for the correctness of all mathematical arguments. MB and TN were supported by NSF grants DMS--2324962 and CIF--2526630. MB was in addition supported by the BSF through grant No. 2022076. FM was supported by NSF grants CCF-2523653 and DMS-2524362. We would like to thank Jesse Geneson for pointing out a mistake in our previous paper~\cite{bauer2025zgromovwassersteindistance}, as we explain in Remark~\ref{rem:separability}.

\section{Preliminaries}\label{sec:preliminaries}

This section recalls some terminology that will be used for the rest of the paper. We also formalize some concepts which were briefly described in the introduction. Before proceeding, we call attention to the following expository choice, which we observe throughout the rest of the paper.

\begin{remark}[Finite $p$]\label{rem:finite_p}
    The paper considers several metrics defined in terms of $L^p$-type integrals for $p \in [1,\infty)$. While these definitions can be extended to the $p=\infty$ case, we choose to restrict our attention to $p < \infty$ throughout the paper. The $p=\infty$ versions of the metrics are topologically and geometrically distinct from their $p < \infty$ counterparts, and extending results to $p=\infty$ generally requires alternative arguments. The subtleties which can arise when $p=\infty$ are illustrated in Remark~\ref{rem:separability}. To streamline exposition, we make a blanket restriction to $p \in [1,\infty)$. 
\end{remark}

\subsection{Wasserstein Distances}

We first recall some fundamental terminology from optimal transport theory, where we refer to the standard references~\cite{villani2008optimal,villani2021topics}. The following definitions flesh out some concepts which were already mentioned in the introduction.

\begin{definition}[Sets of measures]
    The set of Borel probability measures on a topological space $Z$ is denoted $\mathcal{P}(Z)$. For a separable metric space $(Z,d_Z)$ and $p \in [1,\infty)$, we say that $\mu \in \mathcal{P}(Z)$ has \textbf{finite $p$th moment} if 
    \[
    \int_Z d_Z(z_0,z)^p d\mu(z) < \infty
    \]
    holds for some $z_0 \in Z$ (in which case it holds for \emph{every} $z_0 \in Z$, by the triangle inequality).
    We denote the set of Borel probability measures on $Z$ with finite $p$th moment by $\WZ$. 
\end{definition}

The optimal transport distances we consider in this paper are defined in terms of measure couplings. Throughout the rest of the paper, when dealing with a product space $X \times Y$, we generically use the notations $\mathrm{proj}_1$ and $\mathrm{proj}_2$ for the \textbf{projection maps} onto the first and second coordinates, respectively. When $X$ and $Y$ are topological spaces, these maps are Borel measurable. In general, for a measurable map $\phi$, we use $\phi_\ast$ to denote the \textbf{pushforward map} on measures.

\begin{definition}[Couplings]\label{def:coupling}
Let $X$ and $Y$ be topological spaces, and let $\mu_X\in\mathcal{P}(X)$ and $\mu_Y\in\mathcal{P}(Y)$.
A \textbf{coupling} of $\mu_X$ and $\mu_Y$ is a probability measure $\pi\in\mathcal{P}(X\times Y)$ such that
$(\mathrm{proj}_1)_\ast \pi=\mu_X$ and $(\mathrm{proj}_2)_\ast \pi=\mu_Y$.
We denote the set of couplings by $\Pi(\mu_X,\mu_Y)$.
\end{definition}

We now formally define the standard notion of distance between probability measures used in optimal transport theory, which was already described in \eqref{eqn:wasserstein_distance}. In the following definition, and throughout the rest of the paper, we use $\|f\|_{L^p(\nu)}$ to denote the standard $L^p$-norm of a $p$-integrable function $f:N \to \R$ on a measure space $(N,\nu)$. 

\begin{definition}[Wasserstein distance]\label{def:wasserstein}
Let $(Z,d_Z)$ be a separable metric space and let $p\in[1,\infty)$.
For $\mu,\mu'\in\WZ$ the \textbf{Wasserstein $p$-distance} between $\mu$ and $\mu'$ is
\[
\wz(\mu,\mu')
\coloneqq \inf_{\pi \in \Pi(\mu,\mu')} \|d_Z\|_{L^p(\pi)} = \inf_{\pi\in\Pi(\mu,\mu')} \left(\displaystyle\int_{Z\times Z} d_Z(x,y)^p\,d\pi(x,y)\right)^{1/p}.
\]
\end{definition}

\subsection{Metric Space-Valued Lebesgue Spaces}

Next, we introduce some formalism surrounding the distance \eqref{eqn:Lp_distance}. Let $(M,\mu)$ be a Polish probability space, and let $(Z,d_Z)$ be a separable metric space.
For $p \in [1,\infty)$, define an extended pseudometric $\dpz$ (i.e., it may take the value $+\infty$, and it may assign distance zero to distinct maps)  between measurable maps $\omega,\omega':M \to Z$ by
\[
\dpz(\omega,\omega') \coloneqq \left(\int_{M} d_Z\bigl(\omega(x),\omega'(x)\bigr)^p\,d\mu(x)\right)^{1/p}.
\]

\begin{definition}[Metric-space-valued $\Lp$ spaces]\label{def:Lp}
Fix $z_0\in Z$ and write $\bar z_0:M\to Z$ for the constant map $\bar z_0(x)=z_0$.
We define the \textbf{nonlinear Lebesgue space} or \textbf{$Z$-valued Lebesgue space}
\[
\Lp(M,Z)\coloneqq\{\omega:M\to Z\text{ measurable}\mid \dpz(\omega,\bar z_0)<\infty\}/\!\sim,
\]
where $\omega\sim\omega'$ if $\dpz(\omega,\omega')=0$ (i.e., if $\omega = \omega'$ $\mu$-almost everywhere).
Abusing notation, the induced metric on the quotient is still denoted by $\dpz$. 
\end{definition}

\begin{remark}
    Note that the space $\Lp(M,Z)$ is independent of the choice of $z_0$; this follows directly from the triangle inequality and the fact that $\dpz(\bar z_0,\bar z_1) = d_Z(z_0,z_1) <\infty$ for any $z_1\in Z$.
\end{remark}

In this paper, we are mostly concerned with metric properties of nonlinear Lebesgue spaces. The next proposition shows that, from this perspective, the domain space $(M,\mu)$ is not centrally important. This result will, in particular, be applied in the proofs of Theorems \ref{thm:geodesics_main} and \ref{thm:curvature}.

\begin{proposition}\label{prop:nonlinear_Lebesgue_isometry}
    For Polish, nonatomic probability spaces $(M,\mu)$ and $(M',\mu')$, the nonlinear Lebesgue spaces $\Lp(M,Z)$ and $\Lp(M',Z)$ are isometric.
\end{proposition}

\begin{proof}
    By a fundamental result of measure theory (see, e.g., \cite[Theorem 17.41]{kechris2012classical}), there exists a measure-preserving bijection $\phi:M \to M'$ with measure-preserving inverse. Define $\Phi:\Lp(M,Z) \to \Lp(M',Z)$ by 
    \[
    \Phi(\omega)(x) = \omega \circ \phi^{-1}(x).
    \]
    It is easy to check that $\Phi$ is surjective. Moreover, a straightforward calculation involving the change-of-variables formula shows that $\Phi$ preserves the $\Lp$-metrics on these spaces.
\end{proof}

In particular, we frequently reduce statements about a general nonlinear Lebesgue space $\Lp(M,Z)$ to the specific case where $M$ is an interval or a unit square, endowed with Lebesgue measure. As these spaces appear frequently throughout the paper, we introduce the following notation:

\begin{convention}[Unit Interval and Square]\label{conv:unit_square}
    We use $\mathbf{I} = [0,1]$ to denote the unit interval and $\mathscr{L}$ to denote Lebesgue measure, restricted to $\mathbf{I}$. When dealing with $\mathbf{I}$, we always implicitly assume that it is endowed with the measure $\mathscr{L}$, unless specifically stated otherwise. Similarly, the unit square $\mathbf{I}^2 = [0,1] \times [0,1]$ is always assumed to be endowed with Lebesgue measure, which we consider as a product measure $\mathscr{L}^{\otimes 2} = \mathscr{L} \otimes \mathscr{L}$. 
\end{convention}

\subsection{The \texorpdfstring{$Z$}{Z}-Gromov-Wasserstein Distance}\label{sec:GWZ_distance}
Finally, we formalize some definitions regarding the $Z$-Gromov-Wasserstein distances~\eqref{eqn:ZGW_intro}, where we follow the presentation of our previous paper~\cite{bauer2025zgromovwassersteindistance}. We now fix a separable metric space $(Z,d_Z)$. 

\begin{definition}[$(Z,p)$-network]\label{def:znetwork}
A \textbf{$(Z,p)$-network} is a triple $(X,\omega_X,\mu_X)$ where $(X,\mu_X)$ is a Polish probability space
and $\omega_X:X \times X \to Z$ is a $p$-integrable function, with respect to the product measure $\mu_X \otimes \mu_X$.
\end{definition}

\begin{convention}[Notation for $Z$-networks]
 When no confusion will arise, we abuse notation and simply use $X$ to denote the entire structure $(X,\omega_X,\mu_X)$. When the particular value of $p$ is not important, we refer to $X$ simply as a \textbf{$Z$-network}.
\end{convention}

A notion of distance between $Z$-networks is given by the following generalization of the classical Gromov-Wasserstein distance.

\begin{definition}[$Z$-Gromov--Wasserstein distance]\label{def:zgw}
Let $X$ and $Y$ be $(Z,p)$-networks.
For $\pi\in\Pi(\mu_X,\mu_Y)$, we define the \textbf{$Z$-distortion} of $\pi$ by
\begin{align*}
\mathrm{dis}^Z_p(\pi)
&\coloneqq  \bigl\| d_Z\bigl(\omega_X,\omega_Y\bigr)\bigr\|_{L^p(\pi \otimes \pi)} = \left(\displaystyle\int_{X \times Y} \int_{X \times Y} d_Z(\omega_X(x,x'),\omega_Y(y,y'))^p d\pi(x,y)d\pi(x',y') \right)^{1/p}.
\end{align*}
The \textbf{$Z$-Gromov-Wasserstein $p$-distance} between $X$ and $Y$ is defined by
\begin{equation}\label{eqn:ZGW_main_def}
\dgwz(X,Y)
\coloneqq\frac12\inf_{\pi\in\Pi(\mu_X,\mu_Y)} \mathrm{dis}^Z_p(\pi).
\end{equation}
\end{definition}
The following technical result from our previous paper will be useful later on.
\begin{proposition}[{\cite[Theorem 26]{bauer2025zgromovwassersteindistance}}]\label{prop:optimal_couplings}
    For any $(Z,p)$-networks $X$ and $Y$, the infimum in
    \eqref{eqn:ZGW_main_def} is attained by some $\pi\in\Pi(\mu_X,\mu_Y)$.
\end{proposition}

Next, we define an equivalence relation on the space of $Z$-networks, whose relevance is explained in Proposition \ref{prop:metric_properties} below.

\begin{definition}[Weak isomorphism]\label{def:weak_iso}
Let $X$ and $Y$ be $Z$-networks. We define a \textbf{weak isomorphism} from $X$ to $Y$ to be a measurable map $\psi:X \to Y$ such that:
\begin{enumerate}
    \item $\psi$ preserves measure, in the sense that $\psi_\ast \mu_X = \mu_Y$, and
    \item $\psi$ preserves kernel structure, in the sense that $\psi^\ast \omega_Y = \omega_X$ almost everywhere (with respect to $\mu_X\otimes \mu_X$), where
    \[
    \psi^\ast \omega_Y(x,x') \coloneqq \omega_Y(\psi(x),\psi(x')).
    \]
\end{enumerate}
We therefore say that $X$ and $Y$ are \textbf{weakly isomorphic} if there exists a $(Z,p)$-network $W$ and weak isomorphisms $\psi_X:W \to X$ and $\psi_Y:W \to Y$. We write $X\sim Y$ to indicate that $X$ and $Y$ are weakly isomorphic; it is easy to check that $\sim$ defines an equivalence relation, and we denote the equivalence class of $X$ by $[X]$ or $[X,\omega_X,\mu_X]$, depending on what is most convenient.

We use $\MZpNoSim$ to denote the collection of all $(Z,p)$-networks and
\[
\MZp \coloneqq \MZpNoSim/\!\sim
\]
to denote the quotient by weak isomorphism; that is, $\MZp$ is the set of weak isomorphism classes of $(Z,p)$-networks.
By slight abuse of notation, we also write $\dgwz([X],[Y])\coloneqq\dgwz(X,Y)$, which is shown to be well-defined in Proposition \ref{prop:metric_properties}. The pair $(\MZp,\dgwz)$ is referred to as the \textbf{$Z$-Gromov-Wasserstein space}.
\end{definition}

The following result makes precise the metric properties of the $Z$-GW distance and contextualizes the weak isomorphism equivalence relation defined above. 

\begin{proposition}[{\cite[Theorem 29]{bauer2025zgromovwassersteindistance}}]\label{prop:metric_properties}
    The function $\dgwz$ is a pseudometric on $\MZpNoSim$ and induces a metric on the quotient $\MZp$. In particular, for $X,Y \in \MZpNoSim$, $\dgwz(X,Y) = 0$ if and only if $X \sim Y$. 
\end{proposition}

\begin{remark}[Separability of $\MZp$]\label{rem:separability}
    It was shown in \cite[Proposition 36]{bauer2025zgromovwassersteindistance} that the space $\MZp$ is separable (under the running assumption that $Z$ is separable). The proof uses~\cite[Proposition 8]{bauer2025zgromovwassersteindistance}, which claims that the nonlinear Lebesgue space $\Lp(M,Z)$ is  separable, without specifying the restriction $p < \infty$. However, this is not necessarily true for $p=\infty$, even when $Z$ is finite. Indeed, the proof presented in the paper uses the integral formulation of the metric $\dpz$, rather than the supremum norm that would naturally appear for $p = \infty$, so the result implicitly only holds for $p < \infty$. This issue illustrates the subtleties which can arise in the $p=\infty$ case, which is the reason for the $p<\infty$ convention introduced in Remark \ref{rem:finite_p}.
\end{remark}

\subsection{Examples of \texorpdfstring{$Z$}{Z}-Gromov-Wasserstein Distances and Related Work}\label{sec:variants_of_GW}

A significant portion of our previous article~\cite{bauer2025zgromovwassersteindistance} was devoted to explaining how various metrics in the extant literature can be viewed as instances of $Z$-GW distances. To give context for the present paper, we recall a few simple examples here, avoiding details and proofs. Besides the examples below, we additionally show in~\cite{bauer2025zgromovwassersteindistance} how connection graphs~\cite{robertson2025generalization}, spectral GW distances~\cite{memoli2009spectral,memoli2011spectral}, and shape graphs~\cite{guo2022statistical,bal2024statistical}, among others, can be understood within the $Z$-GW framework. 

\begin{example}[Gromov-Wasserstein Distance]
Taking $Z=\R$ with its standard metric clearly recovers the classical Gromov-Wasserstein distance \eqref{eqn:GW_distance_intro}, or, more precisely, the generalized measure network version studied in~\cite{chowdhury2019gromovwassersteindistancenetworksstable}. Choosing different metrics on $Z = \R_{\geq 0}$ (the nonnegative reals) leads to variants of the GW distance, such as the ultrametric-~\cite{memoli2023ultrametric} or $(p,q)$-Gromov-Wasserstein distance~\cite{arya2026gromov}.
\end{example}

\begin{example}[Wasserstein Distance]\label{ex:wasserstein_distance}
Given a probability distribution $\mu$ over Polish $Z$ with finite $p$th moment, one obtains a $(Z,p)$-network $(Z,\mathrm{proj}_1,\mu)$, where $\mathrm{proj}_1:Z \times Z \to Z$ is projection onto the first coordinate. We showed in~\cite{bauer2025zgromovwassersteindistance} that, for two distributions $\mu$ and $\nu$ over $Z$, the $Z$-GW distance between $(Z,\mathrm{proj}_1,\mu)$ and $(Z,\mathrm{proj}_1,\nu)$ is equal to $\frac{1}{2} \wz(\mu,\nu)$. We note that this is used in the proof of Theorem \ref{thm:submetries_for_GW} below (see Proposition \ref{prop:edge_law_section}, specifically).
\end{example}

\begin{example}[Fused Gromov-Wasserstein Distance]
In network science applications, one frequently deals with graphs which are endowed with additional node features (e.g., a social network encodes interactions between users (nodes), as well as user data (vector features on nodes)). In the case that node features are valued in a metric space $(W,d_W)$ and that edges are weighted with real numbers, such a weighted graph structure is modeled as follows. Let $X$ be the node set, let $a_X:X \times X \to \R$ be a weighted adjacency function, and let $f:X \to W$ be an attribution function. Set $Z = \R \times W$, endowed with the metric $d_Z((t,w),(t',w'))^2 = |t-t'|^2 + d_W(w,w')^2$. The weighted graph is encoded as the $Z$-network $(X,\omega_X,\mu_X)$, where  $\mu_X$ is (say) uniform measure, and 
\[
\omega_X(x,x') = \big(a_X(x,x'),f(x)\big).
\]
With these choices, the $Z$-GW distance between two such attributed graph representations agrees with the \emph{Fused Gromov-Wasserstein distance}~\cite{Vayer2020FGW}, which was designed specifically to handle such attributed structures.\footnote{To be precise, the $Z$-GW formalism recovers a  (better behaved) slight variant of the usual Fused GW distance, up to a normalization, and for a particular choice of parameter which balances the relative importance of adjacency kernel structure versus attributes. The differences are described in detail in \cite[Section 3.1.2]{bauer2025zgromovwassersteindistance}.} 

Graphs with attributes attached to \emph{edges} are also important in applications to, e.g.,  social network modeling (where edge attributes encode complex user interactions) and drug design (where a graph represents a molecule, and edge attributes encode bond types). Variants of GW distance designed to compare edge-attributed graphs were introduced in~\cite{yang2024exploiting,kawano2024multi}. We show in~\cite{bauer2025zgromovwassersteindistance} that these variants can also be recovered within the $Z$-GW framework.
\end{example}

\begin{example}[Probabilistic Metric Spaces]
A probabilistic metric space assigns to each pair \(x,x'\in X\) a
probability distribution \(p_X(x,x')\) on \(\mathbb R_{\geq0}\), subject
to axioms generalizing symmetry, vanishing on the diagonal, and the
triangle inequality~\cite{menger1942statistical,schweizer1960statistical,
kramosil1975fuzzy,wald1943statistical}.

Take
\[
    Z=\mathcal P_p(\mathbb R_{\geq0}),
\]
equipped with the \(p\)-Wasserstein metric. After equipping \(X\) with a
probability measure \(\mu_X\), every probabilistic metric space for which
\(p_X\in\mathcal L_p(X^2,Z)\) determines a \(Z\)-network
\((X,p_X,\mu_X)\). Thus probabilistic metric spaces form the subclass of
these \(Z\)-networks whose distribution-valued kernels satisfy the
probabilistic metric axioms, and \(\dgwz\) gives a natural distance
between them.
\end{example}

As was pointed out in our previous work, the $Z$-GW construction builds on similar constructions which were already present in the literature. Namely, the work of Jain and Obermayer~\cite{jain2009structure} considered a distance between metric-space-valued kernels on finite sets of fixed cardinality, and more recent work of Yang et al.~\cite{yang2024exploiting} and Kawano et al.~\cite{kawano2024multi} (already mentioned above) introduced Gromov-Wasserstein-type distances which were specifically designed to handle edge-attributed graphs. Since publishing our previous paper, we have learned about other related work, which we briefly survey below.

In \cite{lovasz2010limits}, a notion that is closely related to that of a $(Z,p)$-network is studied. Given a compact, second-countable, Hausdorff topological space $Z$, a \emph{$Z$-decorated graph} is a symmetric function $[n] \times [n] \to Z$, where $[n]\coloneqq \{1,2,\ldots,n\}$. The main results of \cite{lovasz2010limits} show that a sequence of $Z$-decorated graphs which converges, in a certain technical sense inspired by random graph sampling, can be understood to converge to a \emph{$Z$-graphon}, or a symmetric function $[0,1]\times[0,1] \to \mathcal{P}(Z)$. We observe that $Z$-decorated graphs and $Z$-graphons fall under our framework by endowing the sets $[n]$ and $[0,1]$ with uniform measures, respectively, and by endowing $\mathcal{P}(Z)$ with, say, a Wasserstein distance. The paper \cite{lovasz2010limits} shares a main motivation with our previous work \cite{bauer2025zgromovwassersteindistance}: to establish a general framework for reasoning about many generalized network models simultaneously. Besides listing common examples such as simple or weighted  graphs, \cite{lovasz2010limits} suggests examples that we did not consider in our previous work, such as multigraphs or parallel colored graphs. 

The space of $Z$-graphons is studied in \cite{abraham2023probability}, where several generalized notions of cut distance from graphon theory~\cite{lovasz2012large,Janson2013Graphons} are introduced and shown to generate the same topology. These distances share features with $Z$-Gromov-Wasserstein distances, but they are distinct; cf.~the discussion of cut distance in \cite[Section 3.1.5]{bauer2025zgromovwassersteindistance}.  Along lines similar to \cite{lovasz2010limits}, recent articles  consider similar models for graphs (or graphons) decorated in Banach spaces~\cite{kunszenti2022multigraph}, or for multi-relational graphs (or graphons)~\cite{alvarado2023limits}, with a focus on limiting properties. These works also introduced generalized cut distances in these settings.

The main results of this paper on submetries, geodesicity, and curvature bounds, have not appeared previously for any of the related general models. In fact, the results on geodesics and curvature are new even in the relatively common setting of Fused GW distances. The closest related result for Fused GW distances is from~\cite{zhang2024geometry}, where it was shown that Fused GW has nonnegative Alexandrov curvature when the attribute space is a Hilbert space; our Theorem \ref{thm:curvature} subsumes this result.

\section{Submetry Structure}\label{sec:submetry}

The main result of this section gives a precise relationship between the geometries of the $Z$-spaces 
\[
\left(\Lp(\mathbf{I}^2,Z),\dpz\right), \quad \left(\MZp,\dgwz\right), \quad \mbox{and} \quad \left(\WZ, \wz\right).
\]
Namely, we will show that these geometries are related through the concept of \emph{submetries}. This concept is recalled in the following subsection, and the main result is then formulated as Theorem \ref{thm:submetries_for_GW}. The rest of the section is devoted to proving this theorem.

\subsection{Main Result: Submetries Between \texorpdfstring{$Z$}{Z}-Spaces}\label{sec:submetry_main_result}

Let us now recall the notion of a submetry between metric spaces.

\begin{definition}[Submetry \cite{berestovskii1987submetries}]\label{def:submetry}
    For metric spaces $(M, d_M)$ and $(N, d_N)$, a mapping $f: M \to N$ is called a \textbf{submetry} if for any $x\in M$ and $r > 0$, we have
    \begin{equation}
        f(B_{d_M}(x, r)) = B_{d_N}(f(x), r),
    \end{equation}
    where $B_{d_M}(x, r)$ and $B_{d_N}(f(x), r)$ are open balls in $M$ and $N$ respectively.
\end{definition}
\begin{remark}[Surjectivity of Submetries]
Note that we do not explicitly require a submetry to be surjective; however, this follows automatically from the definition; see \cite[Proposition 3.3]{alexander2024alexandrov}.

\end{remark}

\begin{remark}[Open vs Closed Balls]
    An important technical point is that some authors define submetries to be maps which preserve closed balls, rather than open balls. We note that this is the original definition, appearing in~\cite{berestovskii1987submetries}, and this convention is also taken, e.g., in~\cite{berestovskii2000metric}. We will, however, follow~\cite{alexander2024alexandrov}, which uses the open ball convention of Definition~\ref{def:submetry}, as our main reference for properties of submetries. As is noted in~\cite{kapovitch2022structure}, the definitions coincide when the metric spaces are proper.
\end{remark}

Submetries are a purely metric analogue of submersions in Riemannian geometry---indeed, a submetry in the Riemannian category is exactly a $C^{1,1}$ submersion~\cite{berestovskii2000metric} (this paper uses the  closed ball definition of submetry, but the result applies to the open ball definition if we assume that the manifold is complete). Various metric features are preserved under submetries; for example, geodesicity, completeness, and curvature estimates all descend along submetries, as we will discuss in more detail below.

The main theorem of this section is the following, which uses the notation introduced in Convention \ref{conv:unit_square}.

\begin{theorem}[Submetries Between $Z$-Spaces]\label{thm:submetries_for_GW}
Let $(Z,d_Z)$ be a Polish metric space and let $p \in [1,\infty)$. There is a chain of submetries
\[\begin{tikzcd}[ampersand replacement=\&,cramped]
	{\left(\Lp(\mathbf{I}^2,Z), \frac{1}{2} \cdot \dpz\right)} \&\& {\left(\MZp,\dgwz\right)} \&\& {\left(\WZ,\frac{1}{2} \cdot \wz\right)}
	\arrow["Q", from=1-1, to=1-3]
	\arrow["E", from=1-3, to=1-5]
\end{tikzcd}\]
where 
\begin{itemize}
    \item $Q$ is the quotient map corresponding to an isometric action of the semigroup of measure-preserving self-maps of~$\mathbf{I}$~(see Definition \ref{def:quotient_map}), 
    \item $E$ is the edge law map, taking the equivalence class of a $Z$-network $(X,\omega_X,\mu_X)$ to the probability measure $(\omega_X)_\ast (\mu_X \otimes \mu_X)$ on $Z$ (see Definition \ref{def:edge_law}).
\end{itemize}
Moreover, there is an isometric embedding 
\[
\iota:\big(\WZ,\tfrac{1}{2} \cdot \wz\big) \hookrightarrow \big(\MZp,\dgwz\big)
\]
which provides a section of $E$, in the sense that $E \circ \iota$ is the identity on $\WZ$. 
\end{theorem}

The remainder of this section is devoted to defining the various maps more precisely and proving the theorem. The proof is carried out in parts: Corollary \ref{cor:quotient_map_submetry} shows that $Q$ is a submetry, Proposition \ref{prop:edge_law_section} shows that $E$ admits a section, and Proposition \ref{prop:edge_law_submetry} shows that $E$ is a submetry.

Before moving on to the proofs of the main theorem statements, let us state an easy corollary which illustrates the usefulness of the submetry structures. It relies on the well-known fact that completeness propagates along submetries; see, e.g.,~\cite[Proposition 3.4]{alexander2024alexandrov}. Further applications, relating the existence of geodesics and Alexandrov curvature bounds between the spaces, are provided in Sections \ref{sec:geodesic_properties} and \ref{sec:curvature}, respectively.

\begin{corollary}[Completeness Properties]\label{cor:completeness_shared}
Let $p\in[1,\infty)$ and let $(Z,d_Z)$ be a separable metric space.
The following are equivalent:
\begin{enumerate}[label=(\arabic*)]
\item $(Z,d_Z)$ is complete;
\item $\bigl(\Lp(\mathbf{I}^2,Z),\dpz\bigr)$ is complete;
\item $\bigl(\MZp,\dgwz\bigr)$ is complete;
\item $\bigl(\WZ,\wz\bigr)$ is complete.
\end{enumerate}
\end{corollary}

We note that this result could essentially be pieced together from the existing literature: (1) $\Leftrightarrow$ (2) is shown in \cite{serieys2025nonlinear}, (1) $\Leftrightarrow$ (3) is proved in our previous work~\cite{bauer2025zgromovwassersteindistance}, and (1) $\Rightarrow$ (4) is classical~\cite[Theorem 6.18]{villani2008optimal} (we were unable to find a proper reference for the converse, but it is likely known in the optimal transport community). The point of including this result is to illustrate the utility of the submetry structure in the proof.

\begin{proof}
    We first show that (1) is implied by each of the other conditions. We can consider $Z$ as a closed subspace of $\Lp(\mathbf{I}^2,Z)$ by identifying a point in $Z$ with the associated constant function, so that completeness of the Lebesgue space implies completeness of $Z$; therefore, (2) $\Rightarrow$ (1). The implication (3) $\Rightarrow$ (1) follows from \cite[Theorem 39]{bauer2025zgromovwassersteindistance} (which, in fact, proves equivalence of (1) and (3), but we only use one direction for now). Finally, that (4) $\Rightarrow$ (1) is known in the optimal transport community, but appears to be something of a folklore result (see, e.g., the MathOverflow post~\cite{CompletenessMathOverflow}); an  argument is sketched as follows. The Dirac embedding $z\mapsto \delta_z$ is an isometry from $Z$ into $(\WZ,\wz)$, and if $(z_n)$ is Cauchy in $Z$ then $(\delta_{z_n})$ is Cauchy in Wasserstein distance.
    If $(\WZ,\wz)$ is complete, this sequence converges to some $\mu$; one then checks that $\mu$ must be a Dirac mass $\delta_z$, so that $(z_n)$ converges to $z$ in $Z$.

    Now, we prove that (1) $\Rightarrow$ (2) $\Rightarrow$ (3) $\Rightarrow$ (4). The implication (1) $\Rightarrow$ (2) is classical; an argument is sketched in~\cite[Section 1.1]{korevaar1993sobolev}, and a complete argument is given in~\cite[Proposition 3.3]{sturm2001nonlinear}.
     By~\cite[Proposition 3.4]{alexander2024alexandrov} (on conservation of completeness under submetries), the fact that completeness is invariant under constant rescalings of the metric, and Theorem \ref{thm:submetries_for_GW}, we have (2) $\Rightarrow$ (3) $\Rightarrow$ (4). There is some subtlety here: assuming (2), we already showed that (1) must hold, so that $(Z,d_Z)$ is, in fact, Polish, and the theorem is then applicable to deduce (3). A similar argument gives (3) $\Rightarrow$ (4). 

\end{proof}

\subsection{The Quotient Representation of the \texorpdfstring{$Z$}{Z}-Gromov-Wasserstein Space}\label{sec:quotient_map} 
By extending a construction of Sturm in the classical GW setting~\cite[Theorem 5.10]{sturm2020} and a similar result of Chowdhury--M\'emoli in the $\mathbb{R}$-network setting \cite[Section 2.5.1]{chowdhury2019gromovwassersteindistancenetworksstable}, it was shown in \cite[Proposition 32]{bauer2025zgromovwassersteindistance} that any $Z$-network is weakly isomorphic to a $Z$-network of the form $(\mathbf{I}, \omega, \mathscr{L})$ where (as usual) $\mathscr{L}$ is the Lebesgue measure on $\mathbf{I}$ and $\omega: \mathbf{I}^2 \to Z$ is a measurable function. As this idea plays a major role in the subsequent development, let us briefly recall the construction, for a given $Z$-network $X$: by a standard result of measure theory (see, e.g., \cite[Lemma 4.2]{shioya2016metric}), there exists a measurable map $\rho:\mathbf{I} \to X$ such that $\rho_\ast \mathscr{L} = \mu_X$; then it is easy to see that the $Z$-network $(\mathbf{I},\rho^\ast \omega_X, \mathscr{L})$ is weakly isomorphic to $X$.

Let $\mathrm{Inv}(\mathbf{I})$ denote the set of measurable maps $\phi:\mathbf{I} \to \mathbf{I}$ under which the Lebesgue measure $\mathscr{L}$ is invariant. The composition operation endows $\mathrm{Inv}(\mathbf{I})$ with the structure of a \emph{semigroup}, i.e., a set with an associative binary operation (in fact, the composition operation admits an identity element, so that $\mathrm{Inv}(\mathbf{I})$ is a  \emph{monoid}, but we really only use the semigroup structure here).  We define the action of $\mathrm{Inv}(\mathbf{I})$ on $\Lp(\mathbf{I}^2, Z)$ by
\begin{align}
    (\phi^* \omega)(x, y) = \omega(\phi(x), \phi(y)),
\end{align}
for $\phi \in \mathrm{Inv}(\mathbf{I})$ and $\omega \in \Lp(\mathbf{I}^2, Z)$. This action is by isometries, and we abuse notation and denote by $\dpz$ the induced metric on the quotient space $\Lp(\mathbf{I}^2,Z)/\mathrm{Inv}(\mathbf{I})$. We note that the quotient space consists of equivalence classes of $\Lp$ functions with respect to the equivalence relation 
\[
\omega_0\sim \omega_1\Leftrightarrow \mbox{ there exists } \phi,\psi \in \mathrm{Inv}(\mathbf{I}) \mbox{ such that } \phi^*\omega_0=\psi^*\omega_1, \mathscr{L}^{\otimes 2}-\textrm{a.e.}
\]
Throughout this subsection, we use $[\omega]$ to denote the $\mathrm{Inv}(\mathbf{I})$-equivalence class of $\omega$. The induced metric $\dpz$ is defined as 
\[
\dpz([\omega_0],[\omega_1]) \coloneqq \inf_{\phi,\psi\in \mathrm{Inv}(\mathbf{I})}\dz{p}(\phi^*\omega_0,\psi^*\omega_1).
\]
Then, as shown in \cite{bauer2025zgromovwassersteindistance}, we have the following representation of the $Z$-GW space.

\begin{proposition}[{\cite[Theorem 34]{bauer2025zgromovwassersteindistance}}]
    \label{prop:zg_as_quotient_by_inv}
    Let $Z$ be a separable metric space. The embedding 
    \begin{align*}
        \Lp(\mathbf{I}^2,Z) &\hookrightarrow \MZpNoSim \\
        \omega &\mapsto (\mathbf{I},\omega,\mathscr{L})
    \end{align*}
    induces an isometry 
    \begin{equation}\label{eqn:parameterization_isomorphism}
        \Xi:\Lp(\mathbf{I}^2, Z)/\mathrm{Inv}(\mathbf{I}) \to \MZp, 
    \end{equation}
    from the quotient space $(\Lp(\mathbf{I}^2, Z)/\mathrm{Inv}(\mathbf{I}), \frac{1}{2} \cdot \dpz)$ to the $Z$-Gromov-Wasserstein space $(\MZp, \dgwz)$. That is, $\Xi$ takes the $\mathrm{Inv}(\mathbf{I})$-equivalence class of $\omega \in \Lp(\mathbf{I}^2,Z)$ to the weak isomorphism class of $(\mathbf{I},\omega,\mathscr{L}) \in \MZpNoSim$.
\end{proposition}

This leads to the following definition.

\begin{definition}[Quotient Map $Q$]\label{def:quotient_map}
    The \textbf{quotient map} 
    \[
    Q:\Lp(\mathbf{I}^2,Z) \to \MZp
    \]
    is defined to be the composition
    \[
    \Lp(\mathbf{I}^2,Z) \to \Lp(\mathbf{I}^2,Z)/\mathrm{Inv}(\mathbf{I}) \xrightarrow{\Xi} \MZp
    \]
    of the quotient by the $\mathrm{Inv}(\mathbf{I})$-action
    with the isometry \eqref{eqn:parameterization_isomorphism}.
\end{definition}

The first task in proving Theorem \ref{thm:submetries_for_GW} is to show that the quotient map is a submetry. If $\mathrm{Inv}(\mathbf{I})$ were a group, then this would follow by general principles; as this is not the case, some additional work is required. To this end, we introduce the \emph{automorphism group} $\mathrm{Aut}(\mathbf{I})$ of measure-preserving bijections with measure-preserving inverses; note that this is an honest group, rather than a semigroup. We begin with the following result, which is an analogue of~\cite[Theorem 3.2]{memoli_needham_2024_gw_gm} in the classical Gromov-Wasserstein setting, and whose proof follows the same strategy.

\begin{proposition}
    \label{prop:zg_as_quotient_by_aut}
    For any two $Z$-networks of the form $X = (\mathbf{I}, \omega_X, \mathscr{L})$ and $Y = (\mathbf{I}, \omega_Y, \mathscr{L})$, we have
    \begin{align}   
        \dgwz(X, Y) = \tfrac{1}{2} \cdot \inf_{\phi \in \mathrm{Aut}(\mathbf{I})} \dpz(\omega_X, \phi^* \omega_Y).
    \end{align}
\end{proposition}

The proof uses the following lemma. The idea is based on \cite[Theorem 1.1.(i)]{brenier2003approximation}, which shows that any coupling of the product measure $\mathscr{L}^{\otimes d}$ with itself, for $d \geq 2$, can be approximated by a sequence of measure-preserving maps $\mathbf{I}^d \to \mathbf{I}^d$, in the sense of weak convergence, as in the lemma below. Checking their proof, one observes that the maps they construct are in fact automorphisms, and that the proof strategy still works for $d=1$. However, to get a clean result, our lemma explicitly handles the $d=1$ case, which has a more straightforward proof than the higher-dimensional version of \cite{brenier2003approximation}.

\begin{lemma}\label{lem:map_approximation}
Let $\pi \in \Pi(\mathscr{L},\mathscr{L})$. There exists a sequence of maps $\phi_m \in \mathrm{Aut}(\mathbf{I})$ such that $(\mathrm{id}_\mathbf{I} \times \phi_m)_\ast \mathscr{L}$ converges weakly to $\pi$. 
\end{lemma}

\begin{proof}
    Let $m$ be a fixed integer, and throughout the rest of the proof, let $[2^m] = \{1,\ldots,2^m\}$. First, partition $[0,1)$  into intervals $(\mathbf{I}_i)_{i \in [2^m]}$, with $\mathbf{I}_i = [(i-1)/2^m, i/2^m)$. For each pair $(i,j)$ with $i,j \in [2^m]$, set $\pi_{i,j} = \pi(\mathbf{I}_i \times \mathbf{I}_j)$.

    Next, we refine the partitions. For a fixed $i_0 \in [2^m]$, partition $\mathbf{I}_{i_0}$ into half-open intervals $(\mathbf{J}_{i_0,j})_{j \in [2^m]}$ as follows:
    \[
    \mathbf{J}_{i_0,j} = [a_{j-1},a_{j}), \mbox{ where } a_0 = \frac{(i_0-1)}{2^m}, \, a_j = a_0 + \sum_{k=1}^{j} \pi_{i_0,k} \mbox{ for } j \in [2^m].
    \]
    This produces a partition, since $\sum_{k=1}^{2^m} \pi_{i_0,k} = \mathscr{L}(\mathbf{I}_{i_0})$, and has the property that $\mathscr{L}(\mathbf{J}_{i_0,j}) = \pi_{i_0,j}$.
    Note that $\mathbf{J}_{i_0,j}$ can potentially be the empty set, which we consider as a degenerate half-open interval.
    Similarly, for fixed $j_0 \in [2^m]$, partition $\mathbf{I}_{j_0}$ into half-open intervals $(\mathbf{K}_{i,j_0})_{i \in [2^m]}$ such that $\mathscr{L}(\mathbf{K}_{i,j_0}) = \pi_{i,j_0}$. 

    Finally, define $\phi_m:\mathbf{I} \to \mathbf{I}$ as follows. For each pair $(i,j)$ such that $\pi_{i,j} > 0$, the (nonempty) intervals $\mathbf{J}_{i,j}$ and $\mathbf{K}_{i,j}$ defined above have the same Lebesgue measure  $\pi_{i,j}$, so there is a unique orientation-preserving translation taking  $\mathbf{J}_{i,j}$ onto $\mathbf{K}_{i,j}$, and we define $\phi_m|_{\mathbf{J}_{i,j}}$ to be this translation. Then define $\phi_m(1) = 1$. By construction, $\phi_m \in \mathrm{Aut}(\mathbf{I})$.

    It remains to show that $(\mathrm{id}_\mathbf{I} \times \phi_m)_\ast \mathscr{L} \xrightarrow{m \to \infty} \pi$, weakly. Let $f:\mathbf{I} \times \mathbf{I} \to \R$ be a Lipschitz test function, with Lipschitz constant $L$, with respect to (say) the $\ell^\infty$-product metric on $\mathbf{I} \times \mathbf{I}$ (that is, the metric $\left((q,r),(s,t)\right) \mapsto \max\{|q-s|,|r-t|\}$, which metrizes the standard topology). To prove the claim, it suffices to show that 
    \[
    \int_{\mathbf{I} \times \mathbf{I}} f(s,t) d(\mathrm{id}_\mathbf{I}\times \phi_m)_\ast \mathscr{L}(s,t) \xrightarrow{m \to \infty} \int_{\mathbf{I} \times \mathbf{I}} f(s,t) d\pi(s,t)
    \]
    (since convergence for Lipschitz test functions is sufficient to prove weak convergence on a compact domain). Letting $t_i$ denote the midpoint of interval $\mathbf{I}_i$, we will first approximate the integrals on either side of this expression by the quantity $\sum_{i,j=1}^{2^m} f(t_i,t_j)\pi_{i,j}$.
    We have 
    \begin{align}
        &\left|\int_{\mathbf{I} \times \mathbf{I}} f(s,t) d(\mathrm{id}_\mathbf{I} \times \phi_m)_\ast \mathscr{L}(s,t) - \sum_{i,j=1}^{2^m} f(t_i,t_j)\pi_{i,j} \right| \nonumber \\
        & \qquad \qquad = \left|\sum_{i,j=1}^{2^m} \int_{\mathbf{I}_i \times \mathbf{I}_j} f(s,t) d(\mathrm{id}_\mathbf{I} \times \phi_m)_\ast \mathscr{L}(s,t) - \sum_{i,j=1}^{2^m}  f(t_i,t_j) \int_{\mathbf{I}_i \times \mathbf{I}_j} d(\mathrm{id}_\mathbf{I} \times \phi_m)_\ast \mathscr{L}(s,t) \right| \label{eqn:1d_lemma_1}\\
        & \qquad \qquad \leq \sum_{i,j=1}^{2^m} \int_{\mathbf{I}_i \times \mathbf{I}_j} |f(s,t) - f(t_i,t_j)| d(\mathrm{id}_\mathbf{I} \times \phi_m)_\ast \mathscr{L}(s,t) \nonumber \\
        &\qquad \qquad \leq  \sum_{i,j=1}^{2^m} \int_{\mathbf{I}_i \times \mathbf{I}_j} L \cdot \max\{|s-t_i|,|t-t_j|\} d(\mathrm{id}_\mathbf{I} \times \phi_m)_\ast \mathscr{L}(s,t) \label{eqn:1d_lemma_2} \\
        &\qquad \qquad \leq \frac{L}{2^{m+1}}\sum_{i,j=1}^{2^m} \int_{\mathbf{I}_i \times \mathbf{I}_j}  d(\mathrm{id}_\mathbf{I} \times \phi_m)_\ast \mathscr{L}(s,t) = \frac{L}{2^{m+1}}, \label{eqn:1d_lemma_3}
    \end{align}
    where \eqref{eqn:1d_lemma_1} uses the fact that
    \begin{equation*}
    \pi_{i,j} = \mathscr{L}(\mathbf{J}_{i,j}) = \mathscr{L}\left(\left\{t \in \mathbf{I}_i \mid \phi_m(t) \in \mathbf{I}_j \right\} \right) = (\mathrm{id}_\mathbf{I}\times \phi_m)_\ast \mathscr{L}(\mathbf{I}_i \times \mathbf{I}_j),
    \end{equation*}
    \eqref{eqn:1d_lemma_2} uses the Lipschitz property of $f$, and \eqref{eqn:1d_lemma_3} uses the bound $|s-t_i| \leq 1/2^{m+1}$ for all $s \in \mathbf{I}_i$. By a similar (in fact, slightly more straightforward) argument, we also have 
    \[
    \left|\int_{\mathbf{I} \times \mathbf{I}} f(s,t) d\pi (s,t) - \sum_{i,j=1}^{2^m} f(t_i,t_j)\pi_{i,j} \right| \leq \frac{L}{2^{m+1}}.
    \]
    By the triangle inequality, it follows that 
    \[
    \left| \int_{\mathbf{I} \times \mathbf{I}} f(s,t) d(\mathrm{id}_\mathbf{I} \times \phi_m)_\ast \mathscr{L}(s,t) - \int_{\mathbf{I} \times \mathbf{I}} f(s,t) d\pi (s,t)\right| \leq \frac{L}{2^m} \xrightarrow{m \to \infty} 0,
    \]
    completing the proof.
\end{proof}

\begin{proof}[Proof of Proposition~\ref{prop:zg_as_quotient_by_aut}]
     We clearly have $\dgwz(X, Y) \leq \tfrac{1}{2} \cdot \inf_{\phi \in \mathrm{Aut}(\mathbf{I})} \dpz(\omega_X, \phi^* \omega_Y)$, as any automorphism gives rise to a coupling, so it remains to prove the reverse inequality. By Proposition \ref{prop:optimal_couplings}, there exists a coupling $\pi \in \Pi(\mathscr{L},\mathscr{L})$ which realizes $\dgwz(X,Y)$. By Lemma \ref{lem:map_approximation}, there exists a sequence of maps $\phi_m \in \mathrm{Aut}(\mathbf{I})$ such that $(\mathrm{id}_\mathbf{I} \times \phi_m)_\ast \mathscr{L}$ weakly converges to $\pi$ as $m \to \infty$. In the proof of \cite[Theorem 26]{bauer2025zgromovwassersteindistance}, we show that the $Z$-distortion $\disp^Z:\Pi(\mathscr{L},\mathscr{L}) \to \R$ is weakly continuous. It follows that 
    \[
    \tfrac{1}{2}\cdot \dpz(\omega_X,\phi_m^\ast \omega_Y) = \tfrac{1}{2}\cdot \disp^Z\big((\mathrm{id}_\mathbf{I} \times \phi_m)_\ast \mathscr{L} \big) \xrightarrow{m \to \infty} \tfrac{1}{2}\cdot \disp^Z(\pi) = \dgwz(X,Y),
    \]
    completing the proof.
\end{proof}

Transferring this result via the isometry $\Xi$ between $\big(\Lp(\mathbf{I}^2,Z)/\mathrm{Inv},\tfrac{1}{2}\cdot \dpz\big)$ and $\big(\MZp,\dgwz\big)$  immediately yields the following corollary.

\begin{corollary}\label{cor:quotient_metric_Lp}
    The quotient metric on $\Lp(\mathbf{I}^2,Z)/\mathrm{Inv}(\mathbf{I})$ satisfies 
    \[
    \dpz([\omega],[\omega']) = \inf_{\phi \in \mathrm{Aut}(\mathbf{I})} \dpz(\omega,\phi^\ast \omega').
    \]
\end{corollary}

Finally, we obtain part of Theorem \ref{thm:submetries_for_GW} by showing that the map $Q$ is a submetry.

\begin{corollary}\label{cor:quotient_map_submetry}
    The quotient map $Q:\Lp(\mathbf{I}^2,Z) \to \MZp$ is a submetry with respect to $\tfrac{1}{2} \cdot \dpz$ and $\dgwz$.
\end{corollary}

The following proof strategy is standard when the quotient is by a group action with closed orbits (e.g.,~\cite[Proposition 3.8]{alexander2024alexandrov}). We provide details, since the action in question is only by a semigroup, and the validity of the proof relies on the characterization given in Corollary \ref{cor:quotient_metric_Lp}.

\begin{proof}
    By the construction of $Q$, it suffices to prove that the quotient map $\Lp(\mathbf{I}^2,Z) \to \Lp(\mathbf{I}^2,Z)/\mathrm{Inv}(\mathbf{I})$ is a submetry; that is, we wish  to show 
    \[
    \left[B_{\dpz}(\omega,r)\right] = B_{\dpz}([\omega],r),
    \]
    where we continue to use $[\omega]$ for the equivalence class of $\omega$ under the $\mathrm{Inv}(\mathbf{I})$-action, and we use $\left[B_{\dpz}(\omega,r)\right]$ to denote the image of the ball under the quotient map.
    
    By construction of the quotient metric, one immediately has that the quotient map is $1$-Lipschitz, which further implies the  inclusion
    \[
    \left[B_{\dpz}(\omega,r)\right] \subset B_{\dpz}([\omega],r).
    \]
    It remains to prove the reverse inclusion. Let $[\omega'] \in B_{\dpz}([\omega],r)$, i.e.,
    \[
    r > \dpz([\omega],[\omega']) = \inf_{\phi \in \mathrm{Aut}(\mathbf{I})} \dpz(\omega,\phi^\ast \omega'),
    \]
    where the equality is an application of Corollary \ref{cor:quotient_metric_Lp}. It follows that there exists some $\phi \in \mathrm{Aut}(\mathbf{I})$ such that $\dpz(\omega,\phi^\ast \omega') < r$, so that $\phi^\ast \omega' \in B_{\dpz}(\omega,r)$, and 
    \[
    [\omega'] = [\phi^\ast \omega'] \in \left[B_{\dpz}(\omega,r)\right].
    \]
    This proves that the remaining inclusion holds, so that the submetry condition is verified.
\end{proof}

\subsection{Submetry via the Edge Law}
\label{sec:relation_gw_wasserstein}

Next, we study the properties of the edge law map from Theorem \ref{thm:submetries_for_GW}. We begin with a formal definition.

\begin{definition}[Edge Law]
    \label{def:edge_law} Let $Z$ be a Polish space and let $(X,\omega_X,\mu_X)$ be a $Z$-network. The \textbf{edge law} of $(X,\omega_X,\mu_X)$ is defined as the pushforward measure $(\omega_X)_*(\mu_X \otimes \mu_X)$ on $Z$. We denote the edge law by $\nu_X$ and define the \textbf{edge law map}
    \[
    E:\MZp \to \WZ
    \]
    to be the map taking a weak isomorphism class of $(Z,p)$-networks $[X]$ to its edge law $\nu_X$. 
\end{definition}

\begin{remark}[Edge Law and Distributions of Distance]
    In the setting where $Z = \R$ and the kernel $\omega_X$ is a distance metric on $X$, the measure $\nu_X$ defined above is referred to as the  \emph{(global) distribution of distances} of $X$ (among other names). It has classically been used as a signature for comparing shapes~\cite{osada2002shape} (see also~\cite{pottmann2009integral,manay2006integral,belongie2002shape}, which utilize closely related localized descriptors) and appears in lower bounds on the standard GW distance~\cite{Memoli_2007,Memoli2011GromovWasserstein,memoli2022distance}. 
    
    We note that the edge law construction originally appeared in our previous paper,~\cite[Section 5.1]{bauer2025zgromovwassersteindistance}, but that we did not give the map a specific name. Here, the \emph{edge law} moniker is intended to evoke terminology from random graph models; interpreting the $Z$-network as a graph with edges attributed by elements of $Z$, the edge law governs the probability of attributing any particular element of $Z$ to a `random edge'. Similar random graph interpretations of this invariant appear in the recent article~\cite{gomez2025metrics}. 
    
    It is not hard to show that if $X$ and $Y$ are weakly isomorphic $(Z,p)$-networks then $\nu_X = \nu_Y$. This justifies the fact that the edge law map is well-defined as a function on $\MZp$. The converse does not hold: the edge law is not a complete invariant of isomorphism classes of $Z$-networks, or even of highly restricted subclasses thereof. For example, it is well known that sums of $n$ Dirac measures in a Euclidean space (i.e., point clouds) need not be distinguished by their distance distributions (this goes back at least to a counterexample in \cite{bloom1977counterexample},  although point clouds are generically distinguished, in a certain precise sense~\cite{boutin2004reconstructing}). Distance distributions for other classes of shapes, such as manifolds and 1-dimensional stratified spaces, are systematically studied in~\cite{memoli2022distance}.
\end{remark}

The main goal of this subsection is to prove the remaining statements of Theorem \ref{thm:submetries_for_GW}. We begin by showing that $E$ admits a section, given by the isometric embedding defined as follows:
\begin{align}
    \iota: \WZ &\to \MZp \\
    \nu &\mapsto [Z,\mathrm{proj}_1,\nu],
\end{align}
where $\mathrm{proj}_1:Z \times Z \to Z$ is projection onto the first coordinate (cf.~Example \ref{ex:wasserstein_distance}).

\begin{proposition}[A Section of the Edge Law Mapping]
    \label{prop:edge_law_section}
    The map $\iota:\WZ \to \MZp$ is an isometric embedding with respect to $\frac{1}{2} \cdot \wz$ and $\dgwz$. Moreover, it is a section of the edge law mapping $E:\MZp \to \WZ$, in the sense that $E \circ \iota$ is the identity map on $\WZ$.
\end{proposition}

\begin{proof}
    That $\iota$ is an isometric embedding was proved in \cite[Proposition 14]{bauer2025zgromovwassersteindistance}. The second statement is proved by a direct calculation: for any $\nu \in \WZ$ and any Borel subset $A \subset Z$, 
    \[
    E \circ \iota (\nu)(A) = (\mathrm{proj}_1)_\ast (\nu \otimes \nu) (A) = \nu \otimes \nu (A \times Z) = \nu(A).
    \]
\end{proof}

To complete the proof of Theorem \ref{thm:submetries_for_GW}, it remains to be shown that the edge law mapping is a submetry, which we state as the following proposition.

\begin{proposition}
    \label{prop:edge_law_submetry}
    The edge law mapping $E: (\MZp,\dgwz) \to (\WZ,\frac{1}{2} \cdot \wz)$ is a submetry.
\end{proposition}

\begin{remark}
    This is a new result even in the classical setting of $Z = \R$. It was shown in \cite[Theorem 3.1]{chowdhury2019gromovwassersteindistancenetworksstable} (with the idea going back to \cite[Corollary 6.2]{Memoli2011GromovWasserstein}) that (in the notation of this paper)
    \[
    \frac{1}{2}\cdot \mathrm{W}^\R_p(\nu_X,\nu_Y) \leq \mathrm{GW}^\R_p(X,Y),
    \]
    i.e., that the edge law map is Lipschitz. Our submetry result strengthens this classical result, while also extending it to arbitrary $Z$. 
\end{remark}

The proof will use an additional concept and result.

\begin{definition}[Synchronized Couplings]
    \label{def:synchronizable_coupling} Let $Z$ be a Polish space and let $(X,\omega_X,\mu_X)$ and $ (Y,\omega_Y,\mu_Y)$ be $Z$-networks. A coupling $\pi$ between edge laws $\nu_X$ and $\nu_Y$ is said to be \textbf{synchronized} with respect to $Z$-networks $(X,\omega_X,\mu_X)$ and $(Y,\omega_Y,\mu_Y)$ if and only if there exists a coupling $\gamma$ between $\mu_X$ and $\mu_Y$ such that $\pi = (\omega_X, \omega_Y)_* (\gamma \otimes \gamma)$, where we push forward via the map
    \begin{align*}
    (\omega_X,\omega_Y):(X \times Y) \times (X \times Y) &\to Z \times Z \\
    \big((x,y),(x',y')\big) &\mapsto \big(\omega_X(x,x'),\omega_Y(y,y')\big).
    \end{align*}
    We denote the set of synchronized couplings with respect to $X$ and $Y$ by $\Pi_{s}(\nu_X, \nu_Y)$.
\end{definition}

Said differently, a synchronized coupling $\pi$ is the law of $(\omega_X(x, x'), \omega_Y(y, y'))$ where $(x, y)$ and $(x', y')$ are i.i.d.~samples from $\gamma$. That is, the correspondence between $x$ and $y$ is synchronized with that between $x'$ and $y'$. This concept leads to an alternative formulation of the $Z$-GW distance. While the proof is simple, this observation will be useful in our proof of Proposition \ref{prop:edge_law_submetry}.

\begin{proposition}[$Z$-GW is a Constrained $Z$-Wasserstein Distance]\label{prop:zgw_is_constrained_W} The $Z$-Gromov-Wasserstein distance between $Z$-networks $X$ and $Y$ can be expressed as an optimal transport problem constrained to the synchronized couplings:
    \begin{equation}
        \dgwz(X, Y) =  
        \frac{1}{2}\inf_{\pi \in \Pi_{s}(\nu_X,\nu_Y)} \left(\int_{Z \times Z} d_Z(z, z')^p d\pi(z, z')\right)^{1/p}.
    \end{equation}
\end{proposition}

\begin{proof}
    Unraveling the definition of a synchronized coupling and applying the change-of-variables formula, we obtain
    \begin{align}
    &\frac{1}{2}\inf_{\pi \in \Pi_{s}(\nu_X, \nu_Y)}
    \left(\int_{Z \times Z} d_Z(z, z')^p \, d\pi(z, z')\right)^{1/p}
    \\
    &\qquad \qquad = \frac{1}{2}\inf_{\gamma \in \Pi(\mu_X, \mu_Y)}
    \left(\int_{Z \times Z} d_Z(z, z')^p \, d (\omega_X,\omega_Y)_\ast (\gamma \otimes \gamma)(z, z')\right)^{1/p} \\
    &\qquad \qquad  = \frac{1}{2}\inf_{\gamma \in \Pi(\mu_X, \mu_Y)}
    \left(\int_{(X\times Y)\times (X\times Y)}
        d_Z(\omega_X(x, x'), \omega_Y(y, y'))^p
        \, d(\gamma \otimes \gamma)(x, y, x', y')\right)^{1/p} \\
    &\qquad \qquad = \dgwz(X,Y).
    \end{align}
\end{proof}

We are now prepared to prove the main result of this subsection.

\begin{proof}[Proof of Proposition~\ref{prop:edge_law_submetry}]
    By Proposition \ref{prop:zgw_is_constrained_W}, for any $Z$-networks $X$ and $Y$, we have
    \begin{align}
        \frac{1}{2} \cdot \wz(E(X), E(Y)) &= \frac{1}{2}\inf_{\pi \in \Pi(\nu_X,\nu_Y)} \left(\int_{Z \times Z} d_Z(z, z')^p d\pi(z, z')\right)^{1/p} \\
        &\leq \frac{1}{2}\inf_{\pi \in \Pi_{s}(\nu_X,\nu_Y)} \left(\int_{Z \times Z} d_Z(z, z')^p d\pi(z, z')\right)^{1/p} = \dgwz(X, Y).
    \end{align}
    Therefore, we have $E(B_{\dgwz}([X], r)) \subset B_{\frac{1}{2} \cdot \wz}(E(X), r)$ for any $[X]\in \MZp$ and $r > 0$. 
    
    We will now demonstrate the remaining inclusion, $E(B_{\dgwz}([X], r)) \supset B_{\frac{1}{2} \cdot \wz}(E(X), r)$. To this end, fix a $Z$-network $(X,\omega_X,\mu_X)$, and take any $\nu \in B_{\wz}(E(X), 2r)$. We need to show that there exists a $Z$-network $Y$ such that $E(Y) = \nu$ and $\dgwz(X, Y) < r$. First, take an optimal coupling between $E(X)$ and $\nu$, denoted by $\pi^\star$. Consider the probability kernel corresponding to the conditional distribution of $\pi^\star$ given the first coordinate, denoted by $\pi_z^\star$ which exists by the disintegration theorem \cite[Theorem 3.4]{kallenberg2002foundations}. Recall that we use the notation $\mathbf{I} = [0,1]$ and $\mathscr{L}$ for Lebesgue measure on $\mathbf{I}$. By \cite[Lemma 3.2]{kallenberg2002foundations}, there exists a measurable function $F:Z\times \mathbf{I} \to Z$ such that for any $z\in Z$, $\pi_z^\star = F(z,\cdot)_\ast \mathscr{L}$. Define a $Z$-network $Y = (Y, \omega_Y, \mu_Y)$ by
    \begin{align}
        Y &= X \times [0,1], \\
        \omega_Y((x, u), (x', u')) &= F(\omega_X(x, x'), u'), \\
        \mu_Y &= \mu_X \otimes \mathscr{L}.
    \end{align}
    
    We make the following claims:

    \smallskip
    \noindent {\bf Claim 1:} The edge law of $Y$ is $\nu$.

    \smallskip
    \noindent {\bf Claim 2:} The optimal coupling $\pi^\star$ is synchronized.

    \smallskip
    Deferring the proofs of these claims momentarily, we can complete the proof. By Proposition \ref{prop:zgw_is_constrained_W} and the two claims, we have
    \begin{align}
        \dgwz(X, Y) &= \frac{1}{2}\inf_{\pi \in \Pi_s(E(X), E(Y))}\left(\int_{Z \times Z} d_Z(z, z')^p d\pi(z, z')\right)^{1/p} \\
        &= \frac{1}{2}\inf_{\pi \in \Pi_s(E(X), \nu)}\left(\int_{Z \times Z} d_Z(z, z')^p d\pi(z, z')\right)^{1/p} \\
        &= \frac{1}{2}\left(\int_{Z \times Z} d_Z(z, z')^p d\pi^\star(z, z')\right)^{1/p} = \frac{1}{2}\wz(E(X), \nu) < r.
    \end{align}
    This completes the proof that the edge law map is a submetry, modulo the proofs of the claims.

    The proof of Claim 1 follows by a lengthy calculation. For a Borel set $A \subset Z$, let $1_A:Z \to \R$ denote its indicator function. Then we have
    \begin{align}
        \nu(A) &= \int_Z 1_A(z') d\nu(z') = \int_{Z\times Z} 1_A(z') d\pi^\star (z,z') \\
        &= \int_Z \int_Z 1_A(z') d\pi^\star_z(z') dE(X)(z) \label{eqn:claim_1_1}\\
        &= \int_Z \int_Z 1_A(z') dF(z,\cdot)_\ast \mathscr{L}(z') d(\omega_X)_\ast(\mu_X \otimes \mu_X)(z) \label{eqn:claim_1_2}\\
        &= \int_Z \int_\mathbf{I} 1_A(F(z,u')) d\mathscr{L}(u') d(\omega_X)_\ast (\mu_X \otimes \mu_X)(z) \\
        &= \int_X \int_X \int_\mathbf{I} 1_A(F(\omega_X(x,x'),u')) d\mathscr{L}(u') d\mu_X(x) d\mu_X(x') \\
        &= \int_X \int_\mathbf{I} \int_X \int_{\mathbf{I}} 1_A(F(\omega_X(x,x'),u')) d\mathscr{L}(u') d\mu_X(x') d\mathscr{L}(u) d\mu_X(x)   \\
        &= \int_Y \int_Y 1_A(\omega_Y(y,y')) d\mu_Y(y) d\mu_Y(y') \\
        &= \int_Z 1_A(z) d(\omega_Y)_\ast(\mu_Y \otimes \mu_Y)(z) = E(Y)(A),
    \end{align}
    where \eqref{eqn:claim_1_1} follows by the definition of disintegration, \eqref{eqn:claim_1_2} follows by the defining property of the function $F$, and the remaining lines follow by the definition of the edge law map, Fubini's theorem, and various applications of the change of variables formula.  We remark that this also shows that $\omega_Y$ is an $\Lp$ function and thus the $Z$-network $(Y,\omega_Y,\mu_Y)$ is well-defined. The $\mathcal{L}_p$ property follows by remembering that $\nu$ has a finite $p$th moment by being in the $p$-Wasserstein space, and since the edge law is $\nu$,
    \begin{equation}
        \iint_{Y \times Y} d_Z(\omega_Y,z_0)^pd\mu_Y^{\otimes 2} = \int_{Z}d_Z(z',z_0)^pd\nu(z')<\infty.
    \end{equation}

    To prove Claim 2, consider the map 
    \begin{align*}
        g:Y &\to X \times Y \\
        (x,u) &\mapsto \big(x,(x,u)\big).
    \end{align*}
    It is not hard to show that the pushforward $g_\ast \mu_Y$ defines a coupling of $\mu_X$ and $\mu_Y$. Through another long calculation, using the same tricks as those in the proof of Claim 1, one can show that 
    \[
    \pi^\star = (\omega_X,\omega_Y)_\ast (g_\ast \mu_Y \otimes g_\ast \mu_Y),
    \]
    which proves that $\pi^\star$ is a synchronized coupling. 

\end{proof}

We end this subsection with another consequence of Proposition \ref{prop:zgw_is_constrained_W}. Whereas that result realizes $Z$-GW in terms of a Wasserstein-like distance, the following gives an expression for the $Z$-Wasserstein distance in terms of the $Z$-GW distance.

\begin{corollary}
    \label{cor:zw_is_quotient_gw} The $Z$-Wasserstein distance between probability measures $\nu$ and $\nu'$ in $\WZ$ can be expressed as a quotient of the $Z$-Gromov-Wasserstein distances through the edge law:
    \begin{equation}
        \wz(\nu, \nu') = 2 \cdot \inf\left\{\dgwz(X,Y) \mid \nu_X = \nu, \;\nu_Y = \nu' \right\}.
    \end{equation}
    Moreover, the infimum is realized by $X = \iota(\nu)$ and $Y = \iota(\nu')$, where $\iota$ is the section of the edge law map considered in Proposition \ref{prop:edge_law_section}.
\end{corollary}

\begin{proof}
        As in the proof of Proposition \ref{prop:edge_law_submetry}, it follows by Proposition \ref{prop:zgw_is_constrained_W} that the inequality 
        \begin{align}
            \wz(\nu, \nu') \leq  2 \cdot \dgwz(X, Y),
        \end{align}
        holds for any $Z$-networks $(X,\omega_X, \mu_X)$ and $(Y, \omega_Y, \mu_Y)$ with edge laws $\nu_X = \nu$ and $\nu_Y = \nu'$, respectively. Therefore, taking the infimum over all such $X$ and $Y$ proves one of the desired inequalities. On the other hand, Proposition \ref{prop:edge_law_section} implies that 
        \begin{align}
            \wz(\nu, \nu') &= 2 \cdot \dgwz((Z, \mathrm{proj}_1, \nu), (Z, \mathrm{proj}_1, \nu')) \\
            &\geq 2 \cdot \inf\left\{\dgwz(X,Y) \mid \nu_X = \nu, \;\nu_Y = \nu' \right\},
        \end{align}
        proving the desired equality. This also shows that the infimum is realized by $\iota(\nu)$ and $\iota(\nu')$. 
    \end{proof}

\section{Geodesic Properties}\label{sec:geodesic_properties}

This section is concerned with geodesic properties of the $Z$-Gromov-Wasserstein space and the related nonlinear Lebesgue and Wasserstein spaces. Here, we consider geodesics in the sense of metric geometry: a continuous path $\gamma:[0,1] \to M$ in a metric space $(M,d_M)$ is called a \emph{geodesic} if 
\begin{equation}\label{eqn:geodesic_equation}
d_M(\gamma(s),\gamma(t)) \leq (t-s)d_M(\gamma(0),\gamma(1)) \qquad \forall \quad 0 \leq s \leq t \leq 1,
\end{equation}
and the space $M$ is called \emph{geodesic} if for any $m,m' \in M$ there exists a geodesic $\gamma$ with $\gamma(0) = m$ and $\gamma(1) = m'$. We note that the inequality \eqref{eqn:geodesic_equation}, in fact, implies equality
\[
d_M(\gamma(s),\gamma(t)) = (t-s)d_M(\gamma(0),\gamma(1)) \qquad \forall \quad 0 \leq s \leq t \leq 1
\]
(see, e.g.,~\cite[Lemma 1.3]{chowdhury2018explicit}). Our main results are formulated in the following subsection as Theorem \ref{thm:geodesics_main}, with proofs provided throughout the rest of the section.

\subsection{Main Result: Geodesic Properties of \texorpdfstring{$Z$}{Z}-Spaces}\label{sec:geodesics_main_results}

Our main results of this section are summarized as follows.

\begin{theorem}[Geodesic Properties for $Z$-Spaces]\label{thm:geodesics_main}
Let $(M,\mu)$ be a Polish, nonatomic probability space and let $(Z,d_Z)$ be a Polish metric space.
\begin{enumerate}[label=(\arabic*)]
    \item (Geodesicity for $Z$-Spaces) For $1<p < \infty$ the following are equivalent:
\begin{enumerate}
\item $(Z,d_Z)$ is geodesic;
\item $\bigl(\Lp(M,Z),\dpz\bigr)$ is geodesic;
\item $\bigl(\MZp,\dgwz\bigr)$ is geodesic;
\item $\bigl(\WZ,\wz\bigr)$ is geodesic.
\end{enumerate}
\item (Structure of $Z$-GW Geodesics) If $1<p<\infty$ and $Z$ (thus also $\MZp$) is geodesic then any geodesic in  $(\MZp, \dgwz)$ can be expressed as a generalized interpolation (see Definition \ref{def:generalized_interpolation}). 
\item (Geodesicity for $p=1$) The spaces $\bigl(\Lpp{1}(M,Z),\dz{1}\bigr)$, $\bigl(\MZ{1},\dgw{1}\bigr)$ and $\bigl(\WZp{1},\wzp{1}\bigr)$ are always geodesic. 
\end{enumerate}
\end{theorem}

\begin{remark}[Related Work]\label{remark:related_work_geodesics}
We now provide context for these results. 
\begin{enumerate}[label=(\arabic*)]
    \item Various aspects of Point 1 of the theorem have previously appeared in the literature. The fact that (a) and (d) are equivalent is classical in the optimal transport literature, although one direction is typically considered only in special cases~\cite{villani2008optimal}. Early studies in nonlinear Lebesgue spaces showed that (a) implies (b), where existing results are once again typically stated under additional assumptions~\cite{jost2012nonpositive,sturm2003probability}. A proof of the equivalence of (a) and (b) (in a more analytically-flavored setting, for a more general notion of nonlinear Lebesgue space) appeared very recently in~\cite[Theorem 4.9]{serieys2026nonlinear}; the proof there of the implication (b) $\Rightarrow$ (a) is direct, whereas our proof factors through the descending implications in order to utilize Theorem \ref{thm:submetries_for_GW} to our advantage. We showed in~\cite{bauer2025zgromovwassersteindistance} that (a) implies (c). The full proof of equivalence that we present below is streamlined by the submetry theorem (Theorem~\ref{thm:submetries_for_GW}). 
     \item In the setting of classical Gromov-Wasserstein distance, Sturm gave a characterization of the structure of geodesics as  \emph{interpolations} in~\cite{sturm2020} (as we explain in more detail below---see Definition \ref{def:interpolation}). The proof strategy there has been extended in ad hoc ways to various generalized GW-type distances in~\cite{chowdhury2020gromov,zhang2024geometry,zhang2025topological}. The characterization result for $Z$-GW distances is more complicated, provides a vast generalization of Sturm's result, and the proof requires some genuinely new techniques. A similar characterization of the structure of geodesics for nonlinear Lebesgue spaces---that is, all geodesics therein are \emph{pointwise interpolations} in a certain sense---was recently proved in \cite[Theorem 4.2]{serieys2026nonlinear}. Our proof for the $Z$-GW space is independent, and requires more involved arguments.
    \item The fact that Wasserstein space is always geodesic in the $p=1$ case is shown in~\cite[Theorem 3.16]{memoli2023characterization} (following  from~\cite[Theorem 5.1]{bottou2018geometrical}); our result gives an alternative proof, by showing that the geodesic property of Wasserstein space descends from the same property for the nonlinear Lebesgue space. To our knowledge, this result has not been shown previously for nonlinear Lebesgue spaces. We conjectured in~\cite{bauer2025zgromovwassersteindistance} that the claim holds for $Z$-GW spaces, but were unable to provide a proof.
\end{enumerate}
\end{remark}

Point (2) of the theorem, on geodesic structure, has some easy corollaries. We state them here, but defer the proofs to Section \ref{subsec:geodesic_structure}. In the following statements, let $(Z,d_Z)$ be a fixed geodesic Polish space and let $p \in (1,\infty)$. 

To set up the statement of the first corollary, and to elaborate on the comment above on Sturm's characterization of geodesics in the classical GW space, we recall a construction from~\cite{bauer2025zgromovwassersteindistance}.  

\begin{definition}[Interpolation]\label{def:interpolation}
    Fix $Z$-networks $X$ and $Y$. An \textbf{interpolation} is a path of $Z$-networks of the form $X_t = (X \times Y, \omega_t, \pi)$, $t \in [0,1]$,  where $\pi$ is a $\dgwz$-optimal coupling of $\mu_X$ and $\mu_Y$, and $\omega_t$ has the property that for all $(x,y),(x',y') \in X \times Y$, the map
    \[
    t \mapsto \omega_t((x,y),(x',y'))
    \]
    is a geodesic in $Z$ between $\omega_X(x,x')$ and $\omega_Y(y,y')$. 
\end{definition} 

In \cite[Theorem 45]{bauer2025zgromovwassersteindistance}, it was shown that interpolations induce geodesics in $\MZp$. In the case $Z = \R$, this recovers Sturm's formula from~\cite{sturm2020} for geodesics in the classical GW space, where it is shown that \emph{every} geodesic can be expressed as an interpolation (when spaces are considered up to weak isomorphism). We show in Section \ref{subsec:geodesic_structure} that, for general $Z$, there may exist  geodesics which are \emph{not} of this form. However, the situation is simpler if $Z$ has unique geodesics (e.g., if it is a Hadamard space).

\begin{corollary}\label{cor:uniquely_geodesic}
    If $Z$ is a uniquely geodesic space then every geodesic in $\MZp$ (with $1 < p < \infty$) can be expressed as $[X_t]$, where $X_t$ is an interpolation.
\end{corollary}

The next corollary deals with lifting of geodesics. That is, we say that a submetry $f:M \to N$ of geodesic metric spaces has the \textbf{geodesic lifting property} if for every geodesic $\gamma:\mathbf{I} \to N$, there is a geodesic $\hat{\gamma}:\mathbf{I} \to M$ such that $f \circ \hat{\gamma} = \gamma$. If $M$ is assumed to be proper (i.e., closed metric balls are compact), then every submetry defined on $M$ has this property~\cite[Lemma 2.8]{kapovitch2022structure}. The spaces of interest in this paper are not proper, in general, but we nonetheless have:

\begin{corollary}\label{cor:geodesic_lifting}
Both submetries in the chain 
    \[\begin{tikzcd}[ampersand replacement=\&,cramped]
	{\left(\Lp(\mathbf{I}^2,Z),\frac{1}{2} \cdot \dpz\right)} \&\& {\left(\MZp,\dgwz\right)} \&\& {\left(\WZ,\frac{1}{2} \cdot \wz\right)}
	\arrow["Q", from=1-1, to=1-3]
	\arrow["E", from=1-3, to=1-5]
\end{tikzcd}\]
from Theorem \ref{thm:submetries_for_GW} have the geodesic lifting property.
\end{corollary}

The rest of this section is devoted to proving Theorem \ref{thm:geodesics_main} and its corollaries. Point (1) of the theorem, on equivalence of geodesicity is proved in Section \ref{subsec:geodesic_equivalence}. The characterization of geodesics in the $Z$-GW space (Point (2)) is proved in Section \ref{subsec:geodesic_structure}; therein, the notion of generalized interpolation is defined, and we provide examples which illustrate the necessity for our level of generality. Proofs of the corollaries are provided in Section \ref{subsec:proofs_of_corollaries}. Finally, Point (3) on universal geodesicity in the $p=1$ case is proved in Section \ref{subsec:p_equals_1}.

\subsection{Geodesic Equivalence}\label{subsec:geodesic_equivalence}

We now prove Point (1) of Theorem \ref{thm:geodesics_main}, which uses the following general lemma.

\begin{lemma}\label{lem:submetries_and_geodesics}
  Let $(M,d_M)$ and $(N,d_N)$ be metric spaces and let $f:M \to N$ be a submetry.  Suppose that for all $n_0,n_1 \in N$ there exist $m_0 \in f^{-1}(n_0), m_1 \in f^{-1}(n_1)$ such that $d_M(m_0,m_1) = d_N(n_0,n_1)$. If $M$ is geodesic then so is $N$.
\end{lemma}

\begin{remark}
    In general, one has $d_M(m_0,f^{-1}(n_1)) = d_N(n_0, n_1)$, so the condition of the lemma follows from compactness assumptions on the spaces (see \cite[Lemma 2.1]{kapovitch2022structure}). However, the spaces of interest in Theorem \ref{thm:geodesics_main} do not enjoy any useful compactness properties, so the additional assumption is necessary.
\end{remark}

\begin{proof}
    Given $n_0,n_1 \in N$, choose lifts $m_0,m_1 \in M$ satisfying the condition. Choose a geodesic $\gamma:[0,1] \to M$ joining $m_0$ to $m_1$. Then it is easy to check, using the fact that $f$ is 1-Lipschitz, that $f \circ \gamma$ is a geodesic in $N$ joining $n_0$ to $n_1$. 
\end{proof}

\begin{proof}[Proof of Theorem~\ref{thm:geodesics_main}, Point (1)]
By Proposition \ref{prop:nonlinear_Lebesgue_isometry}, geodesicity of $\Lp(M,Z)$ does not depend on the choice of $M$, so we fix $M = \mathbf{I}^2$ in the proof, without loss of generality. Moreover, in several of the arguments below, we derive the existence of geodesics from submetries. Since the property of being a geodesic space is invariant under constant scalings of the metric, we generally ignore the factors of $1/2$ which are present in the submetry theorem (Theorem \ref{thm:submetries_for_GW}).

(a)$\Rightarrow$(b): Constructions of geodesics in $\Lp(M,Z)$ using those in $Z$ have appeared in~\cite{sturm2003probability,korevaar1993sobolev,jost1994equilibrium}, typically under the assumption that $Z$ is $\mathrm{CAT}(0)$. As shown recently in~\cite[Proposition 4.4]{serieys2026nonlinear}, the ideas work in higher generality, and we sketch a proof for the sake of completeness. 

Given $\omega_0,\omega_1\in \Lp(\mathbf{I}^2,Z)$, we use Aumann's measurable selection theorem \cite[Theorem B.9]{serieys2026nonlinear} to show that there exists a measurable assignment $x\mapsto \gamma_{x}\in C([0,1],Z)$ such that $\gamma_x$ is a constant speed geodesic from $\omega_0(x)$ to $\omega_1(x)$ for $\mathscr{L}^{\otimes2}$-almost every $x$. Now define $\omega_t(x)\coloneqq\gamma_x(t)$. 
Then $t\mapsto \omega_t$ is a constant-speed geodesic in $\Lp(\mathbf{I}^2,Z)$. 

(b)$\Rightarrow$(c): By Theorem \ref{thm:submetries_for_GW}, the quotient map $Q:\Lp(\mathbf{I}^2,Z) \to \MZp$ is a submetry. Moreover, distances are realized in the sense of Lemma \ref{lem:submetries_and_geodesics}. In fact, between any $X, Y\in \MZp$, assumed to be of the form $X=(\mathbf{I},\omega_X,\mathscr{L}),Y=(\mathbf{I},\omega_Y,\mathscr{L})$ without loss of generality, take an optimal coupling $\pi$, and a Borel map $\rho:\mathbf{I}\to X\times Y$ such that $\rho_*\mathscr{L}=\pi$. Now write $\rho=(\rho_X,\rho_Y)$ and notice that $\rho_X,\rho_Y\in \mathrm{Inv}(\mathbf{I})$: e.g., 
\[
(\rho_X)_\ast \mathscr{L} = (\mathrm{proj}_1 \circ \rho)_\ast \mathscr{L} = (\mathrm{proj}_1)_\ast \rho_\ast \mathscr{L} = (\mathrm{proj}_1)_\ast \pi = \mathscr{L}.
\]
Therefore we have $Q(\rho_X^*\omega_X)=X,Q(\rho_Y^*\omega_Y)=Y$ and $\frac{1}{2}\dpz(\rho_X^*\omega_X,\rho_Y^*\omega_Y)=\dgwz(X,Y)$, so the distance is realized.

(c)$\Rightarrow$(d): By Theorem \ref{thm:submetries_for_GW}, the edge law map $E:\MZp \to \WZ$ is a submetry, and there is an isometric section $\iota:\WZ \to \MZp$. Lemma \ref{lem:submetries_and_geodesics} therefore implies that geodesicity descends along the edge law map.

(d)$\Rightarrow$(a): This is generally understood in the optimal transport community, although, to our knowledge, most references prove it in special cases (e.g., \cite[Proposition 2.10]{Sturm2006geometryI} assumes $p=2$, or the result follows from \cite[Corollary 7.22]{villani2008optimal}, which assumes that $Z$ is locally compact). We give a self-contained proof here, for completeness, which relies on the concept of midpoints: recall that a complete metric space $(M,d_M)$  is geodesic if and only if every pair of points $m,m' \in M$ admits a \emph{midpoint} $n \in M$; that is, $d_M(m,n) = d_M(m',n) = d_M(m,m')/2$ (see, e.g., \cite[Theorem 2.4.16]{BBI2001}). 
 
 Let $z_0,z_1\in Z$ and put $r\coloneqq d_Z(z_0,z_1)$.
If $(\WZ,\wz)$ is geodesic then there exists a midpoint $\mu\in\WZ$ between the Dirac measures  $\delta_{z_0}$ and $\delta_{z_1}$, i.e.
\[
\wz(\delta_{z_0},\mu)=\wz(\mu,\delta_{z_1})=\frac12\,\wz(\delta_{z_0},\delta_{z_1})=\frac r2.
\]
For $p>1$ we have $\wz(\delta_{z_0},\mu)^p=\int_Z d_Z(z_0,z)^p\,d\mu(z)$ and likewise for $z_1$.
Using the triangle inequality $r\le d_Z(z_0,z)+d_Z(z,z_1)$, followed by Minkowski's inequality, one obtains
\[
r\le \left(\int_Z \bigl(d_Z(z_0,z)+d_Z(z,z_1)\bigr)^p\,d\mu(z)\right)^{1/p}
\le \wz(\delta_{z_0},\mu)+\wz(\mu,\delta_{z_1})=r,
\]
so equality holds throughout.
The equality cases imply that $\mu$-almost every $z$ satisfies
\[
d_Z(z_0,z)=d_Z(z,z_1)=\frac r2
\qquad\text{and}\qquad
d_Z(z_0,z_1)=d_Z(z_0,z)+d_Z(z,z_1),
\]
i.e.\ $z$ is a midpoint between $z_0$ and $z_1$.
Thus $Z$ is a midpoint space, and since $Z$ is complete it follows that $Z$ is geodesic.
\end{proof}

\subsection{Structure of Geodesics in Gromov-Wasserstein Spaces}\label{subsec:geodesic_structure}
    
    The goal of this section is to completely characterize the geodesics in the $Z$-GW space. We now illustrate by example that there can be geodesics in $\MZp$ which are not given by interpolations (Definition \ref{def:interpolation}). The example uses the concept of the \textbf{size} of a $(Z,p)$-network $X$~\cite[Definition 49]{bauer2025zgromovwassersteindistance}: for fixed $z_0 \in Z$, this is the number 
    \begin{equation}\label{eqn:size}
    \mathrm{size}_{p,z_0}(X) = \|d_Z(\omega_X(\cdot,\cdot),z_0)\|_{L^p(\mu_X \otimes \mu_X)}.
    \end{equation}
    It is not hard to show that if two $(Z,p)$-networks $X$ and $Y$ are weakly isomorphic, then $\mathrm{size}_{p,z_0}(X) = \mathrm{size}_{p,z_0}(Y)$.

\begin{figure}
    \centering
    \includegraphics[width=0.75\linewidth]{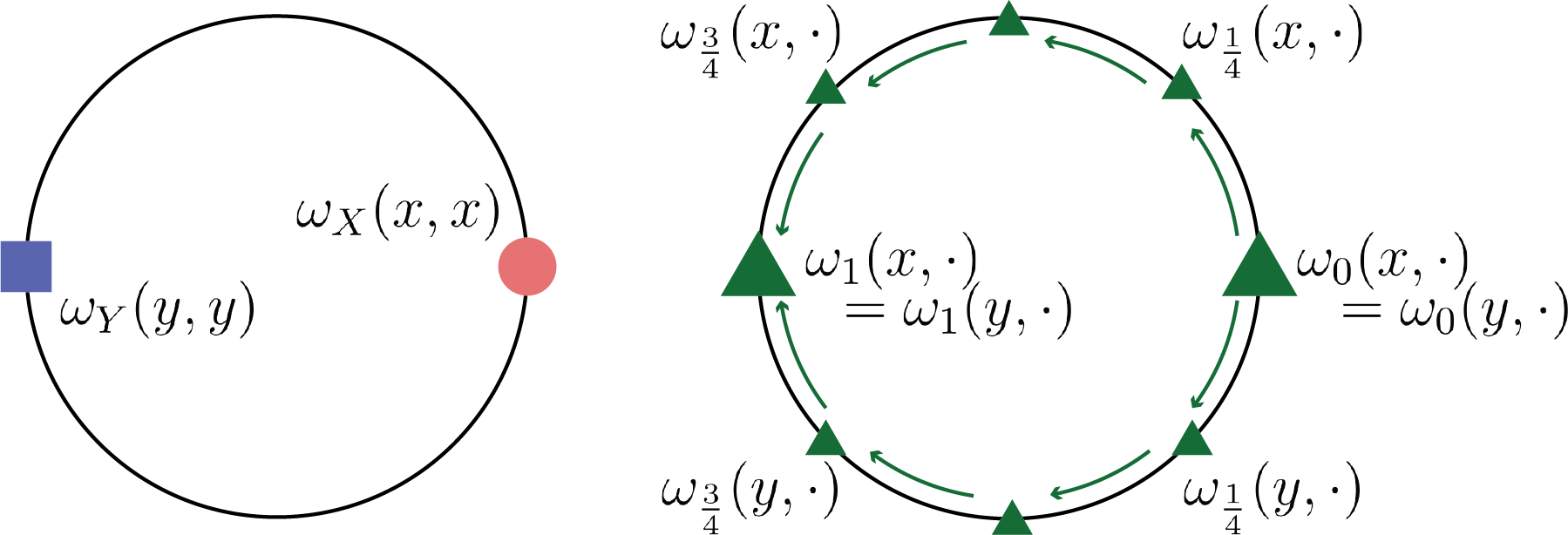}
    \caption{The geodesic path described in Example \ref{ex:circle}. The figure on the left illustrates the $Z$-networks $X$ and $Y$. The right-hand side shows some points along the path $X_t$, which describes a geodesic between $[X]$ and $[Y]$.}
    \label{fig:geodesicExample}
\end{figure}
    
    \begin{example}\label{ex:circle}
        Let $Z=S^1$, the unit circle, endowed with geodesic distance $d_Z$. For convenience of notation, we consider $Z$ as a subset of the complex plane $\mathbb{C}$.  Define $Z$-networks $X$ and $Y$ as follows:
        \begin{align*}
        X &= \{x\}, \quad \omega_X(x,x) = 1, \quad \mu_X = \delta_x, \\
        Y &= \{y\}, \quad \omega_Y(y,y) = -1, \quad \mu_Y = \delta_y.
        \end{align*}
        Consider the path $(X_t,\omega_t,\mu_t)$, $t \in [0,1]$, of $Z$-networks defined by 
        \[
        X_t = \{x,y\}, \quad \omega_t(a,b) = \left\{\begin{array}{cl}
        e^{i \pi t} & \mbox{if $a = x$}, \\
        e^{-i \pi t} & \mbox{if $a = y$},
    \end{array}\right. \quad \mu_{t} = \frac{1}{2} \delta_x + \frac{1}{2} \delta_y.
        \]
    This path is illustrated in Figure \ref{fig:geodesicExample}. It is not hard to show that $X_0$ is weakly isomorphic to $X$ and that $X_1$ is weakly isomorphic to $Y$ (e.g., the only map $X_0 \to X$ is a weak isomorphism). Moreover, the path $t \mapsto [X_t]$ defines a geodesic from $[X]$ to $[Y]$. Indeed, for $0\leq s \leq t \leq 1$, taking the identity coupling between $\mu_s$ and $\mu_t$ yields the estimate
    \begin{align*}
    2^p \dgwz(X_s,X_t)^p &\leq \sum_{a,b \in X_s} d_{Z}(\omega_s(a,b), \omega_t(a,b))^p \cdot  \frac{1}{2} \cdot \frac{1}{2} \\
    &= \frac{1}{4} \cdot 2 \cdot \big(d_Z(e^{i\pi s},e^{i \pi t})^p + d_Z(e^{-i\pi s},e^{-i \pi t})^p\big) \\
    &=  \pi^p (t-s)^p \\
    &= 2^p (t-s)^p \dgwz(X,Y)^p,
    \end{align*}
    where the last equality is an easy direct calculation, and the claim follows. On the other hand, it is straightforward to check directly that, for any $t \in (0,1)$, the path $X_t$ is not weakly isomorphic to an interpolation (in the sense of Definition \ref{def:interpolation}) between $X$ and $Y$. Indeed, setting $a = (x,y)$ for the sake of brevity, such an interpolation is of the form $\widetilde{X}_t = (\{a\},\tilde{\omega}_t,\delta_{a})$, where we either have 
    \[
    \tilde{\omega}_t(a,a) = e^{i\pi t} \; \forall \, t \quad \mbox{or} \quad \tilde{\omega}_t(a,a) = e^{-i\pi t} \; \forall \, t.
    \]
     Without loss of generality, let us assume the former case, and let us also assume that $t \in (0,1/2]$. The other cases can be handled by replacing $i$ by $-i$. We consider the sizes, relative to $z_0 = i$ (see \eqref{eqn:size}):
    \[
    \mathrm{size}_{p,i}(X_t)^p = \frac{1}{4} \sum_{a,b \in X_t} d_Z(\omega_t(a,b),i)^p = \frac{\pi^p}{4} (2 \cdot (1/2 - t)^p + 2 \cdot (1/2 + t)^p) = \frac{\pi^p}{2}((1/2 - t)^p + (1/2 + t)^p)
    \]
    while
    \[
    \mathrm{size}_{p,i}(\widetilde{X}_t)^p =  d_Z(\widetilde{\omega}_t(a,a),i)^p = \pi^p(1/2-t)^p.
    \]
    As the sizes are unequal for $t > 0$, it must be that $X_t$ and $\widetilde{X}_t$ are not weakly isomorphic. 
    \end{example}

The issue in the previous example is that the target space $Z$ admits distinct geodesics between antipodal points. To fully characterize the structure of geodesics in $\MZp$, a slight generalization of the notion of interpolation will suffice.

\begin{definition}[Generalized Interpolation]\label{def:generalized_interpolation}
    Fix $Z$-networks $X$ and $Y$. A \textbf{generalized interpolation} is a path of $Z$-networks of the form  $X_t = (X \times Y \times \mathbf{I}, \omega_t, \pi \otimes \mathscr{L})$, $t \in [0,1]$, 
    such that $\pi$ is a $\dgwz$-optimal coupling of $X$ and $Y$, and $\omega_t$ has the property
    that for a.e.\,$(x,y,s),(x',y',s') \in X \times Y \times \mathbf{I}$, the map
    \[
    t \mapsto \omega_t((x,y,s),(x',y',s'))
    \]
    is a geodesic in $Z$ between $\omega_X(x,x')$ and $\omega_Y(y,y')$.  To be precise, the a.e.\,qualifier means that there exists $A\subset (X\times Y\times \mathbf{I})^2$ with $(\pi \otimes \mathscr{L})^{\otimes 2}(A)=1$ such that for every $((x,y,s),(x',y',s'))\in A$, the map $t\mapsto \omega_t((x,y,s),(x',y',s'))$ is a constant speed geodesic from $\omega_X(x,x')$ to $\omega_Y(y,y')$. Here, $(\pi \otimes \mathscr{L})^{\otimes 2} = (\pi \otimes \mathscr{L})\otimes (\pi \otimes \mathscr{L})$.
\end{definition}

Intuitively, generalized interpolations are similar to interpolations (in the sense of Definition \ref{def:interpolation}), but they allow mass to be split along multiple geodesics in $Z$ in the kernel $\omega_t$. The extra factor of $\mathbf{I}$ in the underlying set $X \times Y \times \mathbf{I}$ of the generalized interpolation is included to keep track of how the mass splits. This idea is illustrated in the following example.

 \begin{example}
        We now return to the path $X_t$ defined in Example \ref{ex:circle}. It was shown that $X_t$ is not weakly isomorphic to an interpolation, in the sense of Definition \ref{def:interpolation}, but we now show that it is weakly isomorphic to a generalized interpolation for each $t \in [0,1]$, which we denote by $\widehat{X}_t$, in order to distinguish it from the previously defined path. The generalized interpolation is given by
        \[
        \widehat{X}_t = \left(\{(x,y)\} \times \mathbf{I}, \widehat{\omega}_t, \delta_{(x,y)} \otimes \mathscr{L}\right),
        \]
        where, setting $a = (x,y)$ for the sake of brevity, we have 
        \[
        \widehat{\omega}_t((a,s),(a,s')) = \left\{\begin{array}{rl}
        e^{i\pi t} &\mbox{if $s < 1/2$} \\
        e^{-i\pi t} &\mbox{if $s \geq 1/2$}.
        \end{array}\right.
        \]
        For fixed $t \in [0,1]$, define a map $\psi:\widehat{X}_t \to X_t$
        by 
        \[
        \psi((a,s)) = \left\{
        \begin{array}{rl}
        x & \mbox{if $s < 1/2$} \\
        y & \mbox{if $s \geq 1/2$.}
        \end{array}\right.
        \]
        Then 
        \[
        \psi_\ast (\delta_{(x,y)} \otimes \mathscr{L})(\{x\}) = \delta_{(x,y)} \otimes \mathscr{L}(\{a\} \times [0,1/2)) = 1/2,
        \]
        and likewise $\psi_\ast (\delta_{(x,y)} \otimes \mathscr{L})(\{y\}) = 1/2$, so that $\psi$ is measure-preserving. Moreover, 
        \[
        \psi^\ast \omega_t((a,s),(a,s')) = \omega_t(\psi(a,s),\psi(a,s')) = \left\{
        \begin{array}{rl}
        e^{i \pi t} & \mbox{if $s < 1/2$} \\
        e^{-i\pi t} & \mbox{if $s \geq 1/2$}
        \end{array}\right. = \widehat{\omega}_t((a,s),(a,s')).
        \]
        This proves that $\psi$ is a weak isomorphism. 
    \end{example}

We now show that generalized interpolations induce geodesics in the $Z$-GW space. 

\begin{proposition}
    A generalized interpolation $X_t$ between $(Z,p)$-networks $X$ and $Y$ induces a geodesic $[X_t]$ between $[X]$ and $[Y]$ in $\MZp$. 
\end{proposition}

\begin{proof}
First, we show that $X_0$ is weakly isomorphic to $X$. Let $\mathrm{proj}_1:X \times Y \times \mathbf{I} \to X$ be the projection map. Then $\mathrm{proj}_1$ is measure-preserving, since $\pi$ is a coupling. Moreover, 
\[
\mathrm{proj}_1^\ast \omega_X((x,y,s),(x',y',s'))= \omega_X(\mathrm{proj}_1(x,y,s),\mathrm{proj}_1(x',y',s')) = \omega_X(x,x') = \omega_0((x,y,s),(x',y',s')).
\]
Therefore,  $\mathrm{proj}_1$ is a weak isomorphism. Similarly, $X_1$ is weakly isomorphic to $Y$.

Now, let $0 \leq s \leq t \leq 1$. Taking the identity coupling between $\pi \otimes \mathscr{L}$ and itself gives the estimate
\begin{align*}
    &2^p \cdot \dgwz(X_s,X_t)^p \\
    &\leq \iint_{(X \times Y \times \mathbf{I})^2} d_Z\big(\omega_s((x,y,u),(x',y',u')), \omega_t((x,y,u),(x',y',u')\big)^p \, d(\pi \otimes \mathscr{L})(x,y,u) d(\pi \otimes \mathscr{L})(x',y',u') \\
    &= \iint_{(X \times Y)^2} (t-s)^p d_Z\big(\omega_X(x,x'),\omega_Y(y,y')\big)^p \, d\pi(x,y) d\pi(x',y') \iint_{\mathbf{I}^2} d\mathscr{L}(u) d\mathscr{L}(u') \\
    &= (t-s)^p \cdot 2^p \cdot \dgwz(X,Y)^p,
\end{align*}
which completes the proof.
\end{proof}

The next result shows that \emph{every} geodesic in the $Z$-GW space is a generalized interpolation. This proves Point (2) of Theorem \ref{thm:geodesics_main}.

\begin{proposition}\label{prop:unique_geodesics}
    Let $X$ and $Y$ be $Z$-networks and let $X_t$ be a path of $Z$-networks which induces a geodesic in $\MZp$ between $[X]$ and $[Y]$. Then there exists a generalized interpolation $\widehat{X}_t$ such that $[X_t] = [\widehat{X}_t]$ for all $t \in \mathbf{I}$. 
\end{proposition}

We use the following generalized version of Minkowski's inequality.

\begin{lemma}
    \label{lem:minkowski_multiple}
    Let $p \in (1,\infty)$, let $(\Omega, \mathcal{F}, \mu)$ be a probability space, and let $f_1, f_2, \cdots, f_n: \Omega \to [0,\infty)$ be nonnegative, $p$-integrable functions. Then we have
    \begin{equation}
        \left\|\sum_{i=1}^n f_i\right\|_{L^p(\mu)} = \sum_{i=1}^n \|f_i\|_{L^p(\mu)}
    \end{equation}
    if and only if there exist constants $c_i \geq 0, i =1,\cdots, n$ and a common nonnegative $p$-integrable function $h:\Omega \to [0, \infty)$ such that $f_i = c_i h$ $\mu$-almost everywhere for every $i \in \{ 1,\cdots, n\}$.
\end{lemma}

\begin{proof}
    We focus on the nontrivial direction of the statement. The proof is by induction. The case $n=2$ is exactly the equality case of Minkowski's inequality. Suppose the statement is true for $n=k$. We will prove it for $n=k+1$. By the triangle inequality and the assumption,
    \begin{equation}
        \left\|\sum_{i=1}^{k+1} f_i\right\|_{L^p(\mu)} = \left\|\sum_{i=1}^{k} f_i + f_{k+1}\right\|_{L^p(\mu)} \leq  \left\|\sum_{i=1}^{k} f_i\right\|_{L^p(\mu)} + \|f_{k+1}\|_{L^p(\mu)} \leq \sum_{i=1}^{k+1} \|f_i\|_{L^p(\mu)} = \left\|\sum_{i=1}^{k+1} f_i\right\|_{L^p(\mu)}
    \end{equation}
    Therefore, all the inequalities above are equalities. By the equality case of Minkowski's inequality, we either have $\sum_{i=1}^{k}f_i=0$ $\mu$-almost everywhere, or there exists a constant $c_1 \geq 0$ and a function $h_1:\Omega \to \mathbb{R}$ such that $\sum_{i=1}^{k} f_i = h_1$ and $f_{k+1} = c_1 h_1$ $\mu$-almost everywhere. We first handle the case $\sum_{i=1}^{k}f_i=0$ $\mu$-almost everywhere. In this case, we have $f_i=0$ almost everywhere since $f_i$ are nonnegative. Therefore, the lemma for $n=k+1$ holds by taking $c_i=0,i=1,\cdots,k, c_{k+1}=1$ and $h=f_{k+1}$. Now for the other case, by the induction hypothesis, there exist constants $c_i \geq 0$ and a function $h_2:\Omega \to [0, \infty)$ such that $f_i = c_i h_2$ $\mu$-almost everywhere for any $i=1,\cdots, k$. Therefore, we have $h_1 = \sum_{i=1}^k f_i = (\sum_{i=1}^{k}c_i) h_2$ $\mu$-almost everywhere and thus $f_{k+1} = c_1(\sum_{i=1}^{k}c_i)h_2$ $\mu$-almost everywhere, so the lemma follows.
\end{proof}

\begin{proof}[Proof of Proposition \ref{prop:unique_geodesics}]
    It will be convenient to use the notation $X_0 = (X, \omega_X, \mu_X)$, $X_1 = (Y, \omega_Y, \mu_Y)$, and $L = \dgwz(X_0, X_1)$. We assume that $L > 0$, as the claim is otherwise trivial. The proof is somewhat lengthy, and we break it into parts. The overall strategy is to construct an auxiliary path
    \[
    \widetilde{X}_t = (\Omega, \tilde{\omega}_t, \Gamma),
    \]
    and to then show that, for every $t \in [0,1]$, 
    \[
    [X_t] = \left[\widetilde{X}_t\right] = \left[\widehat{X}_t\right],
    \]
    where $\widehat{X}_t$ is a generalized interpolation.

\smallskip
    \noindent \textbf{Step 1: Structure of the Geodesic $X_t$.} We begin by establishing some structural properties of the geodesic $[X_t]$. These are mainly technical preliminaries, and will be used later in the proof.

    We first fix an integer $k$ and consider dyadic rationals in the set $T^{(k)} \coloneqq \{t_i \coloneqq i/2^k|i=0,1,\cdots, 2^k\}$. Let $\rho_{t_i,t_{i+1}}$ be an optimal coupling between adjacent dyadic times $t_i,t_{i+1} \in T^{(k)}$, and define $\Gamma^{(k)} \in \mathcal{P}\left(\prod_{i=0}^{2^k} X_{t_{i}}\right)$ to be a gluing of all these optimal couplings. We then define $\pi^{(k)}$ to be the marginal of $\Gamma^{(k)}$ on $X_0 \times X_1$, which is a coupling between $\mu_X$ and $\mu_Y$. Then, by triangle inequality and Minkowski's inequality, we have
    \begin{align}
        2\,L &\leq \|d_Z(\omega_X, \omega_Y)\|_{L^p(\pi^{(k)}\otimes\pi^{(k)}}) = \|d_Z(\omega_X, \omega_Y)\|_{L^p(\Gamma^{(k)} \otimes \Gamma^{(k)})} \\
        &\leq \left\|\sum_{i=0}^{2^k - 1} d_Z(\omega_{t_{i}}, \omega_{t_{i+1}})\right\|_{L^p(\Gamma^{(k)}\otimes \Gamma^{(k)})} \\        
        &\leq \sum_{i=0}^{2^k - 1} \|d_Z(\omega_{t_{i}}, \omega_{t_{i+1}})\|_{L^p(\rho_{t_{i}, t_{i+1}}\otimes \rho_{t_{i}, t_{i+1}})} \\
        &= 2\sum_{i=0}^{2^k - 1} \dgwz(X_{t_{i}}, X_{t_{i+1}})
        = 2\sum_{i=0}^{2^k - 1} \frac{L}{2^k}
        = 2\,L
    \end{align}
    Therefore, all the inequalities above are equalities and $\pi^{(k)}$ is an optimal coupling. In particular, we have
    \begin{equation}\label{eqn:uniqueness1}
       d_Z(\omega_X, \omega_Y) = \sum_{i=0}^{2^k - 1} d_Z(\omega_{t_{i}}, \omega_{t_{i+1}}) \quad \Gamma^{(k)} \otimes \Gamma^{(k)} \textrm{ a.e.} 
    \end{equation}
    and
    \begin{equation}\label{eqn:uniquness2}
        \left\|\sum_{i=0}^{2^k - 1} d_Z(\omega_{t_{i}}, \omega_{t_{i+1}})\right\|_{L^p(\Gamma^{(k)}\otimes \Gamma^{(k)})} = \sum_{i=0}^{2^k - 1} \|d_Z(\omega_{t_{i}}, \omega_{t_{i+1}})\|_{L^p(\rho_{t_{i}, t_{i+1}}\otimes \rho_{t_{i}, t_{i+1}})}.
    \end{equation}
    The equality \eqref{eqn:uniquness2} invokes the equality condition of Lemma \ref{lem:minkowski_multiple}, and it follows that there exists a common nonnegative function $h:\left(\prod_{i=0}^{2^k} X_{t_{i}}\right)^2 \to [0,\infty)$ and scalars $c_i \geq 0$,  $i=0,1,\cdots, 2^k - 1$, such that  
    \begin{equation}
        d_Z(\omega_{t_{i}}(x_{t_i}, x_{t_i'}), \omega_{t_{i+1}}(x_{t_{i+1}}, x_{t_{i+1}}')) = c_i h(x_{t_0}, x_{t_1}, \cdots, x_{t_{2^k}}, x_{t_0}', x_{t_1}', \cdots, x_{t_{2^k}}')
    \end{equation}
    for $\Gamma^{(k)} \otimes \Gamma^{(k)}$-a.e. $(x_{t_0}, x_{t_1}, \cdots, x_{t_{2^k}}, x_{t_0}', x_{t_1}', \cdots, x_{t_{2^k}}')$. Combining this with \eqref{eqn:uniqueness1}, we have
    \begin{equation}\label{eqn:uniqueness3}
        d_Z(\omega_X(x_0, x_0'), \omega_Y(x_1, x_1')) = \left(\sum_{i=0}^{2^k - 1} c_i\right) h(x_{t_0}, x_{t_1}, \cdots, x_{t_{2^k}}, x_{t_0}', x_{t_1}', \cdots, x_{t_{2^k}}')
    \end{equation}
    for $\Gamma^{(k)} \otimes \Gamma^{(k)}$-a.e. $(x_{t_0}, x_{t_1}, \cdots, x_{t_{2^k}}, x_{t_0}', x_{t_1}', \cdots, x_{t_{2^k}}')$. If $\sum_{i=0}^{2^k - 1} c_i = 0$ then $d_Z(\omega_X, \omega_Y) = 0$ $\pi^{(k)} \otimes \pi^{(k)}$-a.e., which means that $X$ and $Y$ are weakly isomorphic, i.e., that we are in the excluded case $L = 0$. Therefore, we assume $\sum_{i=0}^{2^k - 1} c_i > 0$, and we have 
    \begin{equation}\label{eqn:uniqueness4}
        d_Z(\omega_{t_{i}}(x_{t_i}, x_{t_i}'), \omega_{t_{i+1}}(x_{t_{i+1}}, x_{t_{i+1}}')) = \frac{c_i}{\sum_{j=0}^{2^k - 1} c_j} d_Z(\omega_X(x_0, x_0'), \omega_Y(x_1, x_1'))
    \end{equation}
    for $\Gamma^{(k)} \otimes \Gamma^{(k)}$-a.e. $(x_{t_0}, x_{t_1}, \cdots, x_{t_{2^k}}, x_{t_0}', x_{t_1}', \cdots, x_{t_{2^k}}')$. Taking \(L^p(\Gamma^{(k)}\otimes\Gamma^{(k)})\)-norms on both sides, we have
    \begin{equation}
        \dgwz(X_{t_{i}}, X_{t_{i+1}}) =  \frac{c_i}{\sum_{j=0}^{2^k - 1} c_j}L.
    \end{equation}
    On the other hand, $\dgwz(X_{t_{i}}, X_{t_{i+1}}) = (t_{i+1} - t_i)L$, by the assumption that $X_t$ induces a geodesic. Since $L > 0$, we can put all of this together to deduce that   $c_i = (t_{i+1} - t_i)\sum_{j=0}^{2^k - 1} c_j$. Plugging this back into \eqref{eqn:uniqueness4}, we have
    \begin{equation}
        \label{eqn:uniqueness5}
        d_Z(\omega_{t_{i}}(x_{t_i}, x_{t_i}'), \omega_{t_{i+1}}(x_{t_{i+1}}, x_{t_{i+1}}')) = (t_{i+1} - t_i) d_Z(\omega_X(x_0, x_0'), \omega_Y(x_1, x_1'))
    \end{equation}
    for $\Gamma^{(k)} \otimes \Gamma^{(k)}$-a.e. $(x_{t_0}, x_{t_1}, \cdots, x_{t_{2^k}}, x_{t_0}', x_{t_1}', \cdots, x_{t_{2^k}}')$. 
    
    The main conclusion of Step 1 of the proof is that the expression \eqref{eqn:uniqueness5} implies that samples $\{x_{t_i}\}_{i=0}^{2^k}$ and $\{x_{t_i}'\}_{i=0}^{2^k}$ from the probability distribution $\Gamma^{(k)}$ define a path of points $\{\omega_{t_i}(x_{t_i}, x_{t_i}')\}_{i=0}^{2^k}$ along a geodesic between $\omega_X(x_0, x_0')$ and $\omega_Y(x_1, x_1')$.

\smallskip
    \noindent \textbf{Step 2: Definition of $\Omega$ and $\Gamma$.} Our goal is now to construct the auxiliary path $\widetilde{X}_t = (\Omega, \tilde{\omega}_t, \Gamma)$ mentioned at the beginning of the proof. In this step, we define the probability space $(\Omega,\Gamma)$. 
    
    Let $\mathbb{D} \coloneqq \cup_{k \geq 1} T^{(k)}$ denote the set of all dyadic rationals between $0$ and $1$. We define  $\Omega$  to be the product space
    \[
    \Omega \coloneqq \prod_{t \in \mathbb{D}} X_t.
    \]
It remains to construct the measure $\Gamma$, which will be realized as a projective limit. Specifically, we follow the formulation of~\cite[Theorem 8.23]{kallenberg2002foundations}, which we now briefly recall by outlining our strategy. We will construct a family of probability measures $\Gamma_J$ on the spaces $S_J \coloneqq \prod_{t \in J} X_t$, indexed by finite subsets $J \subset \mathbb{D}$. The measures are  required to satisfy the \emph{projective condition}, that is, for all $I \subset J$, $(\mathrm{proj}_{S_J,S_I})_\ast \Gamma_J = \Gamma_I$, where $\mathrm{proj}_{S_J,S_I}:S_J \to S_I$ is the coordinate projection. Given such a family, \cite[Theorem 8.23]{kallenberg2002foundations} guarantees the existence of a measure $\Gamma$ on $\Omega$ with the property that $(\mathrm{proj}_{\Omega,S_J})_\ast \Gamma = \Gamma_J$ for all finite $J \subset \mathbb{D}$.

Fix a finite subset $J \subset \mathbb{D}$ and an integer $k$ such that $J \subset T^{(k)}$. For the measure $\Gamma^{(k)}$ on $S_{T^{(k)}}$ constructed in Step 1, let $\Gamma^{(k)}_J \coloneqq (\mathrm{proj}_{S_{T^{(k)}},S_J})_\ast \Gamma^{(k)}$. We make two useful observations about measures of this form:
\begin{enumerate}[label=(\arabic*)]
    \item[(Obs 1)] For a subset $I \subset J$, the projective condition holds; i.e., $(\mathrm{proj}_{S_J,S_I})_\ast \Gamma^{(k)}_J = \Gamma^{(k)}_I$.
    \item[(Obs 2)] The sequence $(\Gamma_J^{(k)})_{k \geq 1}$ is a tight family of measures on $S_J$. This follows because, for any $t \in J$, the $X_t$-marginal of any such measure is $\mu_t$, and the set of measures with this property (i.e., couplings of the set $\{\mu_t\}_{t \in J}$) is tight. 
\end{enumerate}

The goal is now to remove the dependence on $k$ in our constructed measures. By (Obs 2) and Prokhorov's Theorem on tight families of measures (e.g.,~\cite[Chapter 4]{villani2008optimal}), each $(\Gamma_J^{(k)})_{k}$ has a convergent subsequence, but we need to be a bit more careful to ensure the projective condition for the family; to get subsequences which converge in a consistent way for different choices of $J$, we use the following diagonalization argument. We first enumerate all finite subsets of $\mathbb{D}$ as $J_1,J_2,\ldots$ (this is possible because $\mathbb{D}$ is countable). Then, we define a nested family of subsequences inductively:
\begin{enumerate}
    \item[(0)] Set $N_0 = (1,2,\ldots)$.
    \item[(1)] By Prokhorov's Theorem, there is a convergent subsequence $\Gamma_{J_1}^{(k_n^1)} \to \Gamma_{J_1}$. Set $N_1 = (k_1^1,k_2^1,\ldots)$. 
    \item[(2)] The collection $(\Gamma_{J_2}^{(k)})_{k \in N_1}$ still forms a tight family, and therefore admits a convergent subsequence $\Gamma_{J_2}^{(k_n^2)} \to \Gamma_{J_2}$. Set $N_2 = (k_1^2,k_2^2,\ldots)$. By construction,  $N_2 \subset N_1$, so that $(\Gamma_{J_1}^{(k)})_{k \in N_2}$ still converges to $\Gamma_{J_1}$.
    \item[(i)] The collection $(\Gamma_{J_i}^{(k)})_{k \in N_{i-1}}$  admits a convergent subsequence $\Gamma_{J_i}^{(k_n^i)} \to \Gamma_{J_i}$, and we set $N_i = (k_1^i,k_2^i,\ldots)$. By construction,  $N_{i} \subset N_{i-1}$, so that $(\Gamma_{J_{i-1}}^{(k)})_{k \in N_i}$ still converges to $\Gamma_{J_{i-1}}$. 
\end{enumerate}
This construction defines a nested family of sequences $N_0 \supset N_1 \supset N_2 \supset \cdots$, which, by induction, has the property that $(\Gamma_{J_j}^{(k)})_{k \in N_i}$ converges to $\Gamma_{J_j}$ whenever $i \geq j$. Finally, we diagonalize and set $k_n = k_n^n$, so that $k_n \in N_m$ for all $n \geq m$, and $\Gamma_{J_m}^{(k_n)} \to \Gamma_{J_m}$ for all $m$. We remark that $\Gamma_J^{(k)}$ is only defined when $J\subset T^{(k)}$, so we need to remove finitely many terms from $k_n$ for a rigorous convergence statement.
By (Obs 1) and the continuity of the projection map, it follows that the projective condition holds for the collection $\{\Gamma_J\}_J$ by taking the limit in $k$ along the subsequence defined by $k_n$. We can therefore take $\Gamma$ to be the measure on $\Omega$ guaranteed by \cite[Theorem 8.23]{kallenberg2002foundations}, so that this step of the proof is complete.

    \smallskip
    \noindent \textbf{Step 3: Definition of  $\tilde{\omega}_t$.} Next, we complete the construction of the auxiliary path $\widetilde{X}_t = (\Omega, \tilde{\omega}_t, \Gamma)$ by defining the kernel $\tilde{\omega}_t: \Omega \times \Omega \to Z$, for $t \in [0,1]$. For the rest of the proof, we use the notation $\mathbf{x} = (x_t)_{t \in \mathbb{D}}$ for an element of $\Omega = \prod_{t \in \mathbb{D}} X_t$. For $t\in \mathbb{D}$, we define the kernel straightforwardly as 
    \begin{equation}
    \tilde{\omega}_t(\mathbf{x},\mathbf{x}')\coloneqq \omega_t(x_t,x_t').
    \end{equation}
    
    In order to define $\tilde{\omega}_t$ for $t \in [0,1] \setminus \mathbb{D}$, we will use an abstract extension argument. For our purposes, it suffices to extend the definition $\Gamma \otimes \Gamma$-almost everywhere. Our strategy is to show that, for $\Gamma \otimes \Gamma$-almost every pair $(\mathbf{x},\mathbf{x}')$, it holds that 
    \begin{equation}\label{eqn:Lipschitz_condition_for_omega}
       d_Z(\tilde{\omega}_{s}(\mathbf{x},\mathbf{x}'), \tilde{\omega}_{t}(\mathbf{x},\mathbf{x}')) = |s-t|d_Z(\omega_X(x_0, x_0'), \omega_Y(x_1, x_1')),
    \end{equation}
    for every $s,t \in \mathbb{D}$. 
    Since $Z$ is complete and dyadic rationals are dense in $[0,1]$, we can then  extend $\tilde{\omega}_{t}(\mathbf{x},\mathbf{x}')$ to all $t\in [0,1]$ uniquely via a general Lipschitz extension result~\cite[Proposition 1.5.9]{BBI2001}. 
    
    We now proceed to establish that \eqref{eqn:Lipschitz_condition_for_omega} holds almost everywhere, using \eqref{eqn:uniqueness5} from Step 1 of the proof. We first note that the Minkowski argument in Step 1 showed that the marginal $\pi^{(k)}$ of $\Gamma^{(k)}$ on $X_0\times X_1$ is optimal. The argument applies verbatim to marginals at other times, so for any $r,s\in T^{(k)}$, the marginal $\Gamma_{\{r,s\}}^{(k)}$ on $X_r\times X_s$ is an optimal coupling between $\mu_r$ and $\mu_s$. Since the distortion functional is continuous for $1<p<\infty$ \cite[Lemma 27 and Proof of Theorem 26]{bauer2025zgromovwassersteindistance}, $\Gamma_{\{r,s\}}$ is also optimal. This implies that for any $k\geq 1$, every adjacent two-time marginal of  $\Gamma_{T^{(k)}}$ is optimal by the projective condition, so we can run the Minkowski argument in Step 1 again with $\Gamma^{(k)}$ replaced by $\Gamma_{T^{(k)}}$.

    Therefore, since the $T^{(k)}$-marginal of $\Gamma$ is $\Gamma_{T^{(k)}}$, for $\Gamma \otimes \Gamma$-almost every $(\mathbf{x}, \mathbf{x}')$, we have
    \begin{equation}
        \label{eqn:ae-lipschitz}
        d_Z(\tilde{\omega}_{i/2^k}(\mathbf{x},\mathbf{x}'), \tilde{\omega}_{(i+1)/2^k}(\mathbf{x},\mathbf{x}')) = \frac{1}{2^k}d_Z(\omega_X(x_0, x_0'), \omega_Y(x_1, x_1'))
    \end{equation}
    for all dyadic $i/2^k, (i+1)/2^k \in \mathbb{D}$. Chaining these equalities together, we obtain the desired equality \eqref{eqn:Lipschitz_condition_for_omega}. We note that the set of $(\mathbf{x},\mathbf{x}')$ on which \eqref{eqn:ae-lipschitz} could be different for each $i,k$, and although each set has measure $1$, an arbitrary intersection of full measure sets need not have full measure. However, we only need to take an intersection of a countable collection, in which case the intersection is of full measure, so the set on which \eqref{eqn:ae-lipschitz} is true for all $i,k$ has also measure $1$. Therefore the chaining operation is valid.

We will now extend our $\tilde{\omega}_t$ to all $t\in[0,1]$ by the Lipschitz extension argument earlier, but we will do so in a more careful way to make sure that the defined kernel is measurable and $p$-integrable. We first define $q_k:[0,1]\to  T^{(k)}$ as $q_k(t)=\lfloor t2^{k}\rfloor/2^k$, and then define $f_{k,t}:\Omega\times \Omega\to Z$ by 
\begin{equation}
f_{k,t}(\mathbf x,\mathbf x') = \begin{cases}\tilde{\omega}_{q_k(t)}(\mathbf{x}, \mathbf{x}'), & (\mathbf{x}, \mathbf{x}')\in B \\ z_0, \quad &\textrm{otherwise.} 
\end{cases}
\end{equation}
Here, $B$ is a set on which \eqref{eqn:ae-lipschitz} is true for all $i,k$, and $z_0\in Z$ is any fixed element. Since $B$ is measurable, $f_{k,t}$ is measurable for each $k,t$. Equation \eqref{eqn:ae-lipschitz} shows that for each $(\mathbf x, \mathbf x')\in B$ and $t\in [0,1]$, $k\mapsto f_{k,t}(\mathbf x, \mathbf x')$ is a Cauchy sequence in $Z$ because 
\begin{equation}
    d_Z(f_{n,t}, f_{m,t})=|q_{m}(t)-q_n(t)|d_Z(\omega_X,\omega_Y) \leq (|q_m(t)-t| + |q_n(t)-t|)d_Z(\omega_X,\omega_Y) \leq \left(\frac{1}{2^m} + \frac{1}{2^n}\right)d_Z(\omega_X,\omega_Y).
\end{equation}
(In the above, and in various other calculations below, we omit $(\mathbf x, \mathbf x')$ arguments for the sake of conciseness.) Since $Z$ is assumed to be Polish, $f_{k,t}$ converges, so define 
\begin{equation}
    \tilde{\omega}_{t}(\mathbf{x},\mathbf{x}') = \lim_{k\to \infty} f_{k,t}(\mathbf{x},\mathbf{x}').
\end{equation}
Being a limit of measurable functions, $\tilde{\omega}_t$ is measurable. Moreover, it defines a geodesic from $\omega_X(x_0,x_0')$ to $\omega_Y(x_1,x_1')$ for all $(\mathbf{x},\mathbf{x}')\in B$, as can be seen by taking the limit of
\begin{equation}
    d_Z(f_{k,t},f_{k,s}) = |q_k(s)-q_k(t)|d_Z(\omega_X,\omega_Y),
\end{equation}
which comes from \eqref{eqn:ae-lipschitz}. Finally, for any $t\in [0,1]$, $\tilde{\omega}_t$ is $p$-integrable because 
\begin{equation}
 d_Z(\tilde{\omega}_t,z_0) \leq d_Z(\tilde{\omega}_t,\omega_X)+d_Z(\omega_X,z_0)=td_Z(\omega_X,\omega_Y)+d_Z(\omega_X,z_0) \leq (1+t)d_Z(\omega_X,z_0)+d_Z(\omega_Y,z_0),
\end{equation}
and both $\omega_X$ and $\omega_Y$ are $p$-integrable by assumption, so the claim follows by the elementary estimate $(a+b)^p\leq 2^{p-1}(a^p+b^p)$.  Therefore, $(\Omega, \tilde{\omega}_t,\Gamma)$ defines a valid $Z$-network.
We abuse notation and continue to use $\tilde{\omega}_t$ for the extension.

Before moving on to the next step, we observe that $\widetilde{X}_t$ also defines a geodesic: for $0\leq s \leq t \leq 1$, 
\begin{align*}
2^p\dgwz(\widetilde{X}_s,\widetilde{X}_t)^p
&\leq \iint_{\Omega \times \Omega} d_Z\bigl(\tilde{\omega}_s(\mathbf{x},\mathbf{x}'),\tilde{\omega}_t(\mathbf{x},\mathbf{x}')\bigr)^p d(\Gamma \otimes \Gamma)(\mathbf{x},\mathbf{x}') \\
&= \iint_{\Omega \times \Omega} |t-s|^p\, d_Z\bigl(\omega_X(x_0,x_0'),\omega_Y(x_1,x_1')\bigr)^p d(\Gamma \otimes \Gamma)(\mathbf{x},\mathbf{x}') \\
&= 2^p|t-s|^p\,\dgwz(X,Y)^{p},
\end{align*}
where the last equality uses that $\Gamma$ marginalizes to an optimal coupling.

     \smallskip
    \noindent \textbf{Step 4. The auxiliary path  $\widetilde{X}_t$ is weakly isomorphic to the original geodesic $X_t$.}
    We will now show that $\widetilde{X}_t$ is weakly isomorphic to $X_t$ for any $t\in [0,1]$. First, consider the case where $t\in \mathbb{D}$, and consider the projection map $q_t:\Omega \to X_t, q_t(\mathbf{x})=x_t$. By the construction of $\Gamma$, we have $(q_t)_* \Gamma = \mu_t$, i.e., the map is measure-preserving. Moreover, for $\Gamma^{\otimes 2}-$almost every $\mathbf{x},\mathbf{x}'\in \Omega$, we have
    \begin{equation}
        \tilde{\omega}_t(\mathbf{x},\mathbf{x}') = \omega_t(x_t,x'_t) =  \omega_t(q_t(\mathbf{x}),q_t(\mathbf{x}'))=q_t^*\omega_t(\mathbf{x},\mathbf{x}').
    \end{equation}
    Therefore, $q_t$ is a weak isomorphism from $\widetilde{X}_t$ to $X_t$ (see Definition \ref{def:weak_iso}).

    It remains to consider the case where $t \in [0,1] \setminus \mathbb{D}$. This follows by a standard density argument and properties of the paths. Given $\varepsilon > 0$, choose $t' \in \mathbb{D}$ such that $|t' - t| < \varepsilon$. Then 
    \[
    \dgwz(\widetilde{X}_t,X_t) \leq \dgwz(\widetilde{X}_t,\widetilde{X}_{t'}) + \dgwz(\widetilde{X}_{t'},X_{t'}) + \dgwz(X_{t'},X_t) < 2\varepsilon \cdot \dgwz(X,Y),
    \]
    where we use the facts that $X_t$ and $\widetilde{X}_t$ are geodesics and that $\dgwz(\widetilde{X}_{t'},X_{t'})=0$ (since $t' \in \mathbb{D}$). Since this holds for all $\varepsilon > 0$, it must be that $\dgwz(\widetilde{X}_t,X_t)=0$, so that this part of the proof is complete.

    \smallskip
    \noindent \textbf{Step 5. The auxiliary path  $\widetilde{X}_t$ is weakly isomorphic to a generalized interpolation.} To complete the proof, we will construct a generalized interpolation $\widehat{X}_t = (X \times Y \times \mathbf{I}, \hat{\omega}_t, \pi\otimes \mathscr{L})$. We will then show that the auxiliary path $\widetilde{X}_t$ is weakly isomorphic to $\widehat{X}_t$. 
    
    Let $J = \{0,1\} \subset \mathbb{D}$ and set $\pi = (\mathrm{proj}_{\Omega,S_J})_\ast \Gamma$. We claim that $\pi$ is an optimal coupling of $X$ and $Y$. Indeed, by construction of $\Gamma$, we have that $\pi = \Gamma_J$. Moreover, $\Gamma_J$ is the limit of a subsequence of couplings $\Gamma_J^{(k_j)}$. In the notation of Step 1, $\Gamma_J^{(k_j)} = \pi^{(k_j)}$, which is an optimal coupling between $X$ and $Y$. Therefore
    \[
    \pi = \lim_{j \to \infty} \pi^{(k_j)},
    \]
    and continuity of the $Z$-Gromov-Wasserstein distortion functional implies that $\pi$ must be an optimal coupling.

    Set $\mathbb{D}' = \mathbb{D} \setminus \{0,1\}$. By the disintegration theorem~\cite[Theorem A.1]{BrediesFanzon2020}, there exists a family of probability measures $\{\pi_{x_0,x_1}\}_{x_0\in X, x_1\in Y}$ on $\prod_{t \in \mathbb{D}'} X_t$ satisfying 
    \[
    d\Gamma((x_t)_{t \in \mathbb{D}}) = d\pi_{x_0,x_1}((x_t)_{t \in \mathbb{D}'})d\pi(x_0,x_1).
    \]
    By \cite[Lemma 3.2 (vii)]{kallenberg2002foundations}, there exists a Borel function $G:X_0\times X_1\times \mathbf{I} \to \prod_{t\in \mathbb{D}'} X_t$ with $G(x_0,x_1, \cdot)_* \mathscr{L} = \pi_{x_0,x_1}$. We then define $H: X_0\times X_1\times \mathbf{I} \to \Omega$ by 
    \begin{equation}
        H(x_0,x_1,r)= \left(x_0,\left(G(x_0,x_1,r)_t \right)_{t \in \mathbb{D}'} ,x_1 \right).
    \end{equation}
    This allows us to define a generalized interpolation $\widehat{X}_t = (X \times Y \times \mathbf{I}, \hat{\omega}_t, \pi \otimes \mathscr{L})$, where $\pi$ is as defined above, and $\hat{\omega}_t$ is defined by
    \begin{equation}
        \widehat{\omega}_t((x_0,x_1,r), (x_0',x_1',r')) = \tilde\omega_t(H(x_0,x_1,r), H(x_0',x_1',r')).
    \end{equation}
    Then \eqref{eqn:Lipschitz_condition_for_omega} shows that $t \mapsto \widehat{\omega}_t((x_0,x_1,r), (x_0',x_1',r'))$ is a geodesic for any $((x_0,x_1,r),(x_0',x_1',r'))\in A$ with $A$ defined by 
    \begin{equation}
    A \coloneqq (H\times H)^{-1}(B), \quad H\times H((x_0,x_1,r),(x_0',x_1',r')) = (H(x_0,x_1,r),H(x_0',x_1',r'))
    \end{equation}
    where $B \subset \Omega \times \Omega$ is the set on which $t\mapsto \tilde{\omega}_t$ from Step 3 is a geodesic. Moreover, the set $A$ has measure $1$ because by \eqref{eqn:Lipschitz_condition_for_omega} we have $\Gamma^{\otimes2}(B)=1$, and as we will show in the following paragraph, we have $H_*(\pi\otimes \mathscr{L})=\Gamma$, so 
    \begin{equation}
        (\pi\otimes \mathscr{L})^{\otimes 2}(A)=(H_*(\pi\otimes \mathscr{L}))^{\otimes 2}(B)= \Gamma^{\otimes 2}(B) = 1.
    \end{equation}
    Therefore, $\widehat{X}_t$ is indeed a generalized interpolation.

    Finally, we will show that $H$ defines a weak isomorphism from $\widehat{X}_t$ to $\widetilde{X}_t$, for all $t \in [0,1]$. First, observe that $H^\ast \tilde{\omega}_t = \hat{\omega}_t$, by definition. It only remains to show that $H$ satisfies $H_\ast(\pi\otimes \mathscr{L})=\Gamma$. For a Borel set $A \subset \Omega$, let $1_A$ denote its indicator function, and consider 
    \begin{align}
        H_*(\pi\otimes \mathscr{L})(A) &= \int_{X_0\times X_1}\int_{[0,1]}1_{A}(x_0,G(x_0,x_1,t),x_1)d\mathscr{L}(t) d\pi(x_0,x_1)\\ & =  \int_{X_0\times X_1}\int_{\prod_{t\in \mathbb{D}'} X_t}1_{A}(x_0,s,x_1)d\pi_{x_0,x_1}(s) d\pi(x_0,x_1) \\ &=\int_{\Omega}1_A(x)  d\Gamma(x) \\ 
        &= \Gamma(A).
    \end{align}
    This shows that $H$ is a weak isomorphism, so we have finally proved that 
    \[
    [X_t] = [\widetilde{X}_t] = [\widehat{X}_t],
    \]
    where the latter space is a generalized interpolation. 
\end{proof}

\subsection{Proofs of Corollaries}\label{subsec:proofs_of_corollaries}

We now proceed with the proofs of the corollaries that were stated in Section \ref{sec:geodesics_main_results}, which primarily depend on the geodesic structure characterization for the $Z$-GW space (Point (2) of Theorem \ref{thm:geodesics_main} or Proposition \ref{prop:unique_geodesics}). 

\begin{proof}[Proof of Corollary \ref{cor:uniquely_geodesic}]
    Suppose that $Z$ is uniquely geodesic, let $p \in (1,\infty)$, and let $X$ and $Y$ be $(Z,p)$-networks. Let $X_t$ be a path of $(Z,p)$-networks which induces a geodesic $[X_t]$ between $[X]$ and $[Y]$ in $\MZp$. By Proposition \ref{prop:unique_geodesics}, we can assume without loss of generality that $X_t$ is a generalized interpolation, $X_t = (X \times Y \times \mathbf{I}, \omega_t, \pi \otimes \mathscr{L})$. By assumption, for $(\pi\otimes \mathscr{L})^{\otimes 2}$-a.e. $((x,y,s),(x',y',s'))\in X\times Y \times \mathbf{I}$, 
    \[
    t \mapsto \omega_t((x,y,s),(x',y',s'))
    \]
    is a geodesic between $\omega_X(x,x')$ and $\omega_Y(y,y')$, so it does not depend on $s$ or $s'$, by uniqueness of geodesics. Therefore, if we denote the geodesic between $\omega_X(x,x')$ and $\omega_Y(y,y')$ as $\gamma_{x,x',y,y'}:[0,1]\to Z$, then for any $t$ and for $(\pi\otimes \mathscr{L})^{\otimes 2}$-a.e. pair $((x,y,s),(x',y',s'))$, we have 
    \begin{equation}
        \gamma_{x,x',y,y'}(t) = \omega_t((x,y,s),(x',y',s')).
    \end{equation}
    
    We then define an interpolation (in the sense of Definition \ref{def:interpolation}) $\widetilde{X}_t = (X \times Y, \widetilde{\omega}_t, \pi)$, where 
    \[
    \widetilde{\omega}_t((x,y),(x',y')) = \gamma_{x,x',y,y'}(t).
    \]
    The measurability of $\tilde{\omega}_t$ follows from Aumann's measurable selection theorem \cite[Theorem B.9]{serieys2026nonlinear} similar to the proof of the first point of Theorem \ref{thm:geodesics_main}.
    Let $\psi:X \times Y \times \mathbf{I} \to X \times Y$ be the projection map. Then $\psi$ is measure-preserving, and, for any $t \in [0,1]$, 
    \[
    \psi^\ast \widetilde{\omega}_t((x,y,s),(x',y',s')) = \gamma_{x,x',y,y'}(t) = \omega_t((x,y,s),(x',y',s')), \quad (\pi\otimes \mathscr{L})^{\otimes 2}-a.e.
    \]
    Therefore $\psi$ is a weak isomorphism, so that $[\widetilde{X}_t] = [X_t]$ for all $t$.
\end{proof}

\begin{proof}[Proof of Corollary \ref{cor:geodesic_lifting}]
    Geodesic lifting for the edge law map follows easily by the fact that it admits an isometric section, as we proved in Theorem \ref{thm:submetries_for_GW}. 
    
    It remains to establish the geodesic lifting property for the quotient map $Q:\Lp(\mathbf{I}^2,Z) \to \MZp$. By Proposition \ref{prop:unique_geodesics}, any geodesic in $\MZp$ can be expressed as $[X_t]$, where $X_t = (X \times Y \times \mathbf{I}, \omega_t, \pi \otimes \mathscr{L})$ is a generalized interpolation. Following the construction in the first paragraph of Section \ref{sec:quotient_map}, choose a measurable map $\rho:\mathbf{I} \to X \times Y \times \mathbf{I}$ such that $\rho_\ast \mathscr{L} = \pi \otimes \mathscr{L}$ and set $\widetilde{\omega}_t = \rho^\ast \omega_t$. We claim that $t \mapsto \widetilde{\omega}_t$ is a geodesic in $\Lp(\mathbf{I}^2,Z)$. Indeed, for $0 \leq s \leq t \leq 1$, we have 
    \begin{align*}
        \dpz(\widetilde{\omega}_s,\widetilde{\omega}_t)^p &= \iint_{\mathbf{I}^2} d_Z\left(\widetilde{\omega}_s(u,v),\widetilde{\omega}_t(u,v)\right)^p \, d\mathscr{L}(u) d \mathscr{L}(v) \\
        &= \iint_{\mathbf{I}^2} d_Z(\omega_s(\rho(u),\rho(v)),\omega_t(\rho(u),\rho(v)))^p \, d\mathscr{L}(u) d \mathscr{L}(v) \\
        &= \iint_{(X \times Y \times \mathbf{I})^2} d_Z\left(\omega_s((x,y,r),(x',y',r')),\omega_t((x,y,r),(x',y',r'))\right)^p 
        \\
        &\hspace{2in} \cdot d(\pi \otimes \mathscr{L})(x,y,r) d (\pi \otimes \mathscr{L})(x',y',r') \\
        &\leq \iint_{(X \times Y \times \mathbf{I})^2} (t-s)^p d_Z(\omega_X(x,x'),\omega_Y(y,y'))^p \, d(\pi \otimes \mathscr{L})(x,y,r) d (\pi \otimes \mathscr{L})(x',y',r') \\
        &= \iint_{(X \times Y \times \mathbf{I})^2} (t-s)^p d_Z\left(\omega_0((x,y,r),(x',y',r')),\omega_1((x,y,r),(x',y',r'))\right)^p 
        \\
        &\hspace{2in} \cdot d(\pi \otimes \mathscr{L})(x,y,r) d (\pi \otimes \mathscr{L})(x',y',r') \\
        &= (t-s)^p \iint_{\mathbf{I}^2} d_Z(\omega_0(\rho(u),\rho(v)),\omega_1(\rho(u),\rho(v)))^p \, d\mathscr{L}(u) d \mathscr{L}(v) \\
        &= (t-s)^p \dpz(\widetilde{\omega}_0,\widetilde{\omega}_1)^p,
    \end{align*}
    which proves the claim. Applying the quotient map $Q$ to this geodesic gives 
    \[
    Q(\widetilde{\omega}_t) = [\mathbf{I},\widetilde{\omega}_t,\mathscr{L}] = [\mathbf{I}, \rho^\ast \omega_t, \mathscr{L}] = [X_t],
    \]
    so that $\widetilde{\omega}_t$ is a lift of $[X_t]$. 
\end{proof}

\subsection{Geodesics in the \texorpdfstring{$p=1$}{p=1} Case}\label{subsec:p_equals_1}

Finally, we prove that when $p=1$, all $Z$-spaces under consideration are geodesic (even if $Z$ is not). 

\begin{proof}[Proof of Theorem \ref{thm:geodesics_main}, Point (3)]
    First, we consider the space $\mathcal{L}_1(M,Z)$. Since the property of being geodesic is invariant under isometries, Proposition \ref{prop:nonlinear_Lebesgue_isometry} implies that we can reduce the problem to the convenient case of $\mathcal{L}_1(\mathbf{I},Z)$, where $\mathbf{I}$ is endowed with Lebesgue measure $\mathscr{L}$. 
    Let $\omega_0,\omega_1 \in \mathcal{L}_1(\mathbf{I},Z)$, and assume $\dz{1}(\omega_0,\omega_1) > 0$ (as otherwise the constant path defines a geodesic). We construct a geodesic between $\omega_0$ and $\omega_1$ of the form
    \begin{equation}
        \tilde{\omega}_t(u) = \begin{cases} \omega_1(u), & u\in A_t\\ \omega_0(u), & u\notin A_t,\end{cases}
    \end{equation}
    where $\{A_t\}_{t \in [0,1]}$ is a family of sets that still remains to be determined. We will now discover what properties are desirable for $A_t$, in order to make $\tilde{\omega}_t$ a geodesic. It is natural to impose $A_s\subset A_t$ for  $s<t$, in which case we obtain the expression
    \begin{equation}
        \dz{1}(\tilde{\omega}_s,\tilde{\omega}_t) = \int_{A_t}d_Z(\omega_0,\omega_1)d\mathscr{L} - \int_{A_s}d_Z(\omega_0,\omega_1)d\mathscr{L}.
    \end{equation}
    We want the right hand side to be equal to $(t-s)\dz{1}(\omega_0,\omega_1)$, which will be achieved if $A_t$ satisfies
    \begin{equation}\label{eqn:geodesic_1_condition}
    \int_{A_t}d_Z(\omega_0,\omega_1)du = t\dz{1}(\omega_0,\omega_1).
    \end{equation}
    To enforce the nested property $A_s\subset A_t$, we define $A_t$ to be a sublevel set $A_t = \{u\in \mathbf{I}|F(u)< t\}$ for some function $F: \mathbf{I}\to\mathbf{I}$, to be determined. In light of the desired equality \eqref{eqn:geodesic_1_condition}, this function $F$ should satisfy
    \begin{equation}
    \frac{1}{\dz{1}(\omega_0,\omega_1)}\int_{\{F(u)<t\}}d_Z(\omega_0,\omega_1)du = t,
    \end{equation}
    which implies that the distribution of $F(U)$ where $U$ is sampled from the probability distribution on $[0,1]$ with density $\frac{d_Z(\omega_0,\omega_1)}{\dz{1}(\omega_0,\omega_1)}$, is uniform. It is a standard fact that the CDF of $U$ satisfies such a property, so we should take
    \begin{equation}
        F(u) = \int_{0}^{u}\frac{d_Z(\omega_0,\omega_1)}{\dz{1}(\omega_0,\omega_1)} d\mathscr{L}.
    \end{equation}
    By construction, \eqref{eqn:geodesic_1_condition} holds, which implies the desired geodesic property $\dz{1}(\tilde{\omega}_s,\tilde{\omega}_t) = (t-s)\dz{1}(\omega_0,\omega_1)$. Moreover, 
    \[
    A_0=\{u \in \mathbf{I} \mid F(u) < 0\} = \emptyset,
    \]
    which implies that $\tilde{\omega}_0$ is equal to $\omega_0$. It remains to show that the right endpoint of the geodesic $\tilde{\omega}_1$ agrees with $\omega_1$ almost everywhere. The Lebesgue distance between $\tilde{\omega}_1$ and $\omega_1$ is given by 
    \begin{align*}
        \dz{1}(\tilde{\omega}_1,\omega_1) &= \int_\mathbf{I} d_Z(\tilde{\omega}_1(u),\omega_1(u)) d\mathscr{L}(u) \\
        &=\int_{\mathbf{I} \setminus A_1} d_Z(\omega_0(u),\omega_1(u)) d\mathscr{L}(u) 
        = \int_{\{u \in \mathbf{I} \mid F(u) = 1\}} d_Z(\omega_0(u),\omega_1(u)) d\mathscr{L}(u).
    \end{align*}
    For any $u \in \mathbf{I}$ with $F(u) = 1$, we have 
    \[
    0 = 1- F(u) = \int_u^1 \frac{d_Z(\omega_0,\omega_1)}{\dz{1}(\omega_0,\omega_1)} d\mathscr{L},
    \]
    so that $d_Z(\omega_0(v),\omega_1(v)) = 0$ for $\mathscr{L}$-a.e.\, $v \geq u$. Therefore, $\dz{1}(\tilde{\omega}_1,\omega_1) = 0$, and equality of the right endpoint of the geodesic with $\omega_1$ almost everywhere then follows.
    This completes the proof for the case of $\mathcal{L}_1(\mathbf{I},Z)$, hence for any $\mathcal{L}_1(M,Z)$.

    Geodesicity of the remaining spaces $\bigl(\MZ{1},\dgw{1}\bigr)$ and $\bigl(\WZp{1},\wzp{1}\bigr)$ now follows by general principles. As in the proof of Theorem \ref{thm:geodesics_main}, Point (1) (see Section \ref{subsec:geodesic_equivalence}), the quotient map and edge law submetries satisfy the conditions of Lemma \ref{lem:submetries_and_geodesics}. Since $\Lpp{1}(\mathbf{I}^2,Z)$ is geodesic, it follows that the $Z$-GW and $Z$-Wasserstein spaces are as well.
\end{proof}

\section{Alexandrov Curvature Bounds}\label{sec:curvature}

This section concerns the Alexandrov curvature of $Z$-GW spaces, as it relates to curvature of related spaces, such as $Z$ and the $Z$-Wasserstein space. Mirroring the structures of the earlier sections, we begin by recalling some concepts and summarizing our main results below, with the rest of the section being devoted to proving the theorem piecewise.

\subsection{Main Result: Curvature Bounds for \texorpdfstring{$Z$}{Z}-Spaces}

We begin by recalling some definitions and setting notation.

\begin{definition}[Curvature bounds, \cite{alexander2024alexandrov}]
    Let $(M,d_M)$ be a complete length space and let $K \in \R$. Let $\mathbb{M}_K$ denote the simply connected 2-dimensional space of constant curvature $K$; we denote its metric by $d_K$ and its diameter by $\varpi_K$, namely,
    \begin{equation}\label{eqn:diameter_model_space}
        \varpi_K = \left\{\begin{array}{rl}
        +\infty & \mbox{if $K \leq 0$,}\\
        \pi/\sqrt{K} & \mbox{if $K > 0$.}
        \end{array}\right.
    \end{equation} 
    Given three points $z, x_0, x_1 \in M$, a \textbf{$K$-comparison triangle} is a triple of points $\bar{z}$, $\bar{x_0}$, $\bar{x_1}$ in $\mathbb{M}_K$ with side lengths $d_K(\bar{z},\bar{x_0}) = d_M(z, x_0)$, $d_K(\bar{z},\bar{x_1}) = d_M(z, x_1)$, and $d_K(\bar{x}_0,\bar{x}_1) = d_M(x_0, x_1)$; a comparison triangle exists when $d_M(z,x_0) + d_M(x_0,x_1) + d_M(x_1,z) < 2 \varpi_K$. If a comparison triangle exists, $\angk{z}{x_0, x_1}$ denotes the model angle at $\bar{z}$.

    We say that $M$ has \textbf{curvature bounded below by $K$}, or that $M$ is $\mathrm{CBB}(K)$, if for any four points $z,x_0,x_1,x_2 \in M$, it holds that 
\begin{equation}
    \angk{z}{x_0, x_1} + \angk{z}{x_1, x_2} + \angk{z}{x_2, x_0} \leq 2\pi,
\end{equation}
    or if one of the comparison triangles does not exist. 
    Following \cite[Section 2.2]{Sturm2006geometryI}, we introduce the notation 
    \begin{equation}
        \curv(M, d) = \sup\{K\in \mathbb{R}| \mbox{ $M$ has curvature bounded below by $K$}\}.
    \end{equation}

    We say that $M$ has \textbf{curvature bounded above by $K$}, or that $M$ is $\mathrm{CAT}(K)$ (using the traditional acronym, in honor of Cartan-Aleksandrov-Toponogov), if for any four points $z_0, z_1, x_0, x_1 \in M$,
    \[
    \angk{z_0}{x_0, x_1} \leq \angk{z_0}{z_1, x_1} + \angk{z_0}{x_0, z_1} \quad \mbox{or} \quad \angk{z_1}{x_0, x_1} \leq \angk{z_1}{z_0, x_1} + \angk{z_1}{x_0, z_0},
    \]
    or if one of the comparison triangles does not exist. 
    We introduce the notation 
    \begin{equation}
        \curvupper(M, d) = \inf\{K\in \mathbb{R}| \mbox{ $M$ has curvature bounded above by $K$}\}.
    \end{equation}
\end{definition}

We record some basic facts about curvature, which are immediate from the definitions and will be used frequently throughout this section:
\begin{enumerate}[label=(\arabic*)]
    \item If $(M,d)$ and $(M',d')$ are isometric, then they share the same curvature bounds.
    \item If $(M',d')$ is a complete length space which isometrically embeds into $(M,d)$, then
    \[
    \curv(M,d) \leq \curv(M',d') \quad \mbox{and} \quad \curvupper(M,d) \geq \curvupper(M',d').
    \]
    \item For any $K\in \mathbb{R}$, $\curvupper(M,d)=K \Rightarrow (M,d)$ is $\mathrm{CAT}(K)$, and $\curv(M,d)=K \Rightarrow (M,d)$ is $\mathrm{CBB}(K)$ \cite[Proposition 8.4 and 9.7]{alexander2024alexandrov}. 
\end{enumerate}

We note that Theorem \ref{thm:geodesics_main} implies that when $Z$ is geodesic,   each of the associated spaces $(\Lp(\mathbf{I}^2,Z),\dpz)$, 
$(\MZp, \dgwz)$ and 
$(\WZ,\wz)$ is also geodesic, and it therefore makes sense to consider curvature bounds on the $Z$-spaces. The following theorem completely characterizes Alexandrov curvature bounds on the $Z$-spaces, relative to curvature bounds on $Z$ itself. In the statement, we assume that the space $Z$ has at least two points; otherwise, all spaces under consideration are isomorphic to the one-point metric space.

\begin{theorem}[Curvature Properties for $Z$-Spaces]\label{thm:curvature}
    Let $(M,\mu)$ be a Polish, nonatomic probability space and let $(Z,d_Z)$ be a geodesic Polish space containing at least two points.
    \begin{enumerate}[label=(\arabic*)]
        \item (Existence of Lower Curvature Bounds, $p=2$) The following are equivalent:
        \begin{enumerate}
            \item $\curv(Z,d_Z) \geq 0$
            \item $\curv(\Lpp{2}(M,Z),\dz{2}) = 0$
            \item $\curv(\MZ{2}, \dgw{2}) = 0$
            \item $\curv(\WZp{2},\wzp{2}) = 0$.
        \end{enumerate}
        \item (Nonexistence of Lower Curvature Bounds, $p=2$) The following are equivalent:
        \begin{enumerate}
            \item $\curv(Z,d_Z) < 0$
            \item $\curv(\Lpp{2}(M,Z),\dz{2}) = -\infty$
            \item $\curv(\MZ{2}, \dgw{2}) = -\infty$
            \item $\curv(\WZp{2},\wzp{2}) = -\infty$.
        \end{enumerate}
        \item (Nonexistence of Lower Curvature Bounds, $p \neq 2$) For $p \in [1,2) \cup (2,\infty)$,
        \[
        \curv(\Lp(M,Z),\dpz) = \curv(\MZp, \dgwz) = \curv(\WZ,\wz) = -\infty.
        \]
        \item (Upper Curvature Bounds, $p=2$) The following characterizations of upper curvature bounds hold in the $p=2$ case:
        \begin{enumerate}
            \item $\displaystyle \curvupper(\Lpp{2}(M,Z),\dz{2}) = \left\{\begin{array}{rl}
            0 & \mbox{if $\curvupper(Z,d_Z) \leq 0$}  \\
            +\infty & \mbox{otherwise;}
            \end{array}\right.$
            
            \item  $\curvupper(\MZ{2}, \dgw{2}) = + \infty$, for any $Z$;
            
            \item $\displaystyle \curvupper(\WZp{2},\wzp{2}) = \left\{\begin{array}{rl}
            0 & \mbox{if $Z$ is isometric to an interval}  \\
            +\infty & \mbox{otherwise.}
            \end{array}\right.$ 
        \end{enumerate}
        \item (Nonexistence of Upper Curvature Bounds, $p \neq 2$) For $p \in [1,2) \cup (2,\infty)$,
        \[
        \curvupper(\Lp(M,Z),\dpz) = \curvupper(\MZp, \dgwz) = \curvupper(\WZ,\wz) = +\infty.
        \]
    \end{enumerate}
\end{theorem}

\begin{remark}[Related Work]\label{remark:related_work_curvature}
    We give some context for these results below.

\begin{itemize}
    \item[(1),(2)] These connections between the lower curvature bounds for $Z$ and for the Wasserstein space $\WZ$ were first shown by Sturm~\cite[Proposition 2.10]{Sturm2006geometryI}. We use the $(d) \Rightarrow (a)$ direction from that work in our proof, but our proof of the converse is new, and factors through the submetry theorem (Theorem \ref{thm:submetries_for_GW}). 
    
    Lower curvature bounds for nonlinear Lebesgue spaces were obtained in \cite[Theorem 4.15]{serieys2026nonlinear}. The results there apply to more general nonlinear Lebesgue spaces, but the conclusions agree with ours; our proofs are distinguished in that they utilize the submetry theorem.

    Lower curvature bounds for the $Z$-GW space generalize the result of Sturm that the classical Gromov-Wasserstein space is nonnegatively curved~\cite{sturm2020}. The latter result was  mildly generalized to measure networks in~\cite{chowdhury2020gromov} and to Fused Gromov-Wasserstein distances with node features valued in a Hilbert space in~\cite{zhang2024geometry}; our theorem  subsumes all of these existing results. The lower curvature bounds are especially interesting, as they imply that the $Z$-GW space is an Alexandrov space of non-negative curvature, and therefore admits Riemannian manifold-like structures such as tangent cones, exponential maps, and gradients (see \cite{plaut2001metric,sturm2020}). This property was used by Sturm in the classical GW setting to describe gradient flows of certain functionals on the GW space~\cite[Chapters 7 and 8]{sturm2020}. Similar structures for the GW space of measure networks and for the Fused GW distance with features valued in a Hilbert space were utilized in~\cite{chowdhury2020gromov} and~\cite{zhang2024geometry}, respectively, to compute Fr\'{e}chet means of ensembles of (attributed) networks. A simple corollary of Theorem~\ref{thm:curvature} gives a generalization of the result for Fused GW distances.

    \item[(4)] The statement that if $Z$ is $\mathrm{CAT}(0)$ then so is $\mathcal{L}_2(M,Z)$ has appeared in various places in the literature~\cite{jost1994equilibrium,jost2012nonpositive,sturm2001nonlinear}. As was pointed out in~\cite[Remark 4.12]{serieys2026nonlinear}, there was an apparent gap in the existing proofs, which was resolved in~\cite[Proposition 4.11]{serieys2026nonlinear}. The converse statement was proved in~\cite[Proposition 4.13]{serieys2026nonlinear}. Our result is more detailed, in the sense that it considers arbitrary upper curvature bounds, rather than only the $\mathrm{CAT}(0)$ condition, and uses different proof techniques.

    Upper curvature bounds for the Wasserstein space are likely already known in the community, but we were unable to find a precise statement of our exact result in the literature. We give a proof for completeness, and explain below its connection to work in ~\cite{kloeckner2010geometric,bertrand2012geometric}.

    To our knowledge, the lack of upper curvature bounds for $Z$-GW spaces has not appeared previously in the literature, even in the classical case of $Z=\R$. A similar result for a metric on the space of unlabeled graphs appears in \cite[Theorem 1]{calissano2024populations}, but this is only related in spirit, as the geometries under consideration have no formal connection.

    \item[(3),(5)] The lack of curvature bounds for $p \neq 2$ is intuitive, due to similar phenomena for classical $\ell^p$-spaces (i.e., $\R^n$ endowed with the norm $\|\cdot\|_p$). However, we were surprised that we were unable to find formal statements of these results in the literature, even in the Wasserstein setting. Our proofs exploit the intuition coming from $\ell^p$ spaces by constructing (infinitesimal) embeddings of them into our $Z$-spaces, which yield obstructions to curvature bounds.   
\end{itemize}
\end{remark}

The rest of this section is devoted to proving Theorem \ref{thm:curvature}. The lower curvature bound existence statements from Point (1) are handled in Section \ref{sec:lower_curvature}. The lack of lower curvature bound statements from Points (2) and (3) are treated in Section \ref{sec:no_lower_curvature}. Upper curvature bounds in the $p=2$ case are considered in Section \ref{sec:upper_curvature}, and the lack of upper curvature bounds when $p \neq 2$ is established in Section \ref{sec:upper_curvature_p_not_2}.

\subsection{Lower Curvature Bounds: Existence}\label{sec:lower_curvature}

We first prove Point (1) of Theorem \ref{thm:curvature}, on the existence of lower curvature bounds for $Z$-spaces. The proof uses the following lemmas.

\begin{lemma}[{\cite[Corollary 8.34]{alexander2024alexandrov}}]\label{lem:submetries_preserve_curvature_bounds}
    Submetries preserve lower curvature bounds. That is, if $(M,d_M)$ is a complete length space which is $\mathrm{CBB}(K)$, $(N,d_N)$ is a metric space, and $f:M \to N$ is a submetry, then $(N,d_N)$ is a complete length space which is $\mathrm{CBB}(K)$.
\end{lemma}

The statements of the remaining lemmas assume that $(Z,d_Z)$ is a fixed geodesic Polish space containing at least two points.

The following result on curvature bounds for rescaled metrics is standard, but we include a proof sketch for completeness. It will be used in various places to deal with scaling factors of $1/2$ in the submetries of Theorem \ref{thm:submetries_for_GW}, as well as in certain constructions in the ensuing subsections. For $\epsilon > 0$, let $\epsilon d_Z$ denote the rescaling of the metric of $Z$ by a factor of $\epsilon$. 

\begin{lemma}\label{lem:scaling}
    If $(Z,d_Z)$ is $\mathrm{CBB}(K)$ (respectively, $\mathrm{CAT}(K)$), then $(Z,\epsilon d_Z)$ is $\mathrm{CBB}(K/\epsilon^2)$ (resp., $\mathrm{CAT}(K/\epsilon^2)$) for any $\epsilon > 0$.
\end{lemma}

\begin{proof}
    The model space of curvature $K/\epsilon^2$ is obtained by rescaling the metric of the model space $\mathbb{M}_K$ by a factor of $\epsilon$. The respective claims then follow easily by using the \emph{point-on-side comparison} characterizations of the $\mathrm{CBB}(K)$ condition~\cite[Theorem 8.14(b)]{alexander2024alexandrov}, or the $\mathrm{CAT}(K)$ condition~\cite[Theorem 9.14(b)]{alexander2024alexandrov}.
\end{proof}

\begin{lemma}\label{lem:L2_CBB}
    If $(Z,d_Z)$ is $\mathrm{CBB}(0)$ then so is $(\Lpp{2}(\mathbf{I}^2,Z),\dz{2})$. 
\end{lemma}

\begin{proof}
    As in \cite{sturm1999metric}, a metric space $(M, d_M)$ satisfies $\curv(M, d_M) \geq 0$ if and only if for any $x_0, x_1,\cdots x_n\in M$ and all $\lambda_1,\ldots, \lambda_n \in [0,\infty)$ with $\sum_{i=1}^n\lambda_i = 1$, we have
\begin{equation}\label{eqn:equivalent_curvature}
        \sum_{i,j=1}^n \lambda_i \lambda_j d_M(x_i, x_j)^2 \leq 2\sum_{i=1}^n \lambda_{i}d_M(x_0, x_i)^2.
    \end{equation}
    
    The proof of the lemma follows easily from the definition.  Indeed, let $\omega_0, \omega_1,\ldots,\omega_n \in \Lpp{2}(\mathbf{I}^2,Z)$ and $\lambda_1,\ldots,\lambda_n \in [0,\infty)$ with $\sum_{i}\lambda_i=1$. Then
    \begin{align}
        \sum_{i,j=1}^n \lambda_i \lambda_j \dz{2}(\omega_i,\omega_j)^2 &= \sum_{i,j=1}^n \lambda_i \lambda_j \int_{\mathbf{I}^2} d_Z(\omega_i(s,t),\omega_j(s,t))^2 ds \, dt \\
        &= \int_{\mathbf{I}^2} \left(\sum_{i,j=1}^n \lambda_i \lambda_j d_Z(\omega_i,\omega_j)^2 \right) \; ds \, dt \\
        &\leq \int_{\mathbf{I}^2} \left(2 \sum_{i=1}^n \lambda_i d_Z(\omega_0(s,t),\omega_i(s,t))^2\right) \; ds \, dt \\
        &= 2 \sum_{i=1}^n \lambda_i \int_{\mathbf{I}^2}   d_Z(\omega_0(s,t),\omega_i(s,t))^2 \; ds \, dt = 2 \sum_{i=1}^n \lambda_i \dz{2}(\omega_0,\omega_i)^2.
    \end{align}
\end{proof}

\begin{lemma}\label{lem:CBB_leq_0}
    For any $Z$, 
    \[
    \curv(\Lpp{2}(\mathbf{I}^2,Z),\dz{2}), \quad  \curv(\MZ{2},  \dgw{2}),  \quad \mbox{and} \quad\curv(\WZp{2},\wzp{2})
    \]
    are all non-positive.
\end{lemma}

\begin{proof}
The proof of \cite[Proposition 2.10 (iv)]{Sturm2006geometryI} shows that $\curv(\WZp{2},\wzp{2}) \leq 0$. If it were the case that $\curv(\MZ{2},  \dgw{2}) > 0$, then by definition, there exists $K>0$ so that $(\MZ{2},\dgw{2})$ is $\mathrm{CBB}(K)$. Then Theorem \ref{thm:submetries_for_GW} and Lemma \ref{lem:submetries_preserve_curvature_bounds} would imply that $(\WZp{2},\frac{1}{2}\wzp{2})$ is $\mathrm{CBB}(K)$, hence that $(\WZp{2},\wzp{2})$ is $\mathrm{CBB}(K/4)$ (by Lemma \ref{lem:scaling}), yielding a contradiction. The proof for $(\Lpp{2}(\mathbf{I}^2,Z),\dz{2})$ uses the same reasoning.
\end{proof}

We are now prepared to prove the first part of the theorem.

\begin{proof}[Proof of Theorem \ref{thm:curvature}, Point (1)]
    By Proposition \ref{prop:nonlinear_Lebesgue_isometry}, we can replace the generic probability space $(M,\mu)$ with $(\mathbf{I}^2,\mathscr{L} \otimes \mathscr{L})$ in the statement, which will be more convenient.

        Assume that $(a)$ holds. Then $(\Lpp{2}(\mathbf{I}^2,Z),\dz{2})$ is $\mathrm{CBB}(0)$, by Lemma \ref{lem:L2_CBB}. It follows from Lemma \ref{lem:CBB_leq_0} that $\curv(\Lpp{2}(\mathbf{I}^2,Z),\dz{2}) \leq 0$, and $(b)$ therefore follows. Assuming $(b)$, Theorem \ref{thm:submetries_for_GW} and Lemma \ref{lem:submetries_preserve_curvature_bounds} imply that  $(\MZ{2}, \dgw{2})$ is $\mathrm{CBB}(0)$, and Lemmas \ref{lem:CBB_leq_0} and \ref{lem:scaling} then imply that $(c)$ holds. By similar arguments, we get $(c) \Rightarrow (d)$. It is shown in \cite[Proposition 2.10 (iv)]{Sturm2006geometryI} that $(d) \Rightarrow (a)$, and this completes the proof.
\end{proof}

\subsection{Lower Curvature Bounds: Nonexistence}\label{sec:no_lower_curvature}

We now consider Points (2) and (3) of Theorem \ref{thm:curvature}, both dealing with the nonexistence of lower curvature bounds. The proof of the former is straightforward:

\begin{proof}[Proof of Theorem \ref{thm:curvature}, Point (2)]
    As in the proof of Point (1), we can once again assume that $(M,\mu) = (\mathbf{I}^2,\mathscr{L} \otimes \mathscr{L})$, without loss of generality. Since submetries preserve lower curvature bounds (Lemma \ref{lem:submetries_preserve_curvature_bounds}), and using the scaling lemma (Lemma \ref{lem:scaling}) to deal with factors of $1/2$ in the submetries, we have $(d) \Rightarrow (c) \Rightarrow (b)$. The implication $(b) \Rightarrow (a)$ follows from Lemma \ref{lem:L2_CBB}. Finally, $(a) \Rightarrow (d)$ is proved in \cite[Proposition 2.10 (iv)]{Sturm2006geometryI}.
\end{proof}

The only nontrivial part of Point (3) of Theorem \ref{thm:curvature} is the case of the $Z$-Wasserstein space, which we handle with the following lemma. The lemma uses the concept of an \emph{ultratangent space}, which we now recall, following~\cite[Chapter 13.A]{alexander2024alexandrov}. A \textbf{nonprincipal ultrafilter} is a finitely additive probability measure $\zeta$ on the natural numbers $\mathbb{N}$ taking only the values $0$ and $1$, such that $\zeta(S) = 0$ for any finite subset $S \subset \mathbb{N}$. 
We remark that \cite{alexander2024alexandrov} primarily uses a selective ultrafilter, which is a nonprincipal ultrafilter with an additional condition, but all results we use in this manuscript hold without this additional assumption. We make this distinction because the existence of a selective ultrafilter is not provable in ZFC (the continuum hypothesis is sufficient), while the existence of a nonprincipal ultrafilter follows from the axiom of choice \cite[p. 44] {alexander2024alexandrov}. We will discuss the selectivity assumption again once we have necessary terminology. 

Fixing a nonprincipal ultrafilter $\zeta$, one defines the \textbf{$\zeta$-limit} of a sequence of points $(m_n)_n$ in a metric space $(M,d_M)$ as $\lim_{n \to \zeta} m_n = m$ if and only if 
\[
\zeta\left(\{n \mid m_n \in  B_{d_M}(m,\epsilon)\}\right) = 1 \quad \forall \; \epsilon > 0.
\]
Now fix $m \in M$, and consider sequences $(m_n)_n$ such that the sequence of real numbers $(n \cdot d_M(m_n,m))_n$ is bounded. The set of such sequences may be endowed with a pseudometric $d^\zeta_m$ defined by 
\[
d^\zeta_m((m_n)_n,(m'_n)_n) \coloneqq \lim_{n \to \zeta} n \cdot d_M(m_n,m'_n)
\]
(the $\zeta$-limit being taken in $\R$ with its standard metric). The \textbf{ultratangent space} of $M$ at $m$ is the associated metric space, obtained by quotienting this set of sequences by the $d_m^\zeta = 0$ equivalence relation. 

We will now touch on the selectivity assumption employed in \cite{alexander2024alexandrov} again. The main result we use about the ultratangent space is the fact that a complete length $\mathrm{CBB}(K)$ (respectively, $\mathrm{CAT}(K)$) space has a $\mathrm{CBB}(0)$  (respectively, $\mathrm{CAT}(0)$) ultratangent space~\cite[Theorem 13.1]{alexander2024alexandrov}. The proof of this theorem relies on the continuity of the $\mathrm{CAT}(K)$ and $\mathrm{CBB}(K)$ conditions with respect to $K$
\cite[Proposition 8.4, 9.7]{alexander2024alexandrov}. This result follows from the definition of CAT/CBB and the fact that the usual continuity implies ultrafilter continuity, which follows directly from the definition of the ultrafilter limit. Therefore, although \cite{alexander2024alexandrov} fixes a selective ultrafilter, we do not need to do so.

\begin{lemma}
    \label{lem:Wasserstein_curvature_pneq2}
    For any Polish geodesic metric space $(Z, d_Z)$ that is not a singleton, and any $p \in [1,2) \cup (2,\infty)$, the $p$-Wasserstein space $(\WZ, \wz)$  is not $\mathrm{CBB}(K)$ for any $K\in \mathbb{R}$.
\end{lemma}

\begin{proof}
    Fix a nonprincipal ultrafilter $\zeta$. The idea of the proof is that the $p$-Wasserstein space locally looks like an $\ell^p$ space, i.e., the ultratangent space at some point contains a subset isometric to a rectangle in $\mathbb{R}^2$ with a scaled copy of the $\ell^p$ metric. This rectangle is not $\mathrm{CBB}(0)$ for $p\neq 2$, since a finite-dimensional normed space with a curvature lower bound must be an inner product space \cite[Proposition 1.8]{sturm1999metric}, and any four-point configuration violating the $\mathrm{CBB}$ condition can be translated and scaled into this rectangle. If the Wasserstein space had curvature at least $K$,
    then its ultratangent space would have curvature at least $0$
    \cite[Theorem 13.1(a)]{alexander2024alexandrov}, which would force the
    embedded complete geodesic rectangle (after shrinking it to a closed rectangle if needed for completeness) to be $\mathrm{CBB}(0)$, a
    contradiction.

    We proceed with the proof that an ultratangent space of $\WZ$ contains an isometric copy of a rectangle in $\R^2$, endowed with a scaled  $\ell^p$ distance. Let $a, b \in Z$ be distinct points and let $\gamma:[0, 1] \to Z$ be a geodesic from $a$ to $b$. For small enough $\epsilon>0$, define $\alpha, \beta : (-\epsilon, \epsilon) \to Z$ as $\alpha(s) = \gamma(s + \epsilon), \beta(s) = \gamma(1 - s - \epsilon)$, i.e., a small open geodesic interval around $\gamma(\epsilon)$ and $\gamma(1-\epsilon)$. Now, for any positive integer $n$, define $\mu_n: (-n\epsilon, n\epsilon)^2 \to \WZ$ as 
    \begin{equation}
        \mu_n(s, t) = \frac{1}{2}\delta_{\alpha(s/n)} + \frac{1}{2}\delta_{\beta(t/n)}.
    \end{equation} 
    Also set 
    \[
    \mu = \mu_n(0,0) = \frac{1}{2} \delta_{\gamma(\epsilon)} + \frac{1}{2} \delta_{\gamma(1-\epsilon)}.
    \]
    We will work with the ultratangent space to $\mu$. 

    We first claim that, for each $(s,t) \in (-\epsilon,\epsilon)^2$, the sequence of real numbers
    \[
    \big(n \cdot \wz(\mu_n(s,t),\mu)\big)_n
    \]
    is bounded. Indeed, it is a routine calculation to show that for any $s, s', t, t' \in (-n\epsilon, n\epsilon)$, we have
    \begin{equation}
        \label{eqn:Wasserstein_lp_distance}
        \wz(\mu_n(s, t), \mu_n(s', t')) = \frac{1}{n}d_Z(a, b) 2^{-1/p}\left(|s-s'|^p + |t - t'|^p\right)^{1/p}.
    \end{equation}
    Then, for any $n$, we have that
    \[
    n \cdot \wz(\mu_n(s,t),\mu) = d_Z(a, b) 2^{-1/p}\left(|s|^p + |t|^p\right)^{1/p},
    \]
    i.e., the sequence is constant. The sequences $(\mu_n(s,t))_n$ therefore represent elements of the ultratangent space at $\mu$. 

    We embed the rectangle $(-\epsilon,\epsilon)^2$ into the ultratangent space via the map taking $(s,t)$ to the $d_\mu^\zeta$-equivalence class of $(\mu_n(s,t))_n$. Next, we show that the ultratangent space pseudometric $d_\mu^\zeta$ recovers the $\ell^p$ distance, up to a constant scale factor. Indeed, invoking \eqref{eqn:Wasserstein_lp_distance}, we have  
    \[
    d_\mu^\zeta((\mu_n(s,t))_n,(\mu_n(s',t'))_n) = \lim_{n \to \zeta} n \wz(\mu_n(s,t),\mu_n(s',t')) = d_Z(a,b) 2^{-1/p} \|(s,t) - (s',t')\|_p.
    \]
    The lemma then follows by the argument described in the first paragraph of the proof.

\end{proof}

We now easily prove the third point of the main theorem.

\begin{proof}[Proof of Theorem \ref{thm:curvature}, Point (3)]
    The result follows immediately from the fact that lower curvature bounds are preserved under submetries (Lemma \ref{lem:submetries_preserve_curvature_bounds}), scaling properties of curvature (Lemma \ref{lem:scaling}; applied to deal with factors of $1/2$), the existence of submetries between the $Z$-spaces (Theorem~\ref{thm:submetries_for_GW}), and the lack of curvature bounds for $Z$-Wasserstein spaces (Lemma \ref{lem:Wasserstein_curvature_pneq2}).
\end{proof}

\subsection{Upper Curvature Bounds: The \texorpdfstring{$p=2$}{p=2} Case}\label{sec:upper_curvature}

This section is devoted to proving Point (4) of Theorem \ref{thm:curvature}, concerning upper curvature bounds. Upper bounds are not necessarily preserved by submetries, so the points here require more specialized argumentation than those of the lower curvature bounds in the previous subsections. The full result is broken up into three separate propositions. In the statements below, we assume that $(Z,d_Z)$ is a Polish geodesic space  containing more than one point and that $(M,\mu)$ is a nonatomic Polish probability space.

\begin{proposition}\label{prop:upper_curvature_Lp}
    If  $\curvupper \big(Z,d_Z\big) \leq 0$ then $\curvupper\big(\Lpp{2}(M,Z),\dz{2}\big) = 0$. Otherwise, $(\Lpp{2}(M,Z),\dz{2}\big)$ is not $\mathrm{CAT}(K)$ for any $K \in \R$.
\end{proposition}

The proof uses the following lemma.

\begin{lemma}\label{lem:Z_embedding}
    There exists an isometric embedding of $(Z,\epsilon d_Z)$ into $(\Lpp{2}(M,Z),\dz{2})$ for each $\epsilon \in (0,1]$.
\end{lemma}

\begin{proof}
    If $\epsilon = 1$, then it is easy to check that an isometric embedding is given by the map taking any point $z \in Z$ to the function taking constant value $z$. Let us then consider the case $\epsilon < 1$. Choose $A \subset M$ with $\mu(A) = \epsilon^2$, fix a point $z_0 \in Z$, and define $\phi:Z \to \Lpp{2}(M,Z)$ by
    \[
    \phi(z)(m) = \left\{
    \begin{array}{rl}
    z & m \in A \\
    z_0 & m \not \in A
    \end{array}\right.
    \]
    Then
    \[
    \dz{2}(\phi(z),\phi(z'))^2 = \int_A d_Z(z,z')^2 d\mu(m) + \int_{M \setminus A} d_Z(z_0,z_0)^2 d\mu(m) = \epsilon^2 d_Z(z,z')^2,
    \]
    so that $\phi$ isometrically embeds $\epsilon d_Z$ into $\dz{2}$. 
\end{proof}

\begin{proof}[Proof of Proposition \ref{prop:upper_curvature_Lp}]
    Suppose that $\curvupper \big(Z,d_Z\big) \leq 0$. Then $Z$ is, in particular, $\mathrm{CAT}(0)$, and it is shown in \cite[Proposition 1.2.18]{bacak2014convex} that this implies that $\Lpp{2}(M,Z)$ is $\mathrm{CAT}(0)$. On the other hand, we wish to show that $\Lpp{2}(M,Z)$ is not $\mathrm{CAT}(K)$ for any $K < 0$. This is done by embedding a copy of a flat square; to make the construction more concrete, let us assume without loss of generality that $M=\mathbf{I}$ with Lebesgue measure $\mathscr{L}$. Now, choose a non-constant geodesic path $\gamma:[0,1] \to Z$ in $Z$ and set $A = [0,1/2)$, $B = [1/2, 1]$. Define a map $\Gamma:\mathbf{I}^2 \to \Lpp{2}(\mathbf{I},Z)$ by 
    \[
    \Gamma(s,t)(u) = \left\{\begin{array}{lr}
    \gamma(s) & u \in A \\
    \gamma(t) & u \in B.
    \end{array}\right.
    \]
    Then $\Gamma$ is an isometric embedding, where we endow $\mathbf{I}^2$ with the usual Euclidean distance, rescaled by a factor of $\frac{d_Z(\gamma(0),\gamma(1))}{\sqrt{2}}$:
    \begin{align*}
    \dz{2}\big(\Gamma(s,t),\Gamma(s',t')\big)^2 &= \int_{\mathbf{I}} d_Z\big(\Gamma(s,t)(u),\Gamma(s',t')(u) \big)^2 d\mathscr{L}(u) \\
    &= \int_A d_Z\big(\Gamma(s,t)(u),\Gamma(s',t')(u) \big)^2 d\mathscr{L}(u) + \int_B d_Z\big(\Gamma(s,t)(u),\Gamma(s',t')(u) \big)^2 d\mathscr{L}(u) \\
    &= \int_A d_Z(\gamma(s),\gamma(s'))^2 d\mathscr{L}(u) + \int_B d_Z(\gamma(t),\gamma(t'))^2 d\mathscr{L}(u) \\
    &= \frac{d_Z(\gamma(0),\gamma(1))^2}{2} (|s-s'|^2 + |t-t'|^2).
    \end{align*}
    Since $\Lpp{2}(\mathbf{I},Z)$ contains a copy of a flat space, it cannot be $\mathrm{CAT}(K)$ for any $K < 0$. 

    Let us now prove the converse of the remaining statement; that is, we assume $\Lpp{2}(M,Z)$ is $\mathrm{CAT}(K)$ and prove that $\curvupper(Z,d_Z) \leq 0$. Since $Z$ isometrically embeds into $\Lpp{2}(M,Z)$, by Lemma \ref{lem:Z_embedding}, it must be that $\curvupper(Z,d_Z) = K' \leq K$. If $K' < 0$, we are done, so assume $K' \geq 0$. Then  $\curvupper(Z,\epsilon d_Z) = K'/\epsilon^2$, by Lemma \ref{lem:scaling}. By Lemma \ref{lem:Z_embedding}, $\epsilon d_Z$ also isometrically embeds in $\dz{2}$, which implies that $\curvupper(\Lpp{2}(M,Z),\dz{2}) \geq K'/\epsilon^2$. Taking $\epsilon$ sufficiently small yields a contradiction, unless $K'=0$, and this completes the proof.
\end{proof}

We proceed slightly out of order, and next prove the result pertaining to Wasserstein space.

\begin{proposition}\label{prop:upper_curvature_wasserstein}
    If $Z$ is isometric to an interval, then $\big(\WZp{2},\wzp{2})$ is $\mathrm{CAT}(0)$ and is not $\mathrm{CAT}(K)$ for any $K<0$. Otherwise, it is not $\mathrm{CAT}(K)$ for any $K \in \R$.  
\end{proposition}

This result is potentially known in the community. The fact that the Wasserstein space of the interval is $\mathrm{CAT}(0)$ follows easily from the explicit formula in terms of cumulative distribution functions for the Wasserstein distance~\cite[Chapter 2]{santambrogio2015optimal}---see \cite[Proposition 4.1]{kloeckner2010geometric} for a proof when $Z=\R$, which also applies to intervals. Moreover, it is explained in \cite[Remark 2.10]{bertrand2012geometric} that if $Z$ is $\mathrm{CAT}(0)$ and not isometric to an interval then $\WZp{2}$ is not $\mathrm{CAT}(0)$. In any case, we could not find a precise statement at the level of detail of this proposition in the literature.

The strategy of the proof is to  deal separately with the case that $Z$ is isometric to an interval, a circle, or neither of these. The interesting case turns out to be when $Z$ is isometric to a circle, in which case we wish to show that there exist arbitrarily close points in the Wasserstein space which admit distinct geodesics joining them, which suffices to prove the claim. To achieve this, it turns out to be easier to factor through the following construction.

Consider the space $Z^n/ S_n$, where the quotient is by the symmetric group $S_n$, acting by permuting coordinates. This group action is by isometries on the $\ell^2$-product metric on $Z^n$, so the product metric descends to a metric on the quotient space $Z^n/S_n$. For $p = (p_1,\ldots,p_n) \in Z^n$, let $[p] = [p_1,\ldots,p_n] \in Z^n/S_n$ denote its equivalence class under the action of $S_n$ (this is a potential overload of bracket notation for equivalence classes, but the meaning should be clear from context). Then the (rescaled) quotient metric is given explicitly by 
\begin{equation}\label{eqn:symmetric_product_metric}
d_{Z^n/S_n}([p],[q]) \coloneqq \frac{1}{\sqrt{n}} \min_{\sigma \in S_n} \left(\sum_{i=1}^n d_Z(p_i, q_{\sigma(i)})^2 \right)^{1/2},
\end{equation}
where we have scaled by a factor of $1/\sqrt{n}$, as this normalization will be convenient later on.

We now provide some lemmas regarding the metric space $(Z^n/S_n,d_{Z^n/S_n})$. The first two already exist in the literature.

\begin{lemma}[{\cite[Proposition~20]{needham2023geometric}
    and \cite[Lemma~4.7]{harms2023geometry}}]\label{lem:embedding_symmetric_product}
    For any $n \geq 1$, the metric space $(Z^n/S_n,  d_{Z^n/S_n})$ is a geodesic space which isometrically embeds into $(\WZp{2},\wzp{2})$ via the map
    \[
    [p_1,\ldots, p_n] \mapsto \frac{1}{n} \sum_{i=1}^n \delta_{p_i}.
    \]
\end{lemma}

\begin{lemma}[{\cite[Theorem 22]{needham2023geometric}}]\label{lem:not_isometric_to_interval}
    If $Z$ is not isometric to an interval or a circle (with geodesic distance), then $(Z^n/S_n,d_{Z^n/S_n})$ is not $\mathrm{CAT}(K)$ for any $K \in \R$ and any $n > 1$. 
\end{lemma}

Next, we refine the preceding lemma to provide additional curvature information for the circle.

\begin{lemma}\label{lem:symmetric_product_circle}
    Let $Z$ be isometric to a circle of radius $r > 0$, endowed with geodesic distance. Then $(Z^{2n}/S_{2n},d_{Z^{2n}/S_{2n}})$ is not $\mathrm{CAT}(K)$ for any $K < \frac{4n^2}{r^2}$. 
\end{lemma}

\begin{figure}
    \centering
    \includegraphics[width=0.65\linewidth]{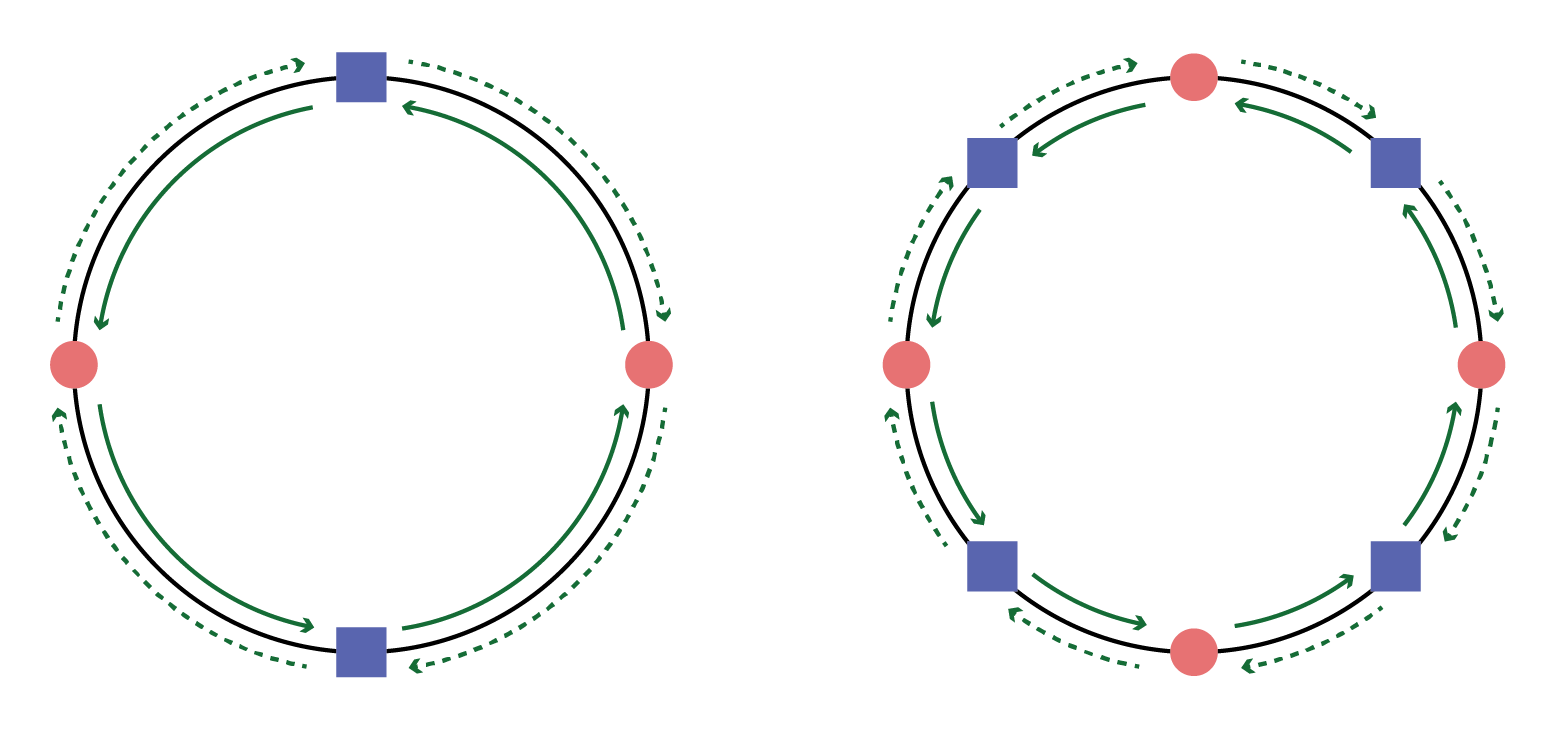}
    \caption{Examples of the points $p$ and $q$ from the proof of Lemma \ref{lem:symmetric_product_circle}, for $n=1$ on the left and $n=2$ on the right. In each case, the point $p$ is represented by $2n$ dots on the circle, and the point $q$ is represented by $2n$ squares. The solid arrows represent the identity matching between $p$ and $q$, and the dotted arrows represent the alternative optimal permutation. The geodesics between $[p]$ and $[q]$ are obtained by interpolating matched points along the arrows.}
    \label{fig:symmetric_product}
\end{figure}

\begin{proof}
        Assume that $Z$ is equal to a circle, which we  parameterize by $[0,2\pi)$ in the standard way---that is, we use $\theta \in [0,2\pi)$ to represent the point $r e^{i \theta} \in Z$, where $Z$ is considered as a subset of the complex plane. Consider the points in $Z^{2n}/S_{2n}$ defined (via the parameterization) by 
        \begin{align*}
        [p] &= \left[0, \frac{\pi}{n}, \frac{2\pi}{n}, \ldots, \frac{(2n-1)\pi}{n}\right] = \left[(k-1) \cdot \frac{\pi}{n} \right]_{k=1}^{2n} \\
        [q] &= \left[\frac{\pi}{2n},\frac{3\pi}{2n}, \frac{5\pi}{2n},\ldots, \frac{(4n-1) \pi}{2n}\right] = \left[(2k-1)\cdot \frac{\pi}{2n} \right]_{k=1}^{2n}
        \end{align*}
     (see Figure~\ref{fig:symmetric_product} for examples). The minimum in \eqref{eqn:symmetric_product_metric} is achieved by both the identity permutation and the cyclic shift permutation $1 \leftrightarrow 2n, 2 \leftrightarrow 1, 3 \leftrightarrow 2, \ldots, 2n \leftrightarrow 2n-1$. Indeed, while this claim is intuitively obvious, it can be proved by embedding into the Wasserstein space via Lemma~\ref{lem:embedding_symmetric_product} then using the explicit formula for Wasserstein distances between discrete distributions on the circle found in~\cite[Theorem 1]{rabin2011transportation}. Either of these optimal permutations then gives the distance between the points as 
     \begin{equation}\label{eqn:distance_between_points_symmetric_product}
     d_{Z^{2n}/S_{2n}}([p],[q]) = \frac{1}{\sqrt{2n}} \left(2n \left(\frac{\pi}{2n} r\right)^2\right)^{1/2} =  \frac{\pi r}{2n}
     \end{equation}
     An explicit formula for geodesics in $Z^{2n}/S_{2n}$ (for an arbitrary geodesic space $Z$) is given in \cite[Proposition 20]{needham2023geometric}; in particular, any optimal permutation $\sigma$ for $[p]$ and $[q]$ yields a geodesic which can be represented as a collection of geodesics in $Z$  between matched pairs $p_i$ and $q_{\sigma(i)}$. In the case of our example, this formula yields two distinct geodesics between $[p]$ and $[q]$, corresponding to the two optimal permutations---the first is induced by a counterclockwise rotation, and the second is induced by a clockwise rotation.

     Recall that $\varpi_K$ denotes the diameter of the model space $\mathbb{M}_K$. By \cite[Proposition 9.8]{alexander2024alexandrov}, since there exist distinct geodesics between $[p]$ and $[q]$, $Z^{2n}/S_{2n}$ is not $\mathrm{CAT}(K)$, provided $d_{Z^{2n}/S_{2n}}([p],[q]) < \varpi_K$. Comparing the distance calculation \eqref{eqn:distance_between_points_symmetric_product} with the formula for diameter \eqref{eqn:diameter_model_space} completes the proof.  

\end{proof}

\begin{proof}[Proof of Proposition \ref{prop:upper_curvature_wasserstein}]
    As stated above, the fact that $\WZp{2}$ is $\mathrm{CAT}(0)$ when $Z$ is isometric to an interval follows from (the proof of)~\cite[Proposition 4.1]{kloeckner2010geometric}. Now we will show that $\WZp{2}$ is not $\mathrm{CAT}(K)$ for any $K<0$, which does not exactly follow from \cite{kloeckner2010geometric}. Since $Z$ is isometric to an interval, we identify $Z$ with an interval $J\subset \mathbb{R}$. By Lemma \ref{lem:embedding_symmetric_product}, the metric space of unordered pairs $Z^{2}/S_{2}$ isometrically embeds into $\WZp{2}$. The space $Z^{2}/S_{2}$ is isometric, via the sorting map, to the closed convex region
    \[
        R \coloneqq \{j=(j_1,j_2)|j_1,j_2\in J,j_1\leq j_2\} \subset \R^2,
    \]
    endowed with the rescaled Euclidean metric
    \[
    d_R(j,j') \coloneqq 2^{-1/2}\sqrt{(j_1-j_1')^2+(j_2-j_2')^2}.
    \]
    Indeed, the sorting map is an isometry because the identity permutation is optimal for sorted pairs:
    \begin{equation}
        (j_1-j_1')^2+(j_2-j_2')^2-\{ (j_1-j_2')^2+(j_2-j_1')^2\} = -2(j_1j_1'+j_2j_2'-j_1j_2'-j_2j_1')=-2(j_1-j_2)(j_1'-j_2')\leq0.
    \end{equation}
    Now, $(R,d_R)$ is convex in the Euclidean plane and its straight line segments are geodesics for $d_R$. Hence, it contains a flat triangle, which is not possible in $\mathrm{CAT}(K)$ spaces for $K<0$.

    Now suppose that $Z$ is not isometric to an interval. If $Z$ is not isometric to a circle, then the isometric embedding $Z^n/S_n \to \WZp{2}$ from Lemma \ref{lem:embedding_symmetric_product}, together with the fact that $Z^n/S_n$ is not $\mathrm{CAT}(K)$ for any $K$ and $n > 1$ (Lemma \ref{lem:not_isometric_to_interval}) implies that $\WZp{2}$ is not $\mathrm{CAT}(K)$ for any $K$. It remains to handle the case that $Z$ is isometric to a circle. For each $n$, Lemma \ref{lem:embedding_symmetric_product} gives an isometric embedding $Z^{2n}/S_{2n} \to \WZp{2}$, so Lemma \ref{lem:symmetric_product_circle} implies that $\WZp{2}$ is not $\mathrm{CAT}(K)$ for any $K < \frac{4n^2}{r^2}$. Since $n$ is arbitrary, this completes the proof.
\end{proof}

Finally, we prove the Gromov-Wasserstein result.

\begin{proposition}\label{prop:upper_curvature_GW}
    For any $Z$, $\big(\MZ{2}, \dgw{2}\big)$ is not $\mathrm{CAT}(K)$ for any $K \in \R$. 
\end{proposition}

\begin{proof}
    By Proposition \ref{prop:upper_curvature_wasserstein} and the fact that Wasserstein space isometrically embeds in Gromov-Wasserstein space, the claim is already proved if $Z$ is not isometric to an interval. Let us prove the remaining case, and assume that $Z \subset \R$ is an interval. Let us assume for simplicity that $Z=\R$, as the general case follows easily by adapting the proof. We will show that there exist arbitrarily Gromov-Wasserstein-close $\R$-valued measure networks admitting distinct geodesics between them, which suffices to prove that the GW space is not $\mathrm{CAT}(K)$ for any $K$, by \cite[Proposition II.1.4]{bridsonhaefliger1999}. 

    The proof strategy described above will be carried out by specific example. Consider the measure networks $(X = \{x_1,x_2\}, \omega_X, \mu_X)$ and $(Y = \{y_1,y_2\}, \omega_Y, \mu_Y)$, where the measures are uniform and where the network functions are encoded by matrices (by ordering the points of the respective spaces by their indices), respectively, as follows:
    \[
    \omega_X = \begin{bmatrix}
        1 & \epsilon \\
        0 & 1
    \end{bmatrix}, \qquad \omega_Y = \begin{bmatrix}
        1+ \epsilon & 0 \\
        0 & 1-\epsilon
    \end{bmatrix},
    \]
    where we assume that $\epsilon > 0$ is sufficiently small (in the context of arguments to come later). Let us pause for a moment to describe the intuition behind this choice. The idea is to start with a measure network with a nontrivial symmetry; the simplest possible choice is represented by the $2 \times 2$ identity matrix. We then choose the simplest nonsymmetric perturbations of this starting point, leading to $\omega_X$ and $\omega_Y$. The goal is to show that there are two optimal couplings of $X$ and $Y$ which lead to distinct geodesics under the explicit formula for geodesics via interpolations. 

An arbitrary coupling $\pi$ of $\mu_X$ and $\mu_Y$ can be expressed in matrix form as 
\[
\pi = \begin{bmatrix}
\alpha & 1/2 - \alpha \\
1/2 - \alpha & \alpha
\end{bmatrix}, \qquad \alpha \in [0,1/2];
\]
that is, $\pi(x_1,y_1) = \alpha$, $\pi(x_1,y_2) = 1/2-\alpha$, and so on. For this small example, one can calculate the squared distortion $\mathrm{dis}_2^\R(\pi)^2$  (recall Definition \ref{def:zgw}) explicitly. In particular, it is quadratic in $\alpha$ with the coefficient of $\alpha^2$  given by $4 \epsilon - 8$. As we are taking $\epsilon$ to be small, this coefficient is negative, meaning that the minima of the distortion function can only occur at the endpoints of the interval of feasible $\alpha$'s; namely, when $\alpha = 0$ or $\alpha = 1/2$. The squared distortions at either of these values are $3\epsilon^2/4$, so that both couplings are optimal. These optimal couplings are induced by maps. That is, letting $\phi:X \to Y$ be the map $\phi:x_i \mapsto y_i$, $i =1,2$, $(\mathrm{id}_X \times \phi)_\ast \mu_X$ gives the coupling with $\alpha = 1/2$.  Similarly, letting $\psi:X \to Y$ denote the other possible bijection,  $(\mathrm{id}_X \times \psi)_\ast \mu_X$ gives the coupling with $\alpha = 0$.

The two optimal couplings described above lead to two geodesics given by interpolations, in the sense of Definition \ref{def:interpolation}. Since the optimal couplings involved in the interpolations are induced by bijective maps, the formula simplifies. Indeed, considered up to weak isomorphism, the geodesic induced by $\phi$ is equivalent to the path of measure networks $(X,\omega_t^\phi,\mu_X)$, with 
\[
\omega_t^\phi = \begin{bmatrix}
    1+ t \epsilon & (1-t)\epsilon \\
    0 & 1-t\epsilon 
\end{bmatrix},
\]
where we represent the network function as a matrix using the same indexing convention as above. Likewise, the interpolation induced by $\psi$ is weakly isomorphic to the simpler representation $(X,\omega_t^\psi,\mu_X)$, with 
\[
\omega_t^\psi = \begin{bmatrix}
    1-t \epsilon & (1-t)\epsilon \\
    0 & 1+t\epsilon
\end{bmatrix}.
\]
The midpoints (i.e., values at $t=1/2$) of the two interpolations are then represented by, respectively,
    \[
    \omega^\phi_{1/2} = \begin{bmatrix}
        1+\epsilon/2 & \epsilon/2 \\
        0 & 1-\epsilon/2
    \end{bmatrix}, \qquad \omega^\psi_{1/2} = \begin{bmatrix}
        1- \epsilon/2 & \epsilon/2 \\
        0 & 1+\epsilon/2
    \end{bmatrix}.
    \]
We claim that these midpoint networks are not weakly isomorphic, so that the induced geodesics joining $[X]$ and $[Y]$ are distinct, from which the desired result follows. To establish the claim, we once again show by brute force calculation that an optimal coupling between them must be induced by one of the two possible bijections from $X$ to itself (in this case, the quadratic term in the squared distortion has coefficient $-\epsilon^2/2 + 4 \epsilon - 8$, which is once again negative for small $\epsilon$). Here, the two bijections give different squared distortion values, with the smaller of the two (i.e., corresponding to the true optimal coupling) being equal to $\epsilon^2/8$. Since this minimal distortion is positive, the GW distance between the midpoint networks is positive, and this proves the claim.
\end{proof}

\subsection{Upper Curvature Bounds: The \texorpdfstring{$p \neq 2$}{p not 2} Case}\label{sec:upper_curvature_p_not_2}

Finally, we prove Point (5) of Theorem \ref{thm:curvature}. To proceed, we recall an alternative characterization of the $\mathrm{CAT}(K)$ condition. Let $\mathbb{M}_K$ denote the model space of constant curvature $K$, with metric denoted $d_K$ and diameter $\varpi_K$.  A space $(M,d_M)$ satisfies the \textbf{$\mathrm{CAT}(K)$ 4-point condition}~\cite[Definition II.1.10]{bridsonhaefliger1999} if for any quadruple $(m_1, m'_1, m_2, m'_2) \in M^4$ with 
\[
d_M(m_1, m'_1) + d_M(m'_1, m_2) + d_M(m_2, m'_2) + d_M(m'_2, m_1) < 2\varpi_{K},
\]
there exists a quadruple $(\bar{m}_1, \bar{m}'_1, \bar{m}_2, \bar{m}'_2) \in \mathbb{M}_{K}^4$ with $d_K(\bar{m}_i, \bar{m}'_j) = d_M(m_i, m'_j)$ for all $i, j \in \{1, 2\}$, such that $d_M(m_1, m_2) \leq d_K(\bar{m}_1, \bar{m}_2)$ and $d_M(m'_1, m'_2) \leq d_K(\bar{m}'_1, \bar{m}'_2)$. 

\begin{proposition}[{\cite[Proposition II.1.11]{bridsonhaefliger1999}}]
    Let $(M,d_M)$ be a complete metric space. The following are equivalent:
    \begin{enumerate}
        \item $M$ is $\mathrm{CAT}(K)$, and
        \item $M$ satisfies the $\mathrm{CAT}(K)$ 4-point condition and any pair of points $m,m' \in M$ with $d_M(m,m') < \varpi_K$ has  \textbf{approximate midpoints}: for all $\delta > 0$, there exists $m_\delta \in M$ such that 
        \[
        \max\{d_M(m,m_\delta),d_M(m',m_\delta)\} < \frac{1}{2} d_M(m,m') + \delta.
        \]
    \end{enumerate}
\end{proposition}

We utilize this characterization to prove the following.

\begin{lemma}\label{lem:isometric_embedding_CAT_K}
    Let $(M, d_M)$ and $(N, d_N)$ be metric spaces. Suppose $M$ is complete and let $N_0$ be a complete, geodesic subset of $N$ such that there exists an isometric embedding $N_0\hookrightarrow M$. If $M$ is $\mathrm{CAT}(K)$, then so is $N_0$.
\end{lemma}

\begin{proof}
        Since $M$ is complete and $\mathrm{CAT}(K)$, $M$ satisfies the $\mathrm{CAT}(K)$ 4-point condition. Given a quadruple of points in $N_0$, we can embed them in $M$ via the assumed isometric embedding, and it follows that $N_0$ also satisfies the $\mathrm{CAT}(K)$ 4-point condition. Moreover, since $N_0$ is geodesic, any two points have a midpoint (which, in particular, implies the existence of approximate midpoints), so $N_0$ is $\mathrm{CAT}(K)$.
\end{proof}

Next, we show that there exists an isometric embedding of a subset of $\MZp$ into $\Lp(\mathbf{I}^2,Z)$, which will allow us to apply Lemma \ref{lem:isometric_embedding_CAT_K}.

\begin{lemma}\label{lem:subset_embedding}
    Let $(Z,d_Z)$ be a separable geodesic metric space containing at least two points. There exists a complete, geodesic subspace $\mathcal{N}_0 \subset \MZp$ and an isometric embedding of  $(\mathcal{N}_0,\dgwz)$ into $(\Lp(\mathbf{I}^2,Z),\frac{1}{2}\dpz)$. Moreover, if \(p\neq 2\), then \(\mathcal{N}_0\) is not \(\mathrm{CAT}(K)\)
    for any \(K\in\mathbb{R}\).
\end{lemma}

\begin{proof}
    Let $a,b \in Z$ be distinct points and let $\gamma:[0,1]\to Z$ be a constant-speed geodesic from $a$ to $b$. Then define $\omega_{s, t}: \mathbf{I}^2\to Z$ for $s, t \in [0,\epsilon]$ for some $\epsilon \in (0,1)$, to be determined, as
    \begin{align}
        \omega_{s,t}(x, y) = \begin{cases}
            \gamma(s), & (x, y) \in A \coloneqq \mathbf{I}\times[0, 1/2],\\
            \gamma(1-t), & (x, y) \in B \coloneqq \mathbf{I}\times(1/2, 1].
        \end{cases}
    \end{align}
    Now define a subset $\mathcal{N}_0$ of $\MZp$ as
    \begin{equation}
        \mathcal{N}_0 = \{[\mathbf{I}, \omega_{s,t}, \mathscr{L}]|s, t \in [0, \epsilon]\}
    \end{equation}
    We claim that the map $\phi: \mathcal{N}_0 \hookrightarrow \Lp(\mathbf{I}^2, Z)$ defined as $\phi([\mathbf{I}, \omega_{s,t}, \mathscr{L}]) = \omega_{s,t}$ is a well-defined isometric embedding. Moreover, we claim that $\mathcal{N}_0$ is a complete, geodesic subset of $\MZp$.

    To prove these claims, we will show that, for any $(s,t),(s',t') \in [0,\epsilon]^2$, 
    \begin{equation}\label{eqn:embedding_distance_formula}
        2 \cdot \dgwz((\mathbf{I}, \omega_{s,t}, \mathscr{L}), (\mathbf{I}, \omega_{s', t'}, \mathscr{L})) = \dpz(\omega_{s,t}, \omega_{s', t'}) = 2^{-1/p} d_Z(a,b) \|(s,t) - (s', t')\|_p.
    \end{equation}
    From this, it follows that if $(s,t) \neq (s',t')$, then $[\mathbf{I},\omega_{s,t},\mathscr{L}] \neq [\mathbf{I}, \omega_{s',t'},\mathscr{L}]$, so that $\phi$ is well-defined, and that $\phi$ is an isometric embedding. Additionally, \eqref{eqn:embedding_distance_formula} proves that the map $(s,t) \mapsto [\mathbf{I},\omega_{s,t},\mathscr{L}]$ defines an isometric embedding of the rectangle $[0,\epsilon]^2 \subset \R^2$, endowed with the scaled $\ell^p$-norm $2^{-1-1/p} \cdot d_Z(a,b) \cdot \|\cdot \|_p$, onto $\mathcal{N}_0$, from which it follows that $\mathcal{N}_0$ is a complete, geodesic subspace.

    Let us now show  \eqref{eqn:embedding_distance_formula}. First, we calculate $\dpz(\omega_{s, t}, \omega_{s', t'})$ for any $s, s', t, t' \in [0, \epsilon]$:
    \begin{align}
        \dpz(\omega_{s,t}, \omega_{s', t'})^p &= \int_{\mathbf{I}^2}d_Z(\omega_{s, t}(x, y), \omega_{s', t'}(x, y))^p dx dy \\
        &= \int_{A}d_Z(\gamma(s), \gamma(s'))^p dx dy + \int_{B}d_Z(\gamma(1-t), \gamma(1-t'))^p dx dy \\
        &= \frac{1}{2} d_Z(a,b)^p |s-s'|^p + \frac{1}{2} d_Z(a,b)^p |t-t'|^p.
    \end{align}
    Therefore, we have 
    \begin{equation}
        \dpz(\omega_{s,t}, \omega_{s', t'}) = 2^{-1/p}d_Z(a,b)\left(|s-s'|^p + |t-t'|^p\right)^{1/p} = 2^{-1/p}d_Z(a,b)\|(s,t) - (s', t')\|_p.
    \end{equation}
     Next, we calculate the $Z$-GW distance between $X_{s,t} \coloneqq (\mathbf{I}, \omega_{s,t}, \mathscr{L})$ and $X_{s', t'} \coloneqq (\mathbf{I}, \omega_{s', t'}, \mathscr{L})$. We note that $\omega_{s,t}(x,y)$ is constant in $x$, so the $Z$-distortion of a coupling $\pi$ can be expressed as 
    \begin{align}
        \mathrm{dis}_p^Z(\pi)^p &= \int_{\mathbf{I}^2}\int_{\mathbf{I}^2}d_Z(\omega_{s,t}(x, y), \omega_{s', t'}(x', y'))^p d\pi(x, x') d\pi (y,y') \\
        &= \int_{\mathbf{I}^2}d_Z(\omega_{s,t}(x, y), \omega_{s', t'}(x', y'))^p d\pi (y,y') \\
        &= d_Z(\gamma(s), \gamma(s'))^p\cdot \pi_1 + d_Z(\gamma(s), \gamma(1-t'))^p\cdot \pi_2 \\   & \qquad + d_Z(\gamma(1-t), \gamma(1-t'))^p\cdot \pi_3 + d_Z(\gamma(1-t), \gamma(s'))^p \cdot \pi_4,
    \end{align}
    where
    \begin{align*}
        \pi_1 &\coloneqq \pi([0,1/2]\times [0,1/2]), \qquad \pi_2 \coloneqq \pi([0,1/2]\times (1/2,1]), \\ 
        \pi_3 &\coloneqq \pi((1/2,1]\times (1/2,1]), \qquad \pi_4 \coloneqq \pi((1/2,1]\times [0,1/2]).
    \end{align*}
    Since $\pi$ is a coupling, we have 
    \[
    \pi_1 + \pi_2 = \pi_3 + \pi_4 = \pi_1 + \pi_4 = \pi_2 + \pi_3 = 1/2.
    \]
    Solving these equations, we obtain $\pi_2 = \pi_4 = 1/2 - \pi_1$ and $\pi_3 = \pi_1$. Therefore, we have 
    \begin{equation}\label{eqn:final_distortion_expression}
    \begin{split}
        \mathrm{dis}_p^Z(\pi)^p &= \pi_1 \cdot \left(d_Z(\gamma(s), \gamma(s'))^p + d_Z(\gamma(1-t), \gamma(1-t'))^p\right)  \\  &\qquad + (1/2 - \pi_1) \cdot \left(d_Z(\gamma(s), \gamma(1-t'))^p + d_Z(\gamma(1-t), \gamma(s'))^p\right).
    \end{split}
    \end{equation}
    Thus, we have that the $Z$-distortion of $\pi$ is completely determined by the value of $\pi_1 = \pi([0,1/2]\times [0,1/2])$. It is straightforward to show that any coupling $\pi$ must have $\pi_1 \in [0,1/2]$, and conversely that each value in $[0,1/2]$ for $\pi_1$ can be obtained by some coupling $\pi$, so the $Z$-GW distance is attained at the minimum of \eqref{eqn:final_distortion_expression} over $\pi_1 \in [0, 1/2]$. Choosing $\epsilon$ to be sufficiently small ensures that 
    \[
    d_Z(\gamma(s), \gamma(s'))^p + d_Z(\gamma(1-t), \gamma(1-t'))^p < d_Z(\gamma(s), \gamma(1-t'))^p + d_Z(\gamma(1-t), \gamma(s'))^p,
    \]
    so that the minimum is achieved at $\pi_1 = 1/2$, hence
    \begin{align}
        2^p \cdot \dgwz(X_{s,t}, X_{s', t'})^p &= \frac{1}{2}\left(d_Z(\gamma(s), \gamma(s'))^p + d_Z(\gamma(1-t), \gamma(1-t'))^p\right) \\
        &= \frac{1}{2} d_Z(a,b)^p |s-s'|^p + \frac{1}{2} d_Z(a,b)^p |t-t'|^p,
    \end{align}
    which establishes  \eqref{eqn:embedding_distance_formula}. 
    
    It remains to prove the final assertion. By \eqref{eqn:embedding_distance_formula}, $\mathcal{N}_0$ is isometric to
    \begin{equation}
\bigl([0,\epsilon]^2,c\|\cdot\|_p\bigr),
    \quad
    c=2^{-1-1/p}d_Z(a,b).
    \end{equation}
    Therefore, we identify $\mathcal{N}_0$ with this space. 
    Choose an interior point \(q\in(0,\epsilon)^2\). We claim that the ultratangent space of this
    rectangle at \(q\) is isometric to $ \bigl(\mathbb{R}^2,c\|\cdot\|_p\bigr).$ Indeed, define a mapping $T$ from $(\mathbb{R}^2,c\|\cdot \|_p)$ to the ultratangent space by 
    \begin{equation}
    T(v) = \left[\left(q+\frac{v}{n}\right)_{n=1}^{\infty}\right]
    \end{equation}
    where the bracket $[\cdot]$ denotes the equivalence class in the ultratangent space; depending on $v$, we replace the first terms of the sequence by $q$, as needed, to make sure that it does not leave $[0,\epsilon]^2$. It is straightforward to show that $T$ preserves the metric, by checking the definition of the ultratangent metric, and it follows that $T$ is injective. Now we show that it is surjective. For any element $[(x_n)]$ of the ultratangent space, $nc\|x_n-q\|_{p}$ is bounded by definition, so $n(x_n-q)$ is a bounded sequence in $\mathbb{R}^2$. Therefore the ultralimit $v\coloneqq \lim_{n\to \zeta}n(x_n-q)$ is in $\mathbb{R}^2$, and satisfies
    \begin{equation}
        \lim_{n\to \zeta}nc\left\|x_n - \left(q+\frac{v}{n}\right)\right\|_p=\lim_{n\to \zeta}c\left\|n(x_n -q)-v\right\|_p=0.
    \end{equation}
    That is, $[(x_n)]$ and $T(v)$ have distance zero in the ultratangent metric, and they are equal as elements of the ultratangent space. Now, if \(\mathcal{N}_0\) were \(\mathrm{CAT}(K)\), then its ultratangent space would be
    \(\mathrm{CAT}(0)\) by \cite[Theorem 13.1]{alexander2024alexandrov}. This is impossible
    when \(p\neq2\), since a normed vector space is \(\mathrm{CAT}(0)\) only
    if its norm is induced by an inner product
    \cite[Proposition II.1.14]{bridsonhaefliger1999}.
\end{proof}

The final lemma shows that the $Z$-Wasserstein space is not $\mathrm{CAT}(K)$ for $p\neq 2$. 

\begin{lemma}
    \label{lem:Wasserstein_upper_curvature_pneq2}
    For any Polish geodesic metric space $(Z, d_Z)$ that is not a singleton, and any $p \in [1,2) \cup (2,\infty)$, the $p$-Wasserstein space $(\WZ, \wz)$  is not $\mathrm{CAT}(K)$ for any $K\in \mathbb{R}$.
\end{lemma}

\begin{proof}
    The proof is essentially the same as that of Lemma \ref{lem:Wasserstein_curvature_pneq2}. The ultratangent space of the Wasserstein space  at some point contains a subset isometric to a rectangle in $\mathbb{R}^2$ endowed with a scaled $\ell^p$ metric, by the construction given in the proof of Lemma \ref{lem:Wasserstein_curvature_pneq2}.
    The rectangle is not $\mathrm{CAT}(0)$ for $p\neq 2$, because a normed
    real vector space which is $\mathrm{CAT}(0)$ must be an inner product
    space~\cite[Proposition II.1.14]{bridsonhaefliger1999}, and any four-point
    configuration violating the $\mathrm{CAT}(0)$ condition can be translated
    and scaled into this rectangle. If the Wasserstein space were
    $\mathrm{CAT}(K)$, then its ultratangent space would be $\mathrm{CAT}(0)$
    by~\cite[Theorem 13.1]{alexander2024alexandrov}; Lemma~\ref{lem:isometric_embedding_CAT_K} would then
    imply that the embedded complete geodesic rectangle (after shrinking it to a closed rectangle to make it complete) is $\mathrm{CAT}(0)$,
    a contradiction.

\end{proof}

We conclude with the proof of the last point of Theorem \ref{thm:curvature}, which now follows easily.

\begin{proof}[Proof of Theorem \ref{thm:curvature}, Point (5)]
    By Lemmas \ref{lem:Wasserstein_upper_curvature_pneq2} and \ref{lem:scaling} (to handle the scaling factor), $(\WZ,\frac{1}{2}\cdot \wz)$ is not $\mathrm{CAT}(K)$ for any $K$. Since $\WZ$ isometrically embeds in $\MZp$, via the section of the edge law map from Theorem \ref{thm:submetries_for_GW}, it follows by Lemma \ref{lem:isometric_embedding_CAT_K} that $\MZp$ is not $\mathrm{CAT}(K)$. Finally, Lemmas \ref{lem:scaling}, \ref{lem:isometric_embedding_CAT_K} and \ref{lem:subset_embedding} imply that $\Lp(\mathbf{I}^2,Z)$ is not $\mathrm{CAT}(K)$, and Proposition \ref{prop:nonlinear_Lebesgue_isometry} shows that this conclusion applies to any $\Lp(M,Z)$.
\end{proof}

\printbibliography

@article{bloom1977counterexample,
  title={A counterexample to a theorem of S. Piccard},
  author={Bloom, Gary S},
  journal={Journal of Combinatorial Theory, Series A},
  volume={22},
  number={3},
  pages={378--379},
  year={1977},
  publisher={Elsevier}
}

@inproceedings{bottou2018geometrical,
  title={Geometrical insights for implicit generative modeling},
  author={Bottou, Leon and Arjovsky, Martin and Lopez-Paz, David and Oquab, Maxime},
  booktitle={Braverman Readings in Machine Learning. Key Ideas from Inception to Current State: International Conference Commemorating the 40th Anniversary of Emmanuil Braverman's Decease, Boston, MA, USA, April 28-30, 2017, Invited Talks},
  pages={229--268},
  year={2018},
  organization={Springer}
}

@article{Memoli2023characterization,
  title={Characterization of Gromov-type geodesics},
  author={M\'{e}moli, Facundo and Wan, Zhengchao},
  journal={Differential Geometry and its Applications},
  volume={88},
  pages={102006},
  year={2023},
  publisher={Elsevier}
}

@article{boutin2004reconstructing,
  title={On reconstructing n-point configurations from the distribution of distances or areas},
  author={Boutin, Mireille and Kemper, Gregor},
  journal={Advances in Applied Mathematics},
  volume={32},
  number={4},
  pages={709--735},
  year={2004},
  publisher={Elsevier}
}

@article{osada2002shape,
  title={Shape distributions},
  author={Osada, Robert and Funkhouser, Thomas and Chazelle, Bernard and Dobkin, David},
  journal={ACM Transactions on Graphics (TOG)},
  volume={21},
  number={4},
  pages={807--832},
  year={2002},
  publisher={ACM New York, NY, USA}
}

@article{schweizer1960statistical,
  title={Statistical metric spaces},
  author={Schweizer, Berthold and Sklar, Abe},
  journal={Pacific J. Math},
  volume={10},
  number={1},
  pages={313--334},
  year={1960}
}

@article{kramosil1975fuzzy,
  title={Fuzzy metrics and statistical metric spaces},
  author={Kramosil, Ivan and Mich{\'a}lek, Ji{\v{r}}{\'\i}},
  journal={Kybernetika},
  volume={11},
  number={5},
  pages={336--344},
  year={1975},
  publisher={Institute of Information Theory and Automation AS CR}
}

@article{wald1943statistical,
  title={On a statistical generalization of metric spaces},
  author={Wald, Abraham},
  journal={Proceedings of the National Academy of Sciences},
  volume={29},
  number={6},
  pages={196--197},
  year={1943}
}

@article{otto2001geometry,
  title={The geometry of dissipative evolution equations: The porous medium equation},
  author={Otto, F},
  journal={Communications in Partial Differential Equations},
  volume={26},
  number={1--2},
  pages={101--174},
  year={2001},
  publisher={Marcel Dekker Inc.}
}

@article{rabin2011transportation,
  title={Transportation distances on the circle},
  author={Rabin, Julien and Delon, Julie and Gousseau, Yann},
  journal={Journal of Mathematical Imaging and Vision},
  volume={41},
  number={1},
  pages={147--167},
  year={2011},
  publisher={Springer}
}

@article{harms2023geometry,
  title={Geometry of sample spaces},
  author={Harms, Philipp and Michor, Peter W and Pennec, Xavier and Sommer, Stefan},
  journal={Differential Geometry and its Applications},
  volume={90},
  pages={102029},
  year={2023},
  publisher={Elsevier}
}

@article{jain2009structure,
  author  = {Brijnesh J. Jain and Klaus Obermayer},
  title   = {Structure Spaces},
  journal = {Journal of Machine Learning Research},
  year    = {2009},
  volume  = {10},
  number  = {93},
  pages   = {2667--2714},
  url     = {http://jmlr.org/papers/v10/jain09a.html}
}

@article{menger1942statistical,
  author    = {Karl Menger},
  title     = {Statistical metrics},
  journal   = {Proceedings of the National Academy of Sciences of the United States of America},
  volume    = {28},
  number    = {12},
  pages     = {535--537},
  year      = {1942},
  publisher = {National Academy of Sciences},
  doi       = {10.1073/pnas.28.12.535}
}

@article{arya2026gromov,
  title={The Gromov--Wasserstein distance between spheres},
  author={Arya, Shreya and Auddy, Arnab and Clark, Ranthony A and Lim, Sunhyuk and M\'{e}moli, Facundo and Packer, Daniel},
  journal={Foundations of Computational Mathematics},
  volume={26},
  number={1},
  pages={75--130},
  year={2026},
  publisher={Springer}
}

@article{robertson2025generalization,
  title={On a generalization of Wasserstein distance and the Beckmann problem to connection graphs},
  author={Robertson, Sawyer Jack and Kohli, Dhruv and Mishne, Gal and Cloninger, Alexander},
  journal={SIAM Journal on Scientific Computing},
  volume={47},
  number={5},
  pages={A2774--A2800},
  year={2025},
  publisher={SIAM}
}

@article{bal2024statistical,
  title={Statistical analysis of complex shape graphs},
  author={Bal, Aditi Basu and Guo, Xiaoyang and Needham, Tom and Srivastava, Anuj},
  journal={IEEE Transactions on Pattern Analysis and Machine Intelligence},
  volume={46},
  number={12},
  pages={8788--8805},
  year={2024},
  publisher={IEEE}
}

@article{guo2022statistical,
  title={Statistical shape analysis of brain arterial networks ({BAN})},
  author={Guo, Xiaoyang and Bal, Aditi Basu and Needham, Tom and Srivastava, Anuj},
  journal={The Annals of Applied Statistics},
  volume={16},
  number={2},
  pages={1130--1150},
  year={2022},
  doi={10.1214/21-AOAS1536},
  publisher={Institute of Mathematical Statistics}
}

@article{memoli2011spectral,
  title={A spectral notion of Gromov--Wasserstein distance and related methods},
  author={M{\'e}moli, Facundo},
  journal={Applied and Computational Harmonic Analysis},
  volume={30},
  number={3},
  pages={363--401},
  year={2011},
  publisher={Elsevier}
}

@inproceedings{memoli2009spectral,
  title={Spectral Gromov-Wasserstein distances for shape matching},
  author={M{\'e}moli, Facundo},
  booktitle={2009 IEEE 12th International Conference on Computer Vision Workshops, ICCV Workshops},
  pages={256--263},
  year={2009},
  organization={IEEE}
}

@article{calissano2024populations,
  title={Populations of unlabelled networks: Graph space geometry and generalized geodesic principal components},
  author={Calissano, Anna and Feragen, Aasa and Vantini, Simone},
  journal={Biometrika},
  volume={111},
  number={1},
  pages={147--170},
  year={2024},
  publisher={Oxford University Press}
}

@article{serieys2026nonlinear,
  title={Nonlinear Lebesgue spaces: Curves and geometry},
  author={S{\'e}rieys, Guillaume},
  journal={arXiv preprint arXiv:2603.09459},
  year={2026}
}

@article{vandistributional,
  title={Distributional Reduction: Unifying Dimensionality Reduction and Clustering with Gromov-Wasserstein},
  author={Van Assel, Hugues and Vincent-Cuaz, C{\'e}dric and Courty, Nicolas and Flamary, R{\'e}mi and Frossard, Pascal and Vayer, Titouan},
  journal={Transactions on Machine Learning Research},
  year = {2025}
}

@inproceedings{clark2025generalized,
  title={Generalized dimension reduction using semi-relaxed Gromov-Wasserstein distance},
  author={Clark, Ranthony A and Needham, Tom and Weighill, Thomas},
  booktitle={Proceedings of the AAAI Conference on Artificial Intelligence},
  volume={39},
  number={15},
  pages={16082--16090},
  year={2025}
}

@article{demetci2022scotv2,
  title={{SCOTv2}: Single-cell multiomic alignment with disproportionate cell-type representation},
  author={Demetci, Pinar and Santorella, Rebecca and Chakravarthy, Manav and Sandstede, Bj{\"o}rn and Singh, Ritambhara},
  journal={Journal of Computational Biology},
  volume={29},
  number={11},
  pages={1213--1228},
  year={2022},
  publisher={SAGE Publications Sage CA: Los Angeles, CA}
}

@inproceedings{chowdhury2021quantized,
  title={Quantized gromov-wasserstein},
  author={Chowdhury, Samir and Miller, David and Needham, Tom},
  booktitle={Joint European Conference on Machine Learning and Knowledge Discovery in Databases},
  pages={811--827},
  year={2021},
  organization={Springer}
}

@article{solomon2016entropic,
  title={Entropic metric alignment for correspondence problems},
  author={Solomon, Justin and Peyr{\'e}, Gabriel and Kim, Vladimir G and Sra, Suvrit},
  journal={ACM Transactions on Graphics (ToG)},
  volume={35},
  number={4},
  pages={1--13},
  year={2016},
  publisher={ACM New York, NY, USA}
}

@article{chowdhury2018explicit,
  title={Explicit geodesics in Gromov-Hausdorff space},
  author={Chowdhury, Samir and M{\'e}moli, Facundo},
  journal={Electronic Research Announcements},
  volume={25},
  pages={48--59},
  year={2018},
  publisher={Electronic Research Announcements}
}

@article{sturm1999metric,
  title={Metric spaces of lower bounded curvature},
  author={Sturm, Karl-Theodor},
  journal={Expositiones Mathematicae},
  volume={17},
  pages={35--48},
  year={1999},
  publisher={SPEKTRUM ACADEMISCHER VERLAG}
}

@book{kechris2012classical,
  title={Classical descriptive set theory},
  author={Kechris, Alexander},
  year={1995},
  publisher={Springer Science \& Business Media}
}

@article{zhang2025topological,
  title={Topological optimal transport for geometric cycle matching},
  author={Zhang, Stephen Y and Stumpf, Michael PH and Needham, Tom and Barbensi, Agnese},
  journal={Journal of Applied and Computational Topology},
  volume={9},
  number={2},
  pages={11},
  year={2025},
  publisher={Springer}
}

@book{jost2012nonpositive,
  title={Nonpositive curvature: geometric and analytic aspects},
  author={Jost, J{\"u}rgen},
  year={1997},
  publisher={Birkh{\"a}user}
}

@article{sturm2001nonlinear,
  title={Nonlinear Markov operators associated with symmetric Markov kernels and energy minimizing maps between singular spaces},
  author={Sturm, Karl-Theodor},
  journal={Calculus of Variations and Partial Differential Equations},
  volume={12},
  number={4},
  pages={317--357},
  year={2001},
  publisher={Springer}
}

@book{shioya2016metric,
  title={Metric measure geometry},
  author={Shioya, Takashi},
  year={2016},
  publisher={European Mathematical Society-EMS-Publishing House GmbH}
}

@article{yang2024exploiting,
  title={Exploiting edge features in graph-based learning with fused network gromov-wasserstein distance},
  author={Yang, Junjie and Labeau, Matthieu and d'Alch{\'e}-Buc, Florence},
  journal={Transactions on Machine Learning Research},
  year={2024}
}

@article{kawano2024multi,
  title={Multi-Dimensional Fused Gromov Wasserstein Discrepancy for Edge-Attributed Graphs},
  author={Kawano, Keisuke and Koide, Satoshi and Shiokawa, Hiroaki and Amagasa, Toshiyuki},
  journal={IEICE TRANSACTIONS on Information and Systems},
  volume={107},
  number={5},
  pages={683--693},
  year={2024},
  publisher={The Institute of Electronics, Information and Communication Engineers}
}

@article{memoli2023ultrametric,
  title={The ultrametric Gromov--Wasserstein distance},
  author={M{\'e}moli, Facundo and Munk, Axel and Wan, Zhengchao and Weitkamp, Christoph},
  journal={Discrete \& Computational Geometry},
  volume={70},
  number={4},
  pages={1378--1450},
  year={2023},
  publisher={Springer}
}

@article{demetci2022scot,
  title={SCOT: single-cell multi-omics alignment with optimal transport},
  author={Demetci, Pinar and Santorella, Rebecca and Sandstede, Bj{\"o}rn and Noble, William Stafford and Singh, Ritambhara},
  journal={Journal of computational biology},
  volume={29},
  number={1},
  pages={3--18},
  year={2022},
  publisher={SAGE Publications Sage CA: Los Angeles, CA}
}

@inproceedings{Memoli_2007,
  booktitle = {Eurographics Symposium on Point-Based Graphics},
  editor    = {M. Botsch and R. Pajarola and B. Chen and M. Zwicker},
  title     = {{On the use of {G}romov-{H}ausdorff Distances for Shape Comparison}},
  author    = {M{\'{e}}moli, Facundo},
  pages     = {81--90},
  year      = {2007},
  publisher = {The Eurographics Association},
  issn      = {1811-7813},
  isbn      = {978-3-905673-51-7},
  doi       = {10.2312/SPBG/SPBG07/081-090}
}

@inproceedings{peyre2016gromov,
  title={Gromov-wasserstein averaging of kernel and distance matrices},
  author={Peyr{\'e}, Gabriel and Cuturi, Marco and Solomon, Justin},
  booktitle={International conference on machine learning},
  pages={2664--2672},
  year={2016},
  organization={PMLR}
}

@article{xu2019scalable,
  title={Scalable Gromov-Wasserstein learning for graph partitioning and matching},
  author={Xu, Hongteng and Luo, Dixin and Carin, Lawrence},
  journal={Advances in neural information processing systems},
  volume={32},
  year={2019}
}

@inproceedings{chowdhury2021generalized,
  title={Generalized spectral clustering via Gromov-Wasserstein learning},
  author={Chowdhury, Samir and Needham, Tom},
  booktitle={International Conference on Artificial Intelligence and Statistics},
  pages={712--720},
  year={2021},
  organization={PMLR}
}

@article{serieys2025nonlinear,
  title={Nonlinear Lebesgue spaces: Dense subspaces, completeness and separability},
  author={S{\'e}rieys, Guillaume and Trouv{\'e}, Alain},
  journal={arXiv preprint arXiv:2512.19208},
  year={2025}
}

@incollection{sturm2003probability,
  author    = {Sturm, Karl-Theodor},
  title     = {Probability measures on metric spaces of nonpositive curvature},
  booktitle = {Heat Kernels and Analysis on Manifolds, Graphs, and Metric Spaces},
  editor    = {Auscher, Pascal and Coulhon, Thierry and Grigor'yan, Alexander},
  series    = {Contemporary Mathematics},
  volume    = {338},
  pages     = {357--390},
  year      = {2003},
  publisher = {American Mathematical Society},
  doi       = {10.1090/conm/338/06080}
}

@article{korevaar1993sobolev,
  title={Sobolev spaces and harmonic maps for metric space targets},
  author={Korevaar, Nicholas J and Schoen, Richard M},
  journal={Communications in Analysis and Geometry},
  volume={1},
  number={4},
  pages={561--659},
  year={1993},
  publisher={International Press of Boston}
}

@article{jost1994equilibrium,
  title={Equilibrium maps between metric spaces},
  author={Jost, J{\"u}rgen},
  journal={Calculus of Variations and Partial Differential Equations},
  volume={2},
  number={2},
  pages={173--204},
  year={1994},
  publisher={Springer}
}

@article{needham2023geometric,
  title={Geometric averages of partitioned datasets},
  author={Needham, Tom and Weighill, Thomas},
  journal={SIAM Journal on Applied Algebra and Geometry},
  volume={7},
  number={1},
  pages={104--132},
  year={2023},
  publisher={SIAM}
}

@article{kloeckner2010geometric,
  title={A geometric study of Wasserstein spaces: Euclidean spaces},
  author={Kloeckner, Beno{\^\i}t},
  journal={Annali della Scuola Normale Superiore di Pisa-Classe di Scienze},
  volume={9},
  number={2},
  pages={297--323},
  year={2010}
}

@article{bertrand2012geometric,
  title={A geometric study of Wasserstein spaces: Hadamard spaces},
  author={Bertrand, J{\'e}r{\^o}me and Kloeckner, Beno{\^\i}t},
  journal={Journal of Topology and Analysis},
  volume={4},
  number={4},
  pages={515--542},
  year={2012},
  publisher={World Scientific}
}

@book{bacak2014convex,
  title={Convex analysis and optimization in Hadamard spaces},
  author={Ba{\v c}{\'a}k, Miroslav},
  volume={22},
  year={2014},
  publisher={Walter de Gruyter GmbH \& Co KG}
}

@incollection{plaut2001metric,
  title={Metric spaces of curvature $\geq k$},
  author={Plaut, Conrad},
  booktitle={Handbook of geometric topology},
  pages={819--898},
  year={2001},
  publisher={Elsevier}
}

@article{zhang2024geometry,
  title={Geometry of the Space of Partitioned Networks: A Unified Theoretical and Computational Framework},
  author={Zhang, Stephen Y and Lan, Fangfei and Zhou, Youjia and Barbensi, Agnese and Stumpf, Michael PH and Wang, Bei and Needham, Tom},
  journal={Information and Inference: A Journal of the IMA},
 year={2026},
 note={To appear; arXiv:2409.06302}
}

@inproceedings{chowdhury2020gromov,
  title={Gromov-Wasserstein averaging in a Riemannian framework},
  author={Chowdhury, Samir and Needham, Tom},
  booktitle={Proceedings of the IEEE/CVF Conference on Computer Vision and Pattern Recognition Workshops},
  pages={842--843},
  year={2020}
}

@article{kapovitch2022structure,
  title={The structure of submetries},
  author={Kapovitch, Vitali and Lytchak, Alexander},
  journal={Geom. Topol},
  volume={26},
  number={6},
  pages={2649--2711},
  year={2022}
}

@book{alexander2024alexandrov,
  title={Alexandrov geometry: foundations},
  author={Alexander, Stephanie and Kapovitch, Vitali and Petrunin, Anton},
  volume={236},
  year={2024},
  publisher={American Mathematical Society}
}

@article{gomez2025metrics,
  title={Metrics for Parametric Families of Networks},
  author={G{\'o}mez, Mario and Ma, Guanqun and Needham, Tom and Wang, Bei},
  journal={arXiv preprint arXiv:2509.22549},
  year={2025}
}

@article{Memoli2022distance,
  title={Distance distributions and inverse problems for metric measure spaces},
  author={M{\'e}moli, Facundo and Needham, Tom},
  journal={Studies in Applied Mathematics},
  volume={149},
  number={4},
  pages={943--1001},
  year={2022},
  publisher={Wiley Online Library}
}

@article{pottmann2009integral,
  title={Integral invariants for robust geometry processing},
  author={Pottmann, Helmut and Wallner, Johannes and Huang, Qi-Xing and Yang, Yong-Liang},
  journal={Computer Aided Geometric Design},
  volume={26},
  number={1},
  pages={37--60},
  year={2009},
  publisher={Elsevier}
}

@article{manay2006integral,
  title={Integral invariants for shape matching},
  author={Manay, Siddharth and Cremers, Daniel and Hong, Byung-Woo and Yezzi, Anthony J and Soatto, Stefano},
  journal={IEEE Transactions on pattern analysis and machine intelligence},
  volume={28},
  number={10},
  pages={1602--1618},
  year={2006},
  publisher={IEEE}
}

@article{belongie2002shape,
  title={Shape matching and object recognition using shape contexts},
  author={Belongie, Serge and Malik, Jitendra and Puzicha, Jan},
  journal={IEEE transactions on pattern analysis and machine intelligence},
  volume={24},
  number={4},
  pages={509--522},
  year={2002},
  publisher={IEEE}
}

@MISC {CompletenessMathOverflow,
    TITLE = {Does complete and separable Wasserstein space imply a complete base space?},
    AUTHOR = {Iosif Pinelis},
    HOWPUBLISHED = {MathOverflow},
    YEAR = {2024},
    NOTE = {URL:https://mathoverflow.net/q/470226 (version: 2024-10-23)},
    EPRINT = {https://mathoverflow.net/q/470226},
    URL = {https://mathoverflow.net/q/470226}
}

@book{villani2021topics,
  title={Topics in optimal transportation},
  author={Villani, C{\'e}dric},
  volume={58},
  year={2003},
  publisher={American Mathematical Soc.}
}

@book{villani2008optimal,
  title={Optimal transport: old and new},
  author={Villani, C{\'e}dric},
  volume={338},
  year={2009},
  publisher={Springer}
}

@article{berestovskii2000metric,
  title={A metric characterization of Riemannian submersions},
  author={Berestovskii, Valerii N and Guijarro, Luis},
  journal={Annals of Global Analysis and Geometry},
  volume={18},
  number={6},
  pages={577--588},
  year={2000},
  publisher={Springer}
}

@article{berestovskii1987submetries,
  title={Submetries of space-forms of negative curvature},
  author={Berestovskii, VN},
  journal={Siberian Mathematical Journal},
  volume={28},
  number={4},
  pages={552--562},
  year={1987},
  publisher={Springer}
}

@book{lovasz2012large,
  title={Large networks and graph limits},
  author={Lov{\'a}sz, L{\'a}szl{\'o}},
  volume={60},
  year={2012},
  publisher={American Mathematical Soc.}
}

@article{alvarado2023limits,
  title={Limits of multi-relational graphs},
  author={Alvarado, Juan and Wang, Yuyi and Ramon, Jan},
  journal={Machine Learning},
  volume={112},
  number={1},
  pages={177--216},
  year={2023},
  publisher={Springer}
}

@article{kunszenti2022multigraph,
  title={Multigraph limits, unbounded kernels, and Banach space decorated graphs},
  author={Kunszenti-Kov{\'a}cs, D{\'a}vid and Lov{\'a}sz, L{\'a}szl{\'o} and Szegedy, Bal{\'a}zs},
  journal={Journal of Functional Analysis},
  volume={282},
  number={2},
  pages={109284},
  year={2022},
  publisher={Elsevier}
}

@article{brenier2003approximation,
  title={Approximation of maps by diffeomorphisms},
  author={Brenier, Yann and Gangbo, Wilfrid},
  journal={Calculus of Variations and Partial Differential Equations},
  volume={16},
  number={2},
  pages={147--164},
  year={2003},
  publisher={Springer}
}

@article{abraham2023probability,
     author = {Abraham, Romain and Delmas, Jean-Fran\c{c}ois and Weibel, Julien},
     title = {Probability-graphons: {Limits} of large dense weighted graphs},
     journal = {Innovations in Graph Theory},
     pages = {25--117},
     year = {2025},
     publisher = {Stichting Innovations in Graph Theory},
     volume = {2},
     doi = {10.5802/igt.7},
     language = {en},
     url = {https://igt.centre-mersenne.org/articles/10.5802/igt.7/}
}

@article{lovasz2010limits,
  title   = {Limits of compact decorated graphs},
  author  = {Lov{\'a}sz, L{\'a}szl{\'o} and Szegedy, Bal{\'a}zs},
  journal = {arXiv preprint arXiv:1010.5155},
  year    = {2010}
}

@book{sturm2020,
  title={The space of spaces: curvature bounds and gradient flows on the space of metric measure spaces},
  author={Sturm, Karl-Theodor},
  volume={290},
  number={1443},
  year={2023},
  publisher={American Mathematical Society}
}

@article{Memoli_needham_2024_gw_gm,
  author  = {Facundo M\'{e}moli and Tom Needham},
  title   = {Comparison results for {G}romov--{W}asserstein and {G}romov--{M}onge distances},
  journal = {ESAIM: Control, Optimisation and Calculus of Variations},
  year    = {2024},
  volume  = {30},
  pages   = {78},
  doi     = {10.1051/cocv/2024063},
  url     = {https://www.esaim-cocv.org/articles/cocv/abs/2024/01/cocv230154/cocv230154.html},
  month   = oct,
  note    = {Article 78, 21 pp.; published online 07 Oct 2024, Open Access}
}

@article{bauer2025zgromovwassersteindistance,
  title={The $Z$-{G}romov-{W}asserstein distance},
  author={Bauer, Martin and M\'{e}moli, Facundo and Needham, Tom and Nishino, Mao},
  journal={Journal of Machine Learning Research},
  volume={26},
  number={291},
  pages={1--57},
  year={2025}
}

@article{Sturm2006geometryI,
  author  = {Sturm, Karl-Theodor},
  title   = {On the geometry of metric measure spaces. {I}},
  journal = {Acta Mathematica},
  volume  = {196},
  number  = {1},
  pages   = {65--131},
  year    = {2006},
  month   = jul,
  doi     = {10.1007/s11511-006-0002-8}
}

@article{Memoli2011GromovWasserstein,
  author  = {M{\'e}moli, Facundo},
  title   = {Gromov--Wasserstein Distances and the Metric Approach to Object Matching},
  journal = {Foundations of Computational Mathematics},
  volume  = {11},
  number  = {4},
  pages   = {417--487},
  year    = {2011},
  doi     = {10.1007/s10208-011-9093-5},
  url     = {https://dx.doi.org/10.1007/s10208-011-9093-5}
}

@article{Vayer2020FGW,
  author    = {Titouan Vayer and Laetitia Chapel and R{\'e}mi Flamary and Romain Tavenard and Nicolas Courty},
  title     = {Fused Gromov--Wasserstein Distance for Structured Objects},
  journal   = {Algorithms},
  year      = {2020},
  volume    = {13},
  number    = {9},
  articleno = {212},
  issn      = {1999-4893},
  doi       = {10.3390/a13090212}
}

@article{chowdhury2019gromovwassersteindistancenetworksstable,
  title={The Gromov--Wasserstein distance between networks and stable network invariants},
  author={Chowdhury, Samir and M{\'e}moli, Facundo},
  journal={Information \& Inference: A Journal of the IMA},
  volume={8},
  number={4},
  pages={757--787},
  year={2019},
  doi={10.1093/imaiai/iaz026}
}

@book{BBI2001,
  author    = {Burago, Dmitri and Burago, Yuri and Ivanov, Sergei},
  title     = {A Course in Metric Geometry},
  series    = {Graduate Studies in Mathematics},
  volume    = {33},
  publisher = {American Mathematical Society},
  address   = {Providence, RI},
  year      = {2001},
  isbn      = {978-0-8218-2129-9},
  doi       = {10.1090/gsm/033},
  mrnumber  = {1835418}
}

@book{kallenberg2002foundations,
  title     = {Foundations of Modern Probability},
  author    = {Kallenberg, Olav},
  year      = {2002},
  edition   = {2nd},
  publisher = {Springer},
  series    = {Probability and Its Applications},
  address   = {New York},
  isbn      = {978-0-387-95313-7},
  doi       = {10.1007/978-1-4757-4015-8}
}

@book{santambrogio2015optimal,
  title     = {Optimal Transport for Applied Mathematicians: Calculus of Variations, PDEs, and Modeling},
  author    = {Santambrogio, Filippo},
  series    = {Progress in Nonlinear Differential Equations and Their Applications},
  volume    = {87},
  year      = {2015},
  publisher = {Birkh{\"a}user, Cham},
  doi       = {10.1007/978-3-319-20828-2},
  isbn      = {978-3-319-20827-5},
  url       = {https://doi.org/10.1007/978-3-319-20828-2}
}

@article{Chowdhury2023HypergraphCOT,
  author  = {Samir Chowdhury and Tom Needham and Ethan Semrad and Bei Wang and Youjia Zhou},
  title   = {Hypergraph Co-Optimal Transport: Metric and Categorical Properties},
  journal = {Journal of Applied and Computational Topology},
  volume  = {8},
  pages   = {1171--1230},
  year    = {2024},
  doi     = {10.1007/s41468-023-00142-9},
  url     = {https://doi.org/10.1007/s41468-023-00142-9}
}

@book{Janson2013Graphons,
  author    = {Svante Janson},
  title     = {Graphons, Cut Norm and Distance, Couplings and Rearrangements},
  series    = {New York Journal of Mathematics Monographs},
  volume    = {4},
  publisher = {State University of New York, University at Albany},
  address   = {Albany, NY},
  year      = {2013}
}

@book{bridsonhaefliger1999,
  author    = {Bridson, Martin R. and Haefliger, Andr{\'e}},
  title     = {Metric Spaces of Non-Positive Curvature},
  series    = {Grundlehren der Mathematischen Wissenschaften},
  volume    = {319},
  publisher = {Springer},
  address   = {Berlin},
  year      = {1999},
  isbn      = {978-3-540-64324-1}
}

@article{BrediesFanzon2020,
  author    = {Bredies, Kristian and Fanzon, Silvio},
  title     = {An optimal transport approach for solving dynamic inverse problems in spaces of measures},
  journal   = {ESAIM: Mathematical Modelling and Numerical Analysis},
  volume    = {54},
  number    = {6},
  pages     = {2351--2382},
  year      = {2020},
  publisher = {EDP-Sciences},
  doi       = {10.1051/m2an/2020056},
  url       = {https://www.numdam.org/articles/10.1051/m2an/2020056/}
}

\end{document}